\documentclass[11pt]{article}
\usepackage[utf8]{inputenc}
\usepackage{amsmath,amsthm,amssymb,amsfonts,amscd}
\usepackage{latexsym,color,graphicx,mathrsfs,cite}
\usepackage[title]{appendix}
\usepackage{enumitem}
\usepackage[colorlinks=true,urlcolor=black,
citecolor=blue,linkcolor=blue,linktocpage,pdfpagelabels,
bookmarksnumbered,bookmarksopen]{hyperref}

\numberwithin{equation}{section}

\newtheorem{theorem}{Theorem}[section]
\newtheorem{corollary}[theorem]{Corollary}
\newtheorem{lemma}[theorem]{Lemma}
\newtheorem{proposition}[theorem]{Proposition}
\newtheorem{definition}[theorem]{Definition}

\newtheorem*{remark*}{Remark}

\DeclareMathOperator{\supp}{supp}
\DeclareMathOperator{\ind}{ind}
\DeclareMathOperator{\spec}{spec}

\DeclareMathOperator{\dist}{dist}

\allowdisplaybreaks

\hypersetup{
  pdftitle={Multi-bubble solutions for an exponential biharmonic Neumann problem in dimension four},
  pdfauthor={Jiaoping Chen, Wenjing Chen and Shengbing Deng},
  pdfsubject={Boundary concentration for an exponential biharmonic Neumann problem},
  pdfkeywords={Biharmonic equation, boundary concentration, Lyapunov--Schmidt reduction, local uniqueness, nondegeneracy, Morse index}
}

\begin{document}
\title{Multi-bubble solutions for an exponential biharmonic Neumann problem in dimension four}
\author{Jiaoping Chen\thanks{E-mail addresses:\, \texttt{cjp123888@email.swu.edu.cn} (J. Chen), \texttt{wjchen@swu.edu.cn} (W. Chen), \texttt{shbdeng@swu.edu.cn} (S. Deng).},
\ Wenjing Chen\ and Shengbing Deng\thanks{Corresponding author: Shengbing Deng.}\\
\footnotesize School of Mathematics and Statistics, Southwest University,\\
\footnotesize Chongqing, 400715, P.R. China}
\date{}
\maketitle

\begin{abstract}
We are concerned with the biharmonic Neumann problem
\[
\Delta^2u+u=0\quad\text{in }\Omega,\qquad
\partial_\nu u=0,\quad
\partial_\nu\Delta u=-\varepsilon h e^{24\pi^2u}
\quad\text{on }\partial\Omega,
\]
where $\Omega\subset\mathbb R^4$ is a smooth bounded domain, $\nu$ is the outward unit normal, $h\in C^5(\partial\Omega)$ is positive, and $\varepsilon>0$ is sufficiently small.
For each fixed integer $m\ge1$, we construct $m$-bubble solutions $u_\varepsilon$ by
Lyapunov--Schmidt reduction. As $\varepsilon\to0$, the solutions
concentrate at $m$ distinct boundary points, each carrying unit mass
for the measure $\varepsilon h e^{24\pi^2u_\varepsilon}\,d\sigma$.
For solutions corresponding to nondegenerate critical points $\xi^*$ of
the interaction function $\mathcal F_m$, we establish local uniqueness
and  nondegeneracy, and prove $\ind(u_\varepsilon)=m+\ind(-D^2\mathcal F_m(\xi^*))$.
In particular, for a generic positive weight $h$, all critical points
of $\mathcal F_m$ are nondegenerate, and the $m$-bubble solutions
associated with global minima of $\mathcal F_m$ have Morse index $4m$.

\vspace{.2cm}
\emph{\bfseries Keywords:} Biharmonic equation,  Lyapunov--Schmidt reduction, local uniqueness, nondegeneracy, Morse index

\vspace{.2cm}
\emph{\bfseries 2020 Mathematics Subject Classification:} 35J30, 35B33, 35B40,  35P30
\end{abstract}

\section{Introduction}
Let \(\Omega\subset\mathbb R^4\) be a smooth bounded domain, and let \(\nu\) denote the outward unit normal to \(\partial\Omega\). For \(h\in C^5(\partial\Omega)\) with \(h>0\) and \(\varepsilon>0\) sufficiently small, we consider
\begin{equation}\label{e1.1}
\begin{cases}
\Delta^2u+u=0 & \text{in }\Omega,\\
\partial_\nu u=0 & \text{on }\partial\Omega,\\
\partial_\nu\Delta u=-\varepsilon h(x)e^{24\pi^2u} & \text{on }\partial\Omega.
\end{cases}
\end{equation}
Problem \eqref{e1.1} is  the Euler-Lagrange equation for the functional
\begin{equation}\label{e1.2}
J_\varepsilon(u)=\frac12\int_\Omega\bigl((\Delta u)^2+u^2\bigr)\,dx
-\frac{\varepsilon}{24\pi^2}\int_{\partial\Omega}he^{24\pi^2u}\,d\sigma,
\end{equation}
defined on
\[
\mathcal H=\{u\in H^2(\Omega):\partial_\nu u=0\text{ in }H^{\frac{1}{2}}(\partial\Omega)\}.
\]
{
The trace map and Sobolev embeddings are continuous:
}
\[
H^2(\Omega)\xrightarrow{\operatorname{Tr}}
H^{\frac32}(\partial\Omega)
\hookrightarrow L^p(\partial\Omega),
\qquad 1\le p<\infty,
\]
whereas $H^{3/2}(\partial\Omega)\not\hookrightarrow L^\infty(\partial\Omega)$.
{
The sharp trace Adams inequality of Chen et al.~\cite{r11} states that
}
\[
\sup_{\substack{u\in \mathcal{H},\ \|\Delta u\|_{L^2(\Omega)}\le1\\
\int_\Omega u\,dx=0}}
\int_{\partial\Omega}e^{12\pi^2u^2}\,d\sigma<\infty.
\]
{
Together with Young's inequality, this ensures that $J_\varepsilon(u)$ is finite for every $u\in\mathcal H$.
}
For related results on Moser–Trudinger, Adams and higher-order trace inequalities, we refer the reader to \cite{r28,r2,r6,r1,r32} and the references therein.

On planar domains, D\'avila et al.~\cite{r13} constructed solutions with any prescribed number of boundary peaks for
\begin{equation*}
\begin{cases}
-\Delta u + u = 0 & \text{in } \Omega,\\[4pt]
\dfrac{\partial u}{\partial \nu} = \varepsilon e^{u} & \text{on } \partial\Omega.
\end{cases}
\end{equation*}
De Marchis et al. \cite{r34} studied the related boundary
condition $\partial_\nu u=|u|^{p-1}u$ and established energy quantization and sharp boundary concentration as
$p\to\infty$ for families of positive solutions satisfying a uniform energy bound.
Qualitative properties of solutions have also been investigated for the planar Dirichlet problems
\[
-\Delta u=\lambda u e^{u^2}\quad(\lambda\to0),
\qquad
-\Delta u=u^p\quad(p\to\infty).
\]
In the Moser–Trudinger setting, Luo et al. \cite{r35} obtained results on the Morse index, nondegeneracy, and uniqueness of low-energy positive solutions, while Chen et al. \cite{r31} established uniqueness of positive solutions on the unit disc. The Lane–Emden equation was studied by Ianni et al. \cite{r40}, who derived Morse index and degree formulas for multi-spike solutions, together with the local uniqueness at nondegenerate critical points.

{
In dimension four, Baraket et al.~\cite{r5} and Clapp et al.~\cite{r12}
constructed interior multi-bubble solutions for
\[
\begin{cases}
\Delta^2u=\rho^4k(x)e^u & \text{in }\Omega,\\
u=\Delta u=0 & \text{on }\partial\Omega.
\end{cases}
\]
The constructions use nondegenerate critical points when $k\equiv1$
and a topological condition for positive $k$, respectively.
Ohtsuka and Takahashi~\cite{OhtsukaTakahashi2015} studied the corresponding
linearized equation and proved local asymptotic nondegeneracy.
For an exponential biharmonic Choquard equation with Navier boundary
conditions, Chen and Deng~\cite{r39} established local uniqueness,
nondegeneracy and a Morse index formula.
For related results on boundary concentration, ground states and
a priori estimates, see~\cite{r4,r10,r41}.

Boundary bubbling also occurs in the resonant prescribed $T$-curvature
problem studied by Ndiaye~\cite{r36} through critical points at infinity
and Morse theory. In \eqref{e1.1}, the term $+u$ makes the boundary
operator strictly positive. We construct boundary multi-bubble
solutions and determine their local uniqueness, linear nondegeneracy
and Morse index.
}
To state the results, we introduce the boundary operator and the interaction function.
For $f\in H^{3/2}(\partial\Omega)$, let $\mathcal Ef\in H^2(\Omega)$ solve
\begin{equation}\label{e1.3}
\Delta^2\mathcal Ef+\mathcal Ef=0\quad\text{in }\Omega,\qquad
\mathcal Ef=f,\quad\partial_\nu\mathcal Ef=0\quad\text{on }\partial\Omega,
\end{equation}
and set
\[
\mathcal Pf=-\partial_\nu\Delta\mathcal Ef\big|_{\partial\Omega}
\in H^{-\frac{3}{2}}(\partial\Omega).
\]
The operator $\mathcal P$, with domain $H^3(\partial\Omega)$, is positive and self-adjoint on $L^2(\partial\Omega,d\sigma)$. Its principal symbol is $2|\zeta|_{g(x)}^3$ for $\zeta\in T_x^*(\partial\Omega)\setminus\{0\}$ (see Lemmas~\ref{t2.1}--\ref{t2.2}), where $g$ is the induced metric on $\partial\Omega$ and \(T_x^*(\partial\Omega)\) denotes the cotangent space at \(x\).
Thus \eqref{e1.1} is equivalent to
\begin{equation}\label{e1.4}
\mathcal Pv=\varepsilon he^{24\pi^2v}\quad\text{on }\partial\Omega,
\qquad v=u|_{\partial\Omega},\quad u=\mathcal Ev.
\end{equation}

The boundary concentration of solutions to \eqref{e1.1}  is governed by the following limiting problem in the upper half-space:
\[
\Delta^2U=0\quad\text{in }\mathbb R^4_+,\qquad
\partial_tU=0,\quad\partial_t\Delta U=e^{24\pi^2U}
\quad\text{on }\partial\mathbb R^4_+,
\]
where $\mathbb R^4_+=\{(y,t):y\in\mathbb R^3,\ t>0\}$. 
In particular, a normalized bubble solution is
\[
U(y,t)=-\frac1{8\pi^2}\log\bigl((1+at)^2+a^2|y|^2\bigr)
+\frac{at}{4\pi^2\bigl((1+at)^2+a^2|y|^2\bigr)},
\qquad a=\left(\frac\pi2\right)^{\frac{2}{3}},
\]
whose boundary trace and density are
\begin{equation}\label{e1.5}
u_0(y)=-\frac1{8\pi^2}\log(1+a^2|y|^2),\qquad
q(y)=e^{24\pi^2u_0(y)}=(1+a^2|y|^2)^{-3}.
\end{equation}
Related results on half-space classification and conformal
boundary problems can be found in \cite{r11,r33,r37,r36}.
The model boundary operator is $2(-\Delta)^{3/2}$,
where $(-\Delta)^{3/2}$ has Fourier symbol $|\zeta|^3$.
For any $f\in L^1_{\mathrm{loc}}(\mathbb R^3)$ satisfying
$f(y)=O(\log(2+|y|))$ as $|y|\to\infty$,
the distribution $(-\Delta)^{3/2}f$ is defined by
\[
\langle (-\Delta)^{\frac{3}{2}}f,\psi\rangle
=\int_{\mathbb R^3}
f(y)(-\Delta)^{\frac{3}{2}}\psi(y)\,dy,
\qquad \psi\in\mathcal S(\mathbb R^3),
\]
where $\mathcal S$ denotes the Schwartz space and $\mathcal S'$ its continuous dual.
{
By~\cite[Theorem~1.1]{r9}, the bounded distributional kernel of $\mathcal L_0=2(-\Delta)^{3/2}-24\pi^2q$ is
}
\begin{equation}\label{e1.6}
\operatorname{span}\{z_0,z_1,z_2,z_3\},\qquad
z_0(y)=\frac{1-a^2|y|^2}{1+a^2|y|^2},\quad
z_k(y)=\frac{2ay_k}{1+a^2|y|^2}\quad(1\le k\le3).
\end{equation}
Here $z_0$ corresponds to dilation, while $z_1,z_2,z_3$ correspond to translations.

Let $G_{\partial\Omega}(x,\xi)$ denote the Green kernel of {$\mathcal P$}, characterized by
$
{\mathcal P_x}G_{\partial\Omega}(x,\xi)=\delta_\xi$, where $\delta_\xi$ is the unit Dirac mass.
Fix a symmetric function $d\ge0$ with $d^2\in C^\infty(\partial\Omega\times\partial\Omega)$, positive off the diagonal and equal to the geodesic distance near the diagonal.
The Green kernel is symmetric, and the Robin function is
\begin{equation}\label{e1.7}
r(\xi)=\lim_{x\to\xi}\left[
G_{\partial\Omega}(x,\xi)+\frac1{4\pi^2}\log d(x,\xi)
\right].
\end{equation}
For $m\ge1$, define
\[
\mathcal M_m=\left\{(\xi_1,\ldots,\xi_m)\in(\partial\Omega)^m:
\xi_i\ne\xi_j\text{ for }i\ne j\right\},
\]
\[
\mathcal F_m(\xi)=\frac12\sum_{i=1}^m r(\xi_i)
+\sum_{1\le i<j\le m}G_{\partial\Omega}(\xi_i,\xi_j)
+\frac1{24\pi^2}\sum_{i=1}^m\log h(\xi_i).
\]
Then
\begin{equation}\label{e1.8}
\nabla_{\xi_i}\mathcal F_m
=\frac12\nabla r(\xi_i)+\sum_{j\ne i}\nabla_xG_{\partial\Omega}(\xi_i,\xi_j)
+\frac1{24\pi^2}\nabla\log h(\xi_i).
\end{equation}
We recall the definition of $C^1$-stable critical set that we will use in the sequel.
\begin{definition}
A compact set $K\subset\mathcal M_m$ is a {$C^1$-stable critical set} of $\mathcal F_m$ if some open set $D\subset\subset\mathcal M_m$ satisfies
\begin{equation}\label{e1.9}
K=\{\xi\in D:\nabla\mathcal F_m(\xi)=0\},\qquad
\nabla\mathcal F_m\ne0\quad\text{on }\partial D,
\end{equation}
and every $\widetilde F$ sufficiently close to $\mathcal F_m$ in $C^1(\overline D)$ has a critical point in $D$.
\end{definition}
In particular, isolated nondegenerate critical points are $C^1$-stable.

For each configuration $\xi\in\mathcal M_m$, let $V_\xi$ denote the approximate $m$-bubble boundary profile constructed in Section~\ref{s3}, with concentration scales $\delta_i(\varepsilon,\xi)$.
{
Our first result concerns the existence of multi-bubble solutions.
}

\begin{theorem}\label{t1.1}
Let $m\ge1$ and let $K\subset\mathcal M_m$ be a $C^1$-stable critical set of $\mathcal F_m$. For  sufficiently small $\varepsilon>0$, problem \eqref{e1.1} has an $m$-bubble solution of the form
\[
u_\varepsilon
=\sum_{i=1}^m\mathcal E W_{\delta_{i,\varepsilon},\xi_{i,\varepsilon}}
+\mathcal E\phi_{\varepsilon,\xi_\varepsilon}
\quad\text{in }\Omega,
\]
with $\xi_\varepsilon=(\xi_{1,\varepsilon},\ldots,\xi_{m,\varepsilon})
\in\mathcal M_m$,  $\delta_{i,\varepsilon}=\delta_i(\varepsilon,\xi_\varepsilon)$,
and a boundary correction $\phi_{\varepsilon,\xi_\varepsilon}\in H^{3/2}(\partial\Omega)\cap L^\infty(\partial\Omega)$ satisfying
\[
\dist(\xi_\varepsilon,K)\rightarrow0,\qquad
\max_{1\le i\le m}\delta_{i,\varepsilon}\rightarrow0,\qquad
\|\phi_{\varepsilon,\xi_\varepsilon}\|_{L^\infty(\partial\Omega)}
+\|\phi_{\varepsilon,\xi_\varepsilon}\|_{H^{\frac{3}{2}}(\partial\Omega)}
\rightarrow0,
\]
where $\delta_i(\varepsilon,\xi_\varepsilon)$ and $W_{\delta_{i,\varepsilon},\xi_{i,\varepsilon}}$ are defined in \eqref{e3.10} and \eqref{defW}.
{
Along any subsequence with $\xi_\varepsilon\to\xi^*=(\xi_1^*,\ldots,\xi_m^*)\in K$,
\begin{equation}\label{e1.mass-convergence}
\varepsilon h e^{24\pi^2u_\varepsilon}\,d\sigma
\rightharpoonup\sum_{i=1}^m\delta_{\xi_i^*}
\end{equation}
as finite measures on $\partial\Omega$.
In particular, every nondegenerate critical point $\xi^*$ of
$\mathcal F_m$ yields such a family with $\xi_\varepsilon\to\xi^*$.
}
\end{theorem}

For $m\ge2$, $\mathcal F_m(\xi)\to+\infty$ as $\min_{i\ne j}d(\xi_i,\xi_j)\to0,$
while $\mathcal M_1=\partial\Omega$ is compact.
Hence $\mathcal F_m$ attains a global minimum for every $m\ge1$.
A maximization argument for the reduced energy then yields the following corollary.
\begin{corollary}\label{t1.2}
For every fixed $m\ge1$, problem \eqref{e1.1} has an $m$-bubble solution of the form given in Theorem~\ref{t1.1} for  sufficiently small $\varepsilon>0$, with
$
\mathcal F_m(\xi_\varepsilon)\rightarrow\min_{\mathcal M_m}\mathcal F_m.
$
\end{corollary}

We next address the local uniqueness and nondegeneracy. 

\begin{theorem}\label{t1.3}
Let $\xi^*$ be a nondegenerate critical point of $\mathcal F_m$, and let $u_\varepsilon$ be the family in Theorem~\ref{t1.1} with $\xi_\varepsilon\to\xi^*$.
For  sufficiently small $\varepsilon>0$, $u_\varepsilon$ is nondegenerate: the linearized problem
\begin{equation}\label{e1.linearized-problem}
\begin{cases}
\Delta^2\Phi+\Phi=0 & \text{in }\Omega,\\
\partial_\nu\Phi=0 & \text{on }\partial\Omega,\\
\partial_\nu\Delta\Phi=-24\pi^2\varepsilon h e^{24\pi^2u_\varepsilon}\Phi
& \text{on }\partial\Omega
\end{cases}
\end{equation}
has only the zero weak solution $\Phi\in H^2(\Omega)$.
{
Let $\Omega_{\star}\subset\subset\mathcal M_m$ be a sufficiently small
neighborhood of $\xi^*$, and fix $A>0$.
If $\widetilde u_\varepsilon$ solves \eqref{e1.1} and satisfies
}
\begin{equation}\label{e1.uniqueness-class}
\bigl\|\widetilde u_\varepsilon|_{\partial\Omega}-V_\xi
\bigr\|_{L^\infty(\partial\Omega)}
\le A\varepsilon^{\frac{1}{3}}|\log\varepsilon|^2
\qquad\text{for some }\xi\in\Omega_{\star}{,}
\end{equation}
{
then $\widetilde u_\varepsilon=u_\varepsilon$ in $\Omega$ for
$\varepsilon>0$ sufficiently small.
}
\end{theorem}

Let $\ind(u_\varepsilon)$ denote the Morse index of $u_\varepsilon$. For a real symmetric matrix $A$, $\ind(A)$ denotes the number of its negative eigenvalues, counted with multiplicity.

\begin{theorem}\label{t1.4}
Under the assumptions of Theorem~\ref{t1.3},
 for sufficiently small $\varepsilon>0$,
\[
\ind(u_\varepsilon)
=m+\ind\bigl(-D^2\mathcal F_m(\xi^*)\bigr)
=4m-\ind\bigl(D^2\mathcal F_m(\xi^*)\bigr).
\]
\end{theorem}

If $\Omega$ is a ball, then $r(\xi)$ is constant. Thus, for $m=1$,
\[
\ind(u_\varepsilon)=
\begin{cases}
1,&\xi^*\text{ is a nondegenerate local maximum of }h,\\
4,&\xi^*\text{ is a nondegenerate local minimum of }h.
\end{cases}
\]
Finally, nondegeneracy of the critical points of $\mathcal F_m$ is generic with respect to $h$.
Set
\[
\mathcal W_+^5(\partial\Omega)
=\{h\in C^5(\partial\Omega):\min_{\partial\Omega}h>0\},
\]
with the $C^5$ topology. {
A subset of $\mathcal W_+^5(\partial\Omega)$ is residual if it contains a countable intersection of open dense subsets.
}

\begin{corollary}\label{t1.5}
There exists a residual subset
$\mathcal W_{\mathrm{reg}}\subset\mathcal W_+^5(\partial\Omega)$
such that, for every $h\in\mathcal W_{\mathrm{reg}}$ and  $m\ge1$,
all critical points of $\mathcal F_m$ in $\mathcal M_m$ are nondegenerate,
and the $m$-bubble
solutions $u_\varepsilon$ associated with global minima of
$\mathcal F_m$ satisfy
$
\operatorname{ind}(u_\varepsilon) = 4m
$
 for
sufficiently small $\varepsilon>0$.
\end{corollary}

{
The reduction requires precise control of the curvature contribution
to the residual and the nondecaying dilation mode.
The leading curvature correction is even, so its
translation moments vanish, also after parameter differentiation;
the reduced equations are therefore governed by $d\mathcal F_m$
(see Proposition~\ref{t2.3} and Lemma~\ref{t5.3}).
A correction by global Green potentials controls the dilation mode
and yields the projected inverse estimate
(see Lemma~\ref{t4.6} and Proposition~\ref{t4.7}).
Estimates in moving coordinates yield the $C^2$ expansion needed for
local uniqueness and nondegeneracy (see Proposition~\ref{t5.7}).
The spectral comparison separates the $m$ negative directions from
the translation contribution $\ind(-D^2\mathcal F_m(\xi^*))$;
the corrected dilation directions are positive
(see Proposition~\ref{t7.5} and Corollary~\ref{t7.6}).
}

The paper is organized as follows. Section \ref{s2} derives the Green kernel expansion, analyzes the model bubble and its linearization, and constructs approximate solutions. Section \ref{s4} carries out the finite-dimensional reduction and establishes the \(C^2\) expansion of the reduced energy. Section \ref{s6} proves Theorem \ref{t1.1} and Corollary \ref{t1.2}.
Finally, Section \ref{s7} establishes the local uniqueness, nondegeneracy, and Morse index formula in Theorems \ref{t1.3}--\ref{t1.4}, as well as the generic nondegeneracy result in Corollary \ref{t1.5}.

{
Throughout, $m$ is fixed, and $C>0$ may vary from line to line but is
independent of $\varepsilon$ and uniform for $\xi$ in each fixed compact
subset of $\mathcal M_m$. We write $A\asymp B$ if
$C^{-1}B\le A\le CB$, and $I$ for the identity on the relevant space.
The notation $O_{C^k}(\cdot)$ and $o_{C^k}(\cdot)$ includes configuration
derivatives up to order $k$. Gradients and Hessians are taken with
respect to the induced metric and its product metric; parameter
derivatives are coordinate derivatives unless a material representative
is specified (see Subsection~\ref{s5.1}).
}
For simplicity, we write
\[
\|\cdot\|_{L^p}=\|\cdot\|_{L^p(\partial\Omega,d\sigma)},\quad
\|\cdot\|_{H^s}=\|\cdot\|_{H^s(\partial\Omega)},\quad
\|\cdot\|_\infty=\|\cdot\|_{L^\infty(\partial\Omega)}.
\]

\section{Preliminaries and approximate solutions}\label{s2}

\subsection{The boundary operator and its Green kernel}\label{s2.1}

After lifting the boundary data in \eqref{e1.3}, the Lax--Milgram theorem
applied to the coercive form $
(u,v)\mapsto\int_\Omega(\Delta u\Delta v+uv)\,dx$
on \(\{u\in H^2(\Omega):u=\partial_\nu u=0\text{ on }\partial\Omega\}\)
provides a unique weak solution \(\mathcal Ef\in H^2(\Omega)\).
Elliptic regularity~\cite{r25,r15} yields
\[
\|\mathcal Ef\|_{H^{s+2}(\Omega)}
\le C_s\|f\|_{H^{s+\frac{3}{2}}(\partial\Omega)},\qquad s\ge0.
\]

\begin{lemma}\label{t2.1}
For every $f,g\in H^{3/2}(\partial\Omega)$,
\begin{equation}\label{e2.1}
\langle\mathcal Pf,g\rangle
=\int_\Omega\bigl(\Delta\mathcal Ef\,\Delta\mathcal Eg
+\mathcal Ef\,\mathcal Eg\bigr)\,dx,
\end{equation}
where $\langle\cdot,\cdot\rangle$ denotes the duality pairing between $H^{-3/2}(\partial\Omega)$ and $H^{3/2}(\partial\Omega)$.
There exist constants $c,C>0$, depending only on $\Omega$, such that
\begin{equation}\label{e2.coercivity}
c\|f\|_{H^{\frac{3}{2}}(\partial\Omega)}^2
\le\langle\mathcal Pf,f\rangle
\le C\|f\|_{H^{\frac{3}{2}}(\partial\Omega)}^2
\qquad\text{for all }f\in H^{\frac{3}{2}}(\partial\Omega).
\end{equation}
\end{lemma}
\begin{proof}
Let $f,g\in C^\infty(\partial\Omega)$. Multiplying \eqref{e1.3} by \(\mathcal Eg\) and integrating by parts twice, we obtain
\begin{equation*}
\langle\mathcal Pf,g\rangle
=-\int_{\partial\Omega}g\,\partial_\nu\Delta\mathcal Ef\,d\sigma=\int_\Omega\bigl(\Delta\mathcal Ef\,\Delta\mathcal Eg
+\mathcal Ef\,\mathcal Eg\bigr)\,dx.
\end{equation*}
Since $\mathcal E:H^{3/2}(\partial\Omega)\to H^2(\Omega)$ is continuous,
\eqref{e2.1} follows for all $f,g\in H^{3/2}(\partial\Omega)$ by density.
Substituting $g=f$ into \eqref{e2.1}, we have
\[
\langle\mathcal Pf,f\rangle
=\|\Delta\mathcal Ef\|_{L^2(\Omega)}^2
+\|\mathcal Ef\|_{L^2(\Omega)}^2
\le C\|f\|_{H^{\frac{3}{2}}(\partial\Omega)}^2.
\]
As $\partial_\nu\mathcal Ef=0$, the trace theorem and the
elliptic estimate for the Neumann problem imply
\begin{equation*}
\|f\|_{H^{\frac{3}{2}}(\partial\Omega)}^2
\le C\|\mathcal Ef\|_{H^2(\Omega)}^2
\le C\bigl(\|\Delta\mathcal Ef\|_{L^2(\Omega)}^2
+\|\mathcal Ef\|_{L^2(\Omega)}^2\bigr)
=C\langle\mathcal Pf,f\rangle.
\end{equation*}
\end{proof}

Lemma~\ref{t2.1} and the Lax--Milgram theorem imply that
$\mathcal P:H^{3/2}(\partial\Omega)\to H^{-3/2}(\partial\Omega)$ is an isomorphism.
To determine the leading singularity of
$G_{\partial\Omega}(x,\xi)$ as $x\to\xi$,
we first compute the principal symbol of $\mathcal P$.

\begin{lemma}\label{t2.2}
The operator $\mathcal P$ is a classical elliptic pseudodifferential operator of order three on $\partial\Omega$, with principal symbol
\[
2|\zeta|_{g(x)}^3,
\quad x\in\partial\Omega,\quad
\zeta\in T_x^*(\partial\Omega)\setminus\{0\}.
\]
\end{lemma}
\begin{proof}
Since \eqref{e1.3} is an elliptic Dirichlet problem, $\mathcal P$ is a classical pseudodifferential operator of order at most three~\cite{r17,r24}.
Fix $x\in\partial\Omega$ and choose boundary normal coordinates $(y,t)$
centered at $x$, with $t>0$ in $\Omega$ and the tangential coordinates
orthonormal at $x$.
The extension operator for the principal part of \eqref{e1.3},
with coefficients frozen at $x$, has tangential Fourier multiplier
$(1+|\zeta|t)e^{-|\zeta|t}$~\cite[Lemma~2.1]{r9},
where $\zeta\in\mathbb R^3$ is the Fourier variable dual to $y$.
Since $\partial_\nu=-\partial_t$ on $t=0$, we have
\[
\left.\partial_t(\partial_t^2-|\zeta|^2)
\bigl[(1+|\zeta|t)e^{-|\zeta|t}\bigr]\right|_{t=0}
=2|\zeta|^3.
\]
Thus the principal symbol of $\mathcal P$ is $2|\zeta|_{g(x)}^3>0$
for $\zeta\ne0$, and $\mathcal P$ is elliptic of order three.
\end{proof}
Lemmas~\ref{t2.1} and~\ref{t2.2} show that \(\mathcal P\) is positive
and self-adjoint in \(L^2(\partial\Omega,d\sigma)\), with domain
\(H^3(\partial\Omega)\)~\cite{r30}.
Its inverse is a classical pseudodifferential operator of order
\(-3\)~\cite{r17,r24}, and
\[
G_{\partial\Omega}(x,\xi)=G_{\partial\Omega}(\xi,x).
\]
Fix $\xi\in\partial\Omega$. In the geodesic normal coordinates $x=\exp_\xi y$, let $g_\xi(y)$ denote the metric matrix.
Since the principal symbol in Lemma~\ref{t2.2} reduces to $2|\zeta|^3$ at $y=0$, the logarithmic singularity of $G_{\partial\Omega}$ is determined by the following fundamental solution of $2(-\Delta)^{3/2}$ (see~\cite{r21}):
\begin{equation}
\Gamma(y)=-\frac1{4\pi^2}\log|y|.\label{e2.2}
\end{equation}
Let \(\mathrm{II}_\xi(X,Y)=\langle D_X\nu,Y\rangle\) be the second
fundamental form, where \(D\) is the Euclidean connection, and let
\(H_\xi^{\partial\Omega}=\operatorname{tr}_g\mathrm{II}_\xi\).
Vector inner products are taken with respect to the induced metric.
When \(\mathrm{II}_\xi\) acts on covectors, we identify them with
tangent vectors using \(g\).
Set
\begin{equation}\label{e2.curvature-kernel}
K_{-4,\xi}(v)=\frac1{32\pi}
\left(\frac{\mathrm{II}_\xi(v,v)}{|v|}-H_\xi^{\partial\Omega}|v|\right),\qquad v\ne0,
\qquad K_{-4,\xi}(0)=0.
\end{equation}
We define $\chi\in C^\infty([0,\infty))$ with
\[
0\le\chi\le1,\qquad \chi=1\text{ on }[0,r_0],\qquad
\chi=0\text{ on }[2r_0,\infty).
\]

\begin{proposition}\label{t2.3}
The Green kernel admits a decomposition $G_{\partial\Omega}=S+H$, where $S$ and $H$ are symmetric and $H\in C^4(\partial\Omega\times\partial\Omega)$.
In normal coordinates, set $S_\xi(y,z)=S(\exp_\xi y,\exp_\xi z),$
$H_\xi(y,z)=H(\exp_\xi y,\exp_\xi z).$
Then
\begin{equation}
\begin{split}
S_\xi(y,z)&=\chi(|y-z|)\bigl[\Gamma(y-z)+K_{-4,\xi}(y-z)\bigr]+E_\xi(y,z),\\
|\partial_\xi^\alpha E_\xi(y,z)|&\le C(|y|+|z|)^2\bigl(1+|\log|y-z||\bigr),\qquad |\alpha|\le2,
\end{split}
\label{e2.6}
\end{equation}
uniformly in $\xi$, for $y,z$ sufficiently close to zero with $y\ne z$.
Here \(\alpha\) is a multi-index. The \(\xi\)-derivatives are taken
in a smooth local orthonormal frame with \(y,z\) fixed.
Moreover, the limit in \eqref{e1.7} exists, and $r$ satisfies
\begin{equation}
r(\xi)=H(\xi,\xi),\quad r\in C^4(\partial\Omega),\quad
\nabla_yH_\xi(0,0)=\nabla_zH_\xi(0,0)=\tfrac12\nabla r(\xi).
\label{e2.10}
\end{equation}
Let $\Phi_\xi$ be a family of boundary diffeomorphisms with uniformly
bounded $C^1$ norms for the maps, their inverses and their parameter
derivatives of order at most two. For $x\ne z$,
\begin{align}
|\partial_\xi^\alpha G_{\partial\Omega}(\Phi_\xi x,\Phi_\xi z)|
&\le C\bigl(1+|\log d(\Phi_\xi x,\Phi_\xi z)|\bigr),
&&|\alpha|\le2,\label{e2.transported-green}\\
|\nabla_z\partial_\xi^\alpha G_{\partial\Omega}(\Phi_\xi x,\Phi_\xi z)|
&\le \frac{C}{d(\Phi_\xi x,\Phi_\xi z)},
&&|\alpha|\le1.\label{e2.transported-green-gradient}
\end{align}
\end{proposition}
\begin{proof}
In geodesic normal coordinates, we have
\[
d\sigma(\exp_\xi y)=\sqrt{\det g_\xi(y)}\,dy,
\quad \det g_\xi(0)=1,
\quad \nabla_y\det g_\xi(0)=0.
\]
For the symbol calculations, we use the change
\[
f(\exp_\xi y)\longmapsto
(\det g_\xi(y))^{\frac{1}{4}}f(\exp_\xi y),
\]
which is an isometry from the local \(L^2(d\sigma)\) space to \(L^2(dy)\).
After this transformation, the Weyl symbols of $\mathcal P$ and $\mathcal P^{-1}$ have the asymptotic expansions $p_3+p_2+\cdots$ and $r_{-3}+r_{-4}+\cdots$, respectively. The subscripts denote the degrees of homogeneity in $\zeta$.
The relation $(\det g_\xi(y))^{1/4}=1+O(|y|^2)$ implies that, at $y=0$, the terms $p_3,p_2,r_{-3},r_{-4}$ coincide with the corresponding symbol terms before the transformation.
We extend $y$ to boundary normal coordinates $(y,t)$, with $t>0$ in $\Omega$.
Using
\begin{equation*}
g_\xi(y)=g_\xi(0)+O(|y|^2),\quad \partial_y g_\xi(0)=0,
\end{equation*}
we compute the first two symbol terms by evaluating the tangential coefficients at $y=0$ and taking the Fourier transform in $y$.
For \(\zeta\ne0\), set \(s=|\zeta|t\). The tangential Fourier transform
of the Laplacian, with coefficients evaluated at \(y=0\), is
\[
|\zeta|^2(\partial_s^2-1)
-|\zeta|\left(H_\xi^{\partial\Omega}\partial_s
+\frac{2s\mathrm{II}_\xi(\zeta,\zeta)}{|\zeta|^2}\right)
+O(1)\partial_s+O(1),
\qquad |\zeta|\to\infty,
\]
uniformly for $s$ in compact subsets of $[0,\infty)$.
Let \(w\) denote the decaying Fourier solution of the transformed
extension problem, normalized by \(w(0)=1\) and \(w'(0)=0\).
By Lemma~\ref{t2.2},

\begin{equation}\label{ws}
w(s)=(1+s)e^{-s}+|\zeta|^{-1}w_1(s)+O(|\zeta|^{-2}).
\end{equation}
Substituting this expansion into the transformed equation in \eqref{e1.3}, we obtain
\begin{equation*}
(\partial_s^2-1)^2w_1
=4\left(H_\xi^{\partial\Omega}
-\frac{3\mathrm{II}_\xi(\zeta,\zeta)}{|\zeta|^2}s\right)e^{-s},\quad
w_1(0)=w_1'(0)=0,\quad \lim_{s\to\infty}w_1(s)=0.
\end{equation*}
Solving this problem yields
\[
w_1(s)=\frac{e^{-s}}2
\left[
\left(H_\xi^{\partial\Omega}
-\frac{3\mathrm{II}_\xi(\zeta,\zeta)}{|\zeta|^2}\right)s^2
-\frac{\mathrm{II}_\xi(\zeta,\zeta)}{|\zeta|^2}s^3
\right].
\]
Substituting this expression into \eqref{ws} and using $\partial_\nu=-\partial_t$, we evaluate $\partial_t\Delta w$ at $t=0$.
Since $\partial_y p_3(0,\zeta)=0$, the standard and Weyl symbols have the same homogeneous terms of degrees three and two at $y=0$. Thus,
\begin{equation}
p_3(0,\zeta)=2|\zeta|^3,
\qquad p_2(0,\zeta)=4\mathrm{II}_\xi(\zeta,\zeta)-2H_\xi^{\partial\Omega}|\zeta|^2.
\label{e2.4}
\end{equation}
Comparing the homogeneous terms of degrees $0$ and $-1$ in the Weyl
composition formula for $\mathcal P\mathcal P^{-1}=I$, we obtain
\begin{equation*}
p_3r_{-3}=1,\qquad
p_3r_{-4}+p_2r_{-3}
+\frac{1}{2i}\sum_{j=1}^3
\bigl(\partial_{\zeta_j}p_3\,\partial_{y_j}r_{-3}
-\partial_{y_j}p_3\,\partial_{\zeta_j}r_{-3}\bigr)=0.
\end{equation*}
Since \(r_{-3}=p_3^{-1}\), the sum in the second identity vanishes.
Thus \(r_{-4}=-p_2p_3^{-2}\), and \eqref{e2.4} implies
\begin{equation}
r_{-3}(0,\zeta)=\frac1{2|\zeta|^3},
\qquad
r_{-4}(0,\zeta)=\frac{H_\xi^{\partial\Omega}|\zeta|^2-2\mathrm{II}_\xi(\zeta,\zeta)}{2|\zeta|^6}.
\label{e2.5}
\end{equation}
We use the inverse Fourier transform $\mathcal F^{-1}f(v)=(2\pi)^{-3}
\int_{\mathbb R^3}e^{iv\cdot\zeta}f(\zeta)\,d\zeta.$
For the homogeneous distributional extensions \cite[Section~2]{r38} of
\(|\zeta|^{-4}\) and \(|\zeta|^{-6}\), we have
\[
\mathcal F^{-1}(|\zeta|^{-4})=-\frac{|v|}{8\pi},
\qquad
\mathcal F^{-1}(|\zeta|^{-6})=\frac{|v|^3}{96\pi},
\]
Consequently, \eqref{e2.5} yields
\begin{equation*}
\mathcal F^{-1}(r_{-4}(0,\cdot))(v)
=-\frac{H_\xi^{\partial\Omega}|v|}{16\pi}
+\frac1{96\pi}\sum_{\alpha,\beta=1}^3
\mathrm{II}_\xi^{\alpha\beta}
\partial_{v_\alpha}\partial_{v_\beta}|v|^3=K_{-4,\xi}(v),
\end{equation*}
where \(\mathrm{II}_\xi^{\alpha\beta}\) are the components of the
second fundamental form in the chosen orthonormal frame.
We construct the local kernels from the Weyl symbol terms $r_{-3},\ldots,r_{-8}$ evaluated at $(y+z)/2$.
To express the kernels with respect to $d\sigma$, we multiply each kernel by the factor
\[
(\det g_\xi(y)\det g_\xi(z))^{-\frac{1}{4}}=1+O(|y|^2+|z|^2).
\]
Using a partition of unity and symmetrization, we write
\[
G_{\partial\Omega}=S^{(0)}+H^{(0)},
\]
where $S^{(0)}$ is supported near the diagonal and $H^{(0)}$
is the kernel of a remainder operator of order $-9$, hence
belongs to $C^4$~\cite{r17,r24}.
The local expansion of \(S^{(0)}\) shows that
\[
c_{\mathrm{reg}}(\xi)
=\lim_{z\to\xi}\left[
S^{(0)}(\xi,z)-\Gamma(\exp_\xi^{-1}z)\right],
\qquad c_{\mathrm{reg}}\in C^4(\partial\Omega).
\]
Choose a symmetric \(A_{\mathrm{reg}}\in C^4(\partial\Omega\times\partial\Omega)\),
supported near the diagonal, such that
\[
A_{\mathrm{reg}}(x,z)
=\tfrac12\bigl(c_{\mathrm{reg}}(x)+c_{\mathrm{reg}}(z)\bigr)
\quad\text{for }x,z\text{ sufficiently close},
\]
and set \(S=S^{(0)}-A_{\mathrm{reg}}\), \(H=H^{(0)}+A_{\mathrm{reg}}\).
The smooth part of the local expansion of \(S\) vanishes on the
diagonal. By symmetry, its first derivatives vanish there as well.
In a fixed coordinate chart, set $s=(u+v)/2$, $w=u-v$, and
$\widetilde G(s,w)=G_{\partial\Omega}(s+w/2,s-w/2)$,
with both arguments understood in that chart. For $0<|w|\ll1$,
the symbol expansion implies
\begin{equation}\label{e2.midpoint-bounds}
|\partial_s^\gamma\partial_w^\eta\widetilde G(s,w)|
\le C\begin{cases}
1+|\log|w||,&|\eta|=0,\\
|w|^{-|\eta|},&|\eta|\ge1,
\end{cases}
\qquad |\gamma|+|\eta|\le2.
\end{equation}
Indeed, a symbol term of degree $-3-j$ yields a kernel homogeneous
of degree $j$ in $w$, up to a term $p_j(s,w)\log|w|$ and a smooth
remainder, where $p_j$ is polynomial and homogeneous of degree $j$
in $w$. Differentiation in $s$ preserves these degrees, while
differentiation in $w$ lowers them by one. The density factors and
the $C^4$ remainder satisfy the same bounds.
In normal coordinates, Taylor's formula yields
\begin{equation*}
r_{-3}\left(\frac{y+z}{2},\zeta\right)-r_{-3}(0,\zeta)
=O\bigl((|y|+|z|)^2\bigr)|\zeta|^{-3},\quad
r_{-4}\left(\frac{y+z}{2},\zeta\right)-r_{-4}(0,\zeta)
=O(|y|+|z|)|\zeta|^{-4}.
\end{equation*}
These bounds and the quadratic vanishing of the density factors
persist under two center derivatives, since $g_\xi(0)=I$ and
$\partial_y g_\xi(0)=0$ for every $\xi$.
Moreover,
\[
y\cdot\nabla_y+z\cdot\nabla_z=s\cdot\nabla_s+w\cdot\nabla_w
\]
preserves the Taylor vanishing orders and the powers of $|w|$.
Thus, for $|\alpha|+j\le2$,
\begin{equation}\label{e2.kernel-dilation}
\begin{aligned}
\left|\partial_\xi^\alpha
(y\cdot\nabla_y+z\cdot\nabla_z)^j E_\xi(y,z)\right|
&\le C(|y|+|z|)^2\bigl(1+|\log|y-z||\bigr)\\
&\quad+C(|y|+|z|)|y-z|
+C|y-z|^2\bigl(1+|\log|y-z||\bigr)\\
&\le C(|y|+|z|)^2\bigl(1+|\log|y-z||\bigr),
\end{aligned}
\end{equation}
uniformly in $\xi$.
This proves \eqref{e2.6}.
Finally, \eqref{e2.6}, $K_{-4,\xi}(y)=O(|y|)$, and $d(\xi,\exp_\xi y)=|y|$ show that the limit in \eqref{e1.7} exists and
\[
r(\xi)=\lim_{y\to0}\bigl[G_{\partial\Omega}(\xi,\exp_\xi y)-\Gamma(y)\bigr]=H(\xi,\xi).
\]
Since $H$ is symmetric and of class $C^4$, we obtain $r\in C^4$ and
\[
\nabla r(\xi)=\nabla_yH_\xi(0,0)+\nabla_zH_\xi(0,0)
=2\nabla_yH_\xi(0,0),
\]
which proves \eqref{e2.10}.
To prove \eqref{e2.transported-green}--\eqref{e2.transported-green-gradient},
write $u=\Phi_\xi x$, $v=\Phi_\xi z$ in a common chart. Uniform
$C^1$ bounds and the mean-value theorem imply
\[
|\partial_\xi^\alpha s|\le C,\qquad
|\partial_\xi^\alpha w|\le C|w|,\qquad 1\le|\alpha|\le2.
\]
Consequently, the chain rule and \eqref{e2.midpoint-bounds} yield
\[
\begin{aligned}
|\partial_\xi^\alpha G_{\partial\Omega}(\Phi_\xi x,\Phi_\xi z)|
&\le C\sum_{|\gamma|+|\eta|\le|\alpha|}
|w|^{|\eta|}|\partial_s^\gamma\partial_w^\eta\widetilde G|
\le C(1+|\log|w||),&&|\alpha|\le2,\\
|\nabla_z\partial_\xi^\alpha G_{\partial\Omega}(\Phi_\xi x,\Phi_\xi z)|
&\le C|w|^{-1},&&|\alpha|\le1.
\end{aligned}
\]
Here $|w|\asymp d(\Phi_\xi x,\Phi_\xi z)$; away from the diagonal,
both estimates follow from smoothness.
\end{proof}
It follows from Proposition \ref{t2.3} that
\begin{equation}\label{e2.green-log-bound}
|G_{\partial\Omega}(x,z)|\le C(1+|\log d(x,z)|),\qquad x\ne z.
\end{equation}
By smooth approximation and H\"older's inequality, the identity
\begin{equation}\label{e2.green-representation}
\mathcal P^{-1}F(x)
=\int_{\partial\Omega}G_{\partial\Omega}(x,z)F(z)\,d\sigma(z),
\qquad x\in\partial\Omega,
\end{equation}
holds for every $F\in L^p(\partial\Omega)$, $p>1$.

\subsection{The model bubble and its linearization}\label{s2.2}

\begin{lemma}
The functions in \eqref{e1.5} and \eqref{e1.6} satisfy
\begin{align}
\int_{\mathbb R^3}q\,dy&=1,
\qquad \int_{\mathbb R^3}qz_k\,dy=0,\quad 0\le k\le3,
\label{e2.16}\\
\int_{\mathbb R^3}qz_kz_l\,dy&=\tfrac14\delta_{kl},\quad 0\le k,l\le3,
\label{e2.17}\\
\int_{\mathbb R^3}q(y)y_\alpha z_\beta(y)\,dy
&=\frac{\pi^2}{12a^4}\delta_{\alpha\beta},\quad 1\le\alpha,\beta\le3.
\label{e2.18}
\end{align}
Here $\delta_{kl}$ and $\delta_{\alpha\beta}$ are Kronecker symbols.
\end{lemma}
\begin{proof}
By oddness, \(\int_{\mathbb R^3}qz_k\,dy=0\) for \(1\le k\le3\) and
$\int_{\mathbb R^3}qz_kz_l\,dy=0$ for $0\le k,l\le3$, $k\ne l$.
Using polar coordinates with $r=a|y|$, we obtain
\[
\begin{aligned}
\int_{\mathbb R^3}q\,dy
=&\frac{4\pi}{a^3}\int_0^\infty\frac{r^2}{(1+r^2)^3}\,dr
=\frac{\pi^2}{4a^3}=1,\quad
\int_{\mathbb R^3}qz_0\,dy
=\frac{4\pi}{a^3}\int_0^\infty
\frac{r^2(1-r^2)}{(1+r^2)^4}\,dr=0,\\
\int_{\mathbb R^3}qz_k^2\,dy
&=\frac{16\pi}{3a^3}\int_0^\infty\frac{r^4}{(1+r^2)^5}\,dr
=\frac14,\quad 1\le k\le3.
\end{aligned}
\]
Since $\sum_{k=0}^3z_k^2=1$, we have
$\int_{\mathbb R^3}qz_0^2\,dy=1/4$.
This proves \eqref{e2.16}--\eqref{e2.17}.
Finally, radial symmetry yields
\[
\int_{\mathbb R^3}q(y)y_\alpha z_\beta(y)\,dy
=\frac{8\pi\delta_{\alpha\beta}}{3a^4}
\int_0^\infty\frac{r^4}{(1+r^2)^4}\,dr
=\frac{\pi^2}{12a^4}\delta_{\alpha\beta},
\qquad 1\le\alpha,\beta\le3.
\]
\end{proof}

Let \(\pi_0:\mathbb R^3\to\mathbb S^3\setminus\{N\}\) be the inverse
stereographic projection from \(N=(0,0,0,1)\), and set
\(\pi_a(y)=\pi_0(ay)\). Then
\begin{equation}\label{hengdengshi}
|\pi_a(y)-\pi_a(z)|
=\frac{2a|y-z|}{\sqrt{1+a^2|y|^2}\sqrt{1+a^2|z|^2}}.
\end{equation}
Under the change of variables $\vartheta=\pi_a(y)$, we have
\begin{equation}
q(y)\,dy=\frac1{2\pi^2}\,d\vartheta,
\label{e2.19}
\end{equation}
where $d\vartheta$ denotes the standard volume measure on $\mathbb S^3$.
On $L^2(\mathbb R^3,q\,dy)$, define
\begin{equation}\label{e2.model-operators}
(\mathscr Kf)(y)=-6\int_{\mathbb R^3}\log|y-z|q(z)f(z)\,dz,
\qquad \Pi f=f-\int_{\mathbb R^3}qf\,dy.
\end{equation}
By \eqref{e2.16}, $\Pi$ is the orthogonal projection onto the mean-zero
subspace. Since $q(y)\le C(1+|y|)^{-6}$,
\[
\iint_{\mathbb R^3\times\mathbb R^3}
q(y)q(z)|\log|y-z||^2\,dy\,dz
\le C\int_{\mathbb R^3}q(y)(1+\log(2+|y|))^2\,dy<\infty,
\]
so $\mathscr K$ is Hilbert--Schmidt and self-adjoint.
For $\theta=\pi_a(y)$ and $\vartheta=\pi_a(z)$, \eqref{hengdengshi} gives
\[
\log|y-z|=\log|\theta-\vartheta|-\log(2a)
+\tfrac12\log(1+a^2|y|^2)+\tfrac12\log(1+a^2|z|^2).
\]
Projection by $\Pi$ on both sides removes the terms depending on at
most one variable. The Funk--Hecke calculation
in~\cite[proof of Theorem~5]{CarlenLoss1992}, with measure
$d\vartheta/(2\pi^2)$, therefore yields
\begin{equation}\label{e2.25}
\Pi\mathscr K\Pi=6P_3^{-1},\qquad
\spec(\Pi\mathscr K\Pi)\setminus\{0\}
=\left\{\frac{6}{\ell(\ell+1)(\ell+2)}:\ell\ge1\right\},
\end{equation}
on the mean-zero subspace, under the identification induced by $\pi_a$.
Here $P_3$ is the self-adjoint operator on $L^2(\mathbb S^3)$ with
domain $H^3(\mathbb S^3)$ and eigenvalues $\ell(\ell+1)(\ell+2)$
on spherical harmonics of degree $\ell\ge0$; $P_3^{-1}$ is taken
on the mean-zero subspace.
\begin{lemma}\label{t2.8}
For \(\Gamma\) defined in \eqref{e2.2}, it holds

\begin{equation}\label{e2.model-potential}
\int_{\mathbb R^3}\Gamma(y-z)q(z)\,dz
=u_0(y)+\frac1{4\pi^2}\log a.
\end{equation}
\end{lemma}
\begin{proof}
By \eqref{hengdengshi} and \eqref{e2.19}, since the measure on \(\mathbb S^3\) is rotational invariant, we have
\[
\int_{\mathbb R^3}\log|y-z|q(z)\,dz
=\frac12\log(1+a^2|y|^2)
+\int_{\mathbb R^3}\log|z|q(z)\,dz.
\]
Using $\int_{\mathbb R^3}q\,dz=1$ and polar coordinates with $r=a|z|$,
we obtain
\[
\int_{\mathbb R^3}\log|z|q(z)\,dz
=-\log a+\frac{4\pi}{a^3}
\int_0^\infty\frac{r^2\log r}{(1+r^2)^3}\,dr
=-\log a.
\]
The conclusion follows from the definitions of \(\Gamma\) (see \eqref{e2.2}) and \(u_0\) (see \eqref{e1.5}).
\end{proof}

\subsection{Construction of the approximate solution}\label{s3}


Let $\exp_\xi$ be the Riemannian exponential map at $\xi\in\partial\Omega$ and identify $T_\xi(\partial\Omega)$ with $\mathbb R^3$ through a local orthonormal frame.
Set $B_r(\xi)=\{x\in\partial\Omega:d(x,\xi)<r\}$.
On each fixed compact subset of $\mathcal M_m$, choose $r_0>0$ so that
$8r_0$ is smaller than the injectivity radius, the balls
$B_{4r_0}(\xi_i)$ are pairwise disjoint, and $d$ coincides with the
geodesic distance on each $B_{4r_0}(\xi_i)\times B_{4r_0}(\xi_i)$.
For $\delta>0$, set
\begin{equation}\label{e3.1}
\widetilde q_{\delta,\xi}(x)
=\chi(d(x,\xi))\delta^{-3}q\left(\frac{\exp_\xi^{-1}x}{\delta}\right),\quad
b_{\delta,\xi}
=\left(\int_{\partial\Omega}\widetilde q_{\delta,\xi}\,d\sigma\right)^{-1},
\quad q_{\delta,\xi}=b_{\delta,\xi}\widetilde q_{\delta,\xi},
\end{equation}
with $\widetilde q_{\delta,\xi}=0$ in  $\partial\Omega\setminus B_{2r_0}(\xi)$.
Radial symmetry of $q$ makes $\widetilde q_{\delta,\xi}$ and $q_{\delta,\xi}$ independent of the frame, and
$\int_{\partial\Omega}q_{\delta,\xi}\,d\sigma=1$.
The projected bubble is
\begin{equation}\label{defW}
W_{\delta,\xi}(x)
=\bigl(\mathcal P^{-1}q_{\delta,\xi}\bigr)(x)
=\int_{\partial\Omega}G_{\partial\Omega}(x,z)q_{\delta,\xi}(z)\,d\sigma(z),
\qquad x\in\partial\Omega.
\end{equation}
For $\xi\in\mathcal M_m$, set
\[
\mathcal R_i(\xi)=r(\xi_i)+\sum_{j\ne i}G_{\partial\Omega}(\xi_i,\xi_j),
\]
\begin{equation}\label{e3.10}
\delta_i(\varepsilon,\xi)
=a^2\varepsilon^{\frac{1}{3}}h(\xi_i)^{\frac{1}{3}}e^{8\pi^2\mathcal R_i(\xi)},
\qquad V_\xi=\sum_{i=1}^mW_{\delta_i(\varepsilon,\xi),\xi_i}.
\end{equation}


\begin{lemma}\label{t3.1}
As $\delta\to0$, the following estimates hold uniformly for
$\xi\in\partial\Omega$:
\begin{equation}\label{e3.2}
b_{\delta,\xi}=1+O(\delta^2),\quad
\int_{B_{2r_0}(\xi)}
\exp_\xi^{-1}(x)q_{\delta,\xi}(x)\,d\sigma(x)
=O(\delta^2),\quad
\int_{\partial\Omega}d(x,\xi)^2q_{\delta,\xi}(x)\,d\sigma(x)
=O(\delta^2).
\end{equation}
\end{lemma}

\begin{proof}
In normal coordinates, we have
\[
d\sigma(\exp_\xi(\delta y))
=\delta^3\bigl(1+O(\delta^2|y|^2)\bigr)\,dy,
\qquad |y|<\frac{2r_0}{\delta},
\]
uniformly in $\xi$.
Since $1-\chi(\delta|y|)\le C\delta^2|y|^2$,
$\int_{\mathbb R^3}q\,dy=1$ and
$\int_{\mathbb R^3}|y|^2q\,dy<\infty$, we obtain
$b_{\delta,\xi}^{-1}=1+O(\delta^2)$ and hence
$b_{\delta,\xi}=1+O(\delta^2)$.
Changing variables $x=\exp_\xi(\delta y)$, we obtain
\[
\int_{\partial\Omega}d(x,\xi)^2q_{\delta,\xi}(x)\,d\sigma(x)
\le C\delta^2\int_{\mathbb R^3}|y|^2q(y)\,dy
=O(\delta^2).
\]
Finally, by radial symmetry,
$\int_{\mathbb R^3}y\chi(\delta|y|)q(y)\,dy=0$, so
\begin{equation*}
\left|\int_{B_{2r_0}(\xi)}
\exp_\xi^{-1}(x)q_{\delta,\xi}(x)\,d\sigma(x)\right|
\le C\delta^3\int_{|y|<2r_0/\delta}|y|^3q(y)\,dy\le C\delta^2\int_{\mathbb R^3}|y|^2q(y)\,dy
=O(\delta^2).
\end{equation*}
\end{proof}
Define $B_\xi(y)=\int_{\mathbb R^3}K_{-4,\xi}(y-s)q(s)\,ds.$
Since \(K_{-4,\xi}\) and \(q\) are even, we obtain
\[
B_\xi(-y)=B_\xi(y).
\]
Formula \eqref{e2.curvature-kernel}, together with
\(\int_{\mathbb R^3}(1+|s|)q(s)\,ds<\infty\), gives
\begin{equation}\label{e3.curvature-bound}
|\partial_\xi^\alpha B_\xi(y)|
+(1+|y|)|\nabla_y\partial_\xi^\alpha B_\xi(y)|
\le C(1+|y|),\qquad |\alpha|\le2,
\end{equation}
where \(\xi\)-derivatives are taken in a smooth local orthonormal
frame with \(y\) fixed.
\begin{proposition}\label{t3.2}
For any fixed $L>0$ and $|y|\le L\delta^{-1/2}$, we have
\begin{equation}\label{e3.4}
\begin{split}
W_{\delta,\xi}(\exp_\xi(\delta y))
={}&u_0(y)-\frac{1}{4\pi^2}\log\delta
+\frac{1}{4\pi^2}\log a+r(\xi)\\
&+\delta B_\xi(y)
+\frac{\delta}{2}\langle\nabla r(\xi),y\rangle
+O_L\bigl(\delta^2(1+|y|^2)|\log\delta|\bigr).
\end{split}
\end{equation}
Moreover, for every fixed $r>0$,
\begin{equation}\label{e3.6}
W_{\delta,\xi}
=G_{\partial\Omega}(\cdot,\xi)+O_r(\delta^2)
\quad\text{in }C^3\bigl(\{x\in\partial\Omega:d(x,\xi)\ge r\}\bigr).
\end{equation}
Both estimates hold uniformly for $\xi\in\partial\Omega$ as $\delta\to0$.
\end{proposition}

\begin{proof}
Fix $L>0$ and assume $|y|\le L\delta^{-1/2}$.
After reducing \(r_0\) if necessary, \eqref{e2.6} applies on the
support of \(q_{\delta,\xi}\) for sufficiently small \(\delta\).
We set \(z=\exp_\xi(\delta s)\), \(|s|<2r_0/\delta\).
By \eqref{e2.model-potential} and the scaling identity
$K_{-4,\xi}(\delta(y-s))=\delta K_{-4,\xi}(y-s)$, we obtain
\begin{equation}\label{n1}
\int_{\mathbb R^3}\Gamma(\delta(y-s))q(s)\,ds
=u_0(y)-\frac1{4\pi^2}\log\delta+\frac1{4\pi^2}\log a,\quad
\int_{\mathbb R^3}K_{-4,\xi}(\delta(y-s))q(s)\,ds
=\delta B_\xi(y).
\end{equation}
In view of \eqref{e2.6}, we have
\begin{equation}\label{n2}
\int_{|s|<\frac{2r_0}{\delta}}|E_\xi(\delta y,\delta s)|q(s)\,ds
\le C\delta^2\int_{\mathbb R^3}(|y|+|s|)^2
\bigl(1+|\log\delta|+|\log|y-s||\bigr)q(s)\,ds
\le C_L\delta^2(1+|y|^2)|\log\delta|.
\end{equation}
To estimate the cutoff errors, we use the decay of $q$ and obtain
\begin{align}\label{n3}
\left|\int_{\mathbb R^3}
\bigl[1-\chi(\delta|s|)\chi(\delta|y-s|)\bigr]
\Gamma(\delta(y-s))q(s)\,ds\right|
&\le C_L\delta^3|\log\delta|,\\
\left|\int_{\mathbb R^3}
\bigl[1-\chi(\delta|s|)\chi(\delta|y-s|)\bigr]
K_{-4,\xi}(\delta(y-s))q(s)\,ds\right|
&\le C_L\delta^3.
\end{align}
From \eqref{e3.2} and $\sqrt{\det g_\xi(\delta s)}=1+O(\delta^2|s|^2)$, we get
\begin{equation}\label{n4}
\left|\int_{|s|<\frac{2r_0}{\delta}}
\chi(\delta|s|)
\bigl[b_{\delta,\xi}\sqrt{\det g_\xi(\delta s)}-1\bigr]
S_\xi(\delta y,\delta s)q(s)\,ds\right|\le C_L\delta^2(1+|y|^2)|\log\delta|.
\end{equation}
Taylor's formula, \eqref{e3.2}, and \eqref{e2.10} yield
\begin{equation}\label{n5}
\int_{\partial\Omega}H(\exp_\xi(\delta y),z)q_{\delta,\xi}(z)\,d\sigma(z)
=r(\xi)+\frac{\delta}{2}\langle\nabla r(\xi),y\rangle
+O_L\bigl(\delta^2(1+|y|^2)\bigr).
\end{equation}
Combining \eqref{n1}--\eqref{n5}, we arrive at \eqref{e3.4}.

To prove \eqref{e3.6}, let $\chi_{r,\xi}$ be a smooth cutoff supported in $B_{r/2}(\xi)$ and equal to one on $B_{r/4}(\xi)$. We write
\[
W_{\delta,\xi}=\mathcal P^{-1}(\chi_{r,\xi} q_{\delta,\xi})
+\mathcal P^{-1}((1-\chi_{r,\xi})q_{\delta,\xi}).
\]
For $0<\alpha<1$, the decay of $q$ and the boundedness of $\mathcal P^{-1}:C^\alpha\to C^{3,\alpha}$ imply
\[
\|(1-\chi_{r,\xi})q_{\delta,\xi}\|_{C^\alpha(\partial\Omega)}\le C_r\delta^3,
\qquad
\|\mathcal P^{-1}((1-\chi_{r,\xi})q_{\delta,\xi})\|_{C^{3,\alpha}(\partial\Omega)}\le C_r\delta^3.
\]
If $d(x,\xi)\ge r$, the kernel $G_{\partial\Omega}(x,z)$ is smooth
for $z\in\supp\chi_{r,\xi}$. Taylor's formula, the moment estimates
in \eqref{e3.2}, and the bound on $(1-\chi_{r,\xi})q_{\delta,\xi}$ yield
\[
\begin{aligned}
\mathcal P^{-1}(\chi_{r,\xi} q_{\delta,\xi})(x)
&=G_{\partial\Omega}(x,\xi)
\int_{\partial\Omega}\chi_{r,\xi} q_{\delta,\xi}\,d\sigma
+O_r(\delta^2)\\
&=G_{\partial\Omega}(x,\xi)
\left(1-\int_{\partial\Omega}(1-\chi_{r,\xi})
q_{\delta,\xi}\,d\sigma\right)+O_r(\delta^2)\\
&=G_{\partial\Omega}(x,\xi)+O_r(\delta^2)
\end{aligned}
\]
in $C^3(\{x\in\partial\Omega:d(x,\xi)\ge r\})$.
Therefore,
\[
\|W_{\delta,\xi}-G_{\partial\Omega}(\cdot,\xi)\|
_{C^3(\{x\in\partial\Omega:d(x,\xi)\ge r\})}
\le C_r\delta^2+C_r\delta^3\le C_r\delta^2,
\]
which proves \eqref{e3.6}.
\end{proof}

We next establish pointwise estimates for $W_{\delta,\xi}$
and its derivatives on $\partial\Omega$.

\begin{lemma}\label{t3.3}
There exist constants $\delta_0>0$ and $C>0$  such that
for every $0<\delta<\delta_0$ and $x,\xi\in\partial\Omega$,
\begin{align}
&\left|W_{\delta,\xi}(x)+\frac1{4\pi^2}\log(\delta+d(x,\xi))\right|\le C,
\label{e3.7}\\
&|\nabla_x^jW_{\delta,\xi}(x)|\le\frac C{(\delta+d(x,\xi))^j},
\qquad 1\le j\le3,
\label{e3.8}\\
&e^{24\pi^2W_{\delta,\xi}(x)}\le\frac C{(\delta+d(x,\xi))^6}.
\label{e3.9}
\end{align}
\end{lemma}

\begin{proof}
For $d(x,\xi)\ge r_0$, it follows from \eqref{e3.6} that
\[
|W_{\delta,\xi}(x)|+\sum_{j=1}^3|\nabla_x^jW_{\delta,\xi}(x)|\le C.
\]
Combining this bound with $\delta+d(x,\xi)\asymp1$, we obtain
\eqref{e3.7}--\eqref{e3.8} for $d(x,\xi)\ge r_0$.
It remains to consider \(d(x,\xi)<r_0\). Set
\(s=\delta+d(x,\xi)\).
From \eqref{e3.1},
\begin{equation}\label{n7}
|\nabla_z^k q_{\delta,\xi}(z)|
\le\frac{C\delta^3}{(\delta+d(z,\xi))^{6+k}},
\qquad k=0,1.
\end{equation}
Let $Y=\delta^{-1}\exp_\xi^{-1}(x)$. Using Proposition~\ref{t2.3},
the normalization $\int_{\partial\Omega}q_{\delta,\xi}\,d\sigma=1$, and the change of variables
$z=\exp_\xi(\delta y)$, we obtain
\begin{equation}\label{n6}
\left|W_{\delta,\xi}(x)+\frac1{4\pi^2}\log s\right|
\le C+C\int_{\mathbb R^3}q(y)
\left|\log\frac{|Y-y|}{1+|Y|}\right|\,dy.
\end{equation}
We split the integral at $|Y-y|=(1+|Y|)/2$ and estimate
\[
\begin{aligned}
\int_{\mathbb R^3}q(y)\left|\log\frac{|Y-y|}{1+|Y|}\right|\,dy
\le\frac{C}{(1+|Y|)^6}
\int_0^{\frac{1+|Y|}{2}}t^2\log\frac{1+|Y|}{t}\,dt+C\int_{\mathbb R^3}q(y)\bigl(1+\log(1+|y|)\bigr)\,dy
\le C,
\end{aligned}
\]
which together with \eqref{n6} establishes \eqref{e3.7}.
Choose $c>0$ sufficiently small. It follows from  \eqref{n7}  that
\[
\|q_{\delta,\xi}\|_{L^\infty(B_{cs}(x))}
+s\|\nabla q_{\delta,\xi}\|_{L^\infty(B_{cs}(x))}
\le Cs^{-3}.
\]
The symbol expansion in the proof of Proposition~\ref{t2.3},
with leading terms \eqref{e2.5}, implies
\begin{equation}\label{e3.green-derivatives}
|\nabla_x^jG_{\partial\Omega}(x,z)|\le Cd(x,z)^{-j},
\qquad x\ne z,\quad 1\le j\le3.
\end{equation}
Hence, for \(j=1,2\),

\begin{equation*}
|\nabla_x^jW_{\delta,\xi}(x)|
\le Cs^{-3}\int_0^{cs}t^{2-j}\,dt
+Cs^{-j}\int_{\partial\Omega\setminus B_{cs}(x)}q_{\delta,\xi}\,d\sigma\le Cs^{-j}.
\end{equation*}
For $j=3$, let $z=\exp_xv$ and recall $\Gamma$ from \eqref{e2.2}.
Using the homogeneous symbol terms of $\mathcal P^{-1}$ in
\eqref{e2.5}, we express the local distribution kernel of
$\nabla_x^3\mathcal P^{-1}$ with respect to $dv$ as $-\operatorname{p.v.}D^3\Gamma(v)
+O\bigl(|v|^{-2}(1+|\log|v||)\bigr).$
Here $\operatorname{p.v.}$ denotes the Cauchy principal value.
Since the third Euclidean derivative $D^3\Gamma$ is odd and
homogeneous of degree $-3$,
\[
|D^3\Gamma(v)|\le C|v|^{-3},\qquad
\operatorname{p.v.}\int_{|v|<cs}D^3\Gamma(v)\,dv=0.
\]
Since $|q_{\delta,\xi}(\exp_xv)-q_{\delta,\xi}(x)|\le Cs^{-4}|v|$
for $|v|<cs$, we have
\[
\begin{aligned}
\left|\operatorname{p.v.}\int_{|v|<cs}
D^3\Gamma(v)q_{\delta,\xi}(\exp_xv)\,dv\right|
&=\left|\int_{|v|<cs}D^3\Gamma(v)
\bigl[q_{\delta,\xi}(\exp_xv)-q_{\delta,\xi}(x)\bigr]\,dv\right|\\
&\quad\le Cs^{-4}\int_{|v|<cs}|v|^{-2}\,dv
\le Cs^{-3}.
\end{aligned}
\]
Therefore, $
|\nabla_x^3W_{\delta,\xi}(x)|
\le Cs^{-3}
+Cs^{-3}\int_0^{cs}(1+|\log t|)\,dt
+Cs^{-3}\int_{\partial\Omega\setminus B_{cs}(x)}q_{\delta,\xi}\,d\sigma
\le Cs^{-3}.$
This proves \eqref{e3.8}. Finally, \eqref{e3.9} follows by exponentiating \eqref{e3.7}.
\end{proof}

Let \(K\) and \(D\) be as in \eqref{e1.9}. Fix a compact parameter set
\(\mathcal K\subset\mathcal M_m\) such that
$
K\subset D\subset\overline D\subset\operatorname{int}\mathcal K.
$
For $\xi\in\mathcal K$, write $\delta_i=\delta_i(\varepsilon,\xi)$ and
$V=V_\xi$ as in \eqref{e3.10}, and set
\[
A_i(\xi)=\frac12\nabla r(\xi_i)
+\sum_{j\ne i}\nabla_xG_{\partial\Omega}(\xi_i,\xi_j),
\qquad 1\le i\le m.
\]
Fix \(L>0\). For sufficiently small \(\varepsilon\),
\(\exp_{\xi_i}(\delta_i y)\in B_{r_0}(\xi_i)\) when
\(|y|\le L\delta_i^{-1/2}\).
Since the balls \(B_{4r_0}(\xi_i)\) are disjoint,
equations \eqref{e3.4} and \eqref{e3.6}, together with Taylor expansion of the interaction terms, yield
\begin{equation}
V(\exp_{\xi_i}(\delta_i y))
=u_0(y)-\frac1{4\pi^2}\log\delta_i+\frac1{4\pi^2}\log a+\mathcal R_i(\xi)+\delta_iB_{\xi_i}(y)+\delta_i\langle A_i(\xi),y\rangle+E_{i,0}(y),
\label{e3.11}
\end{equation}
where
\[
|E_{i,0}(y)|\le C_L\delta_i^2(1+|y|^2)|\log\varepsilon|,
\qquad |y|\le L\delta_i^{-\frac{1}{2}},
\]
uniformly for $\xi\in\mathcal K$.
Using the definition of $\delta_i$ in \eqref{e3.10}, we rewrite
\eqref{e3.11} as
\begin{equation}
V(\exp_{\xi_i}(\delta_i y))
=u_0(y)-\frac1{8\pi^2}\log\delta_i-\frac1{24\pi^2}\log(\varepsilon h(\xi_i))+\delta_iB_{\xi_i}(y)+\delta_i\langle A_i(\xi),y\rangle+E_{i,0}(y).
\label{e3.12}
\end{equation}

Set
\[
\omega(x)=1+\sum_{i=1}^m
\frac{\delta_i^{-3}}{(1+\delta_i^{-1}d(x,\xi_i))^4},
\qquad \delta_{\min}=\min_{1\le i\le m}\delta_i,
\]
and define the residual by $
\mathcal R_\varepsilon(V)=\mathcal PV-\varepsilon he^{24\pi^2V}.$
For \(f\in L^\infty(\partial\Omega)\), let $
\|f\|_* =\|\omega^{-1}f\|_{L^\infty(\partial\Omega)}.$
For sufficiently small $\varepsilon>0$, we have $0<\delta_{\min}<1$ and
\begin{equation}\label{e3.weight-bounds}
\|\omega\|_{L^1(\partial\Omega)}
\le C\left(1+\int_0^\infty\frac{r^2}{(1+r)^4}\,dr\right)\le C,
\qquad
\|\omega\|_{L^\infty(\partial\Omega)}\le C\delta_{\min}^{-3}.
\end{equation}
Using \eqref{e2.green-log-bound}, \eqref{e3.weight-bounds}, and \eqref{e3.10}, we obtain
\begin{equation}\label{e2.14}
\begin{aligned}
\sup_{x\in\partial\Omega}
\int_{\partial\Omega}|G_{\partial\Omega}(x,z)|\omega(z)\,d\sigma(z)
&\le C\delta_{\min}^{-3}
\int_0^{\delta_{\min}}(1+|\log r|)r^2\,dr
+C(1+|\log\delta_{\min}|)\\
&\le C(1+|\log\delta_{\min}|)
\le C|\log\varepsilon|.
\end{aligned}
\end{equation}
We write $e_{i,1},e_{i,2},e_{i,3}$ for the smooth local orthonormal
frame associated with the normal coordinates at $\xi_i$.
For $1\le i\le m$ and $0\le k\le3$, define
\begin{equation}
Z_{i,k}(x)=\chi(d(x,\xi_i))z_k\left(\frac{\exp_{\xi_i}^{-1}x}{\delta_i}\right).
\label{e3.13}
\end{equation}
The functions $Z_{i,k}$ are extended by zero in $\partial\Omega\setminus B_{2r_0}(\xi_i)$.
Here \(d_i\mathcal F_m(\xi)\) denotes the differential with respect
to the \(i\)th center, so that
\[
d_i\mathcal F_m(\xi)[e_{i,k}]
=\langle\nabla_{\xi_i}\mathcal F_m(\xi),e_{i,k}\rangle,
\qquad 1\le k\le3.
\]
\begin{proposition}\label{t3.4}
Let $L>0$ be fixed. For sufficiently small $\varepsilon>0$,
\begin{equation}\label{e3.14}
\|\mathcal R_\varepsilon(V)\|_*
\le C\delta_*|\log\varepsilon|,
\qquad
\delta_*=\max_{1\le i\le m}\delta_i.
\end{equation}
For $1\le i\le m$, we have
\begin{equation}\label{e3.15}
\delta_i^3\mathcal R_\varepsilon(V)
\bigl(\exp_{\xi_i}(\delta_i y)\bigr)
=-24\pi^2\delta_iq(y)
\bigl[B_{\xi_i}(y)
+\langle\nabla_{\xi_i}\mathcal F_m(\xi),y\rangle\bigr]
+E_i(y),\qquad |y|\le L\delta_i^{-\frac{1}{2}},
\end{equation}
where
\begin{equation}\label{e3.17}
|E_i(y)|
\le C_L\delta_i^2|\log\varepsilon|(1+|y|)^{-4}.
\end{equation}
For \(1\le i\le m\) and \(1\le k\le3\),

\begin{equation}\label{e3.18}
\int_{\partial\Omega}
\mathcal R_\varepsilon(V)Z_{i,k}\,d\sigma
=
-\frac{2\pi^4}{a^4}\delta_i
d_i\mathcal F_m(\xi)[e_{i,k}]
+O\bigl(\delta_i^2|\log\varepsilon|\bigr).
\end{equation}
All estimates are uniform for $\xi\in\mathcal K$.
\end{proposition}
\begin{proof}
By \eqref{defW}, $\mathcal PV=\sum_{j=1}^m q_{\delta_j,\xi_j}$.
Recall from \eqref{e1.8} that
\[
A_i+\frac{1}{24\pi^2}\nabla\log h(\xi_i)
=\nabla_{\xi_i}\mathcal F_m(\xi).
\]
For $|y|\le L\delta_i^{-1/2}$, \eqref{e3.12} and Taylor
expansion give
\begin{equation*}
\delta_i^3\varepsilon h\bigl(\exp_{\xi_i}(\delta_i y)\bigr)
e^{24\pi^2V(\exp_{\xi_i}(\delta_i y))}=q(y)\Bigl\{
1+24\pi^2\delta_i
\bigl[B_{\xi_i}(y)
+\langle\nabla_{\xi_i}\mathcal F_m(\xi),y\rangle\bigr]
+O_L\bigl(\delta_i^2(1+|y|^2)|\log\varepsilon|\bigr)
\Bigr\}.
\end{equation*}
On this region, \(\chi(\delta_i|y|)=1\) and $q_{\delta_j,\xi_j}\big(\exp_{\xi_i}(\delta_i y)\big) = 0$ for $j\neq i.$ Hence
\[
\delta_i^3\mathcal PV\bigl(\exp_{\xi_i}(\delta_i y)\bigr)
=b_{\delta_i,\xi_i}q(y).
\]
Using \(b_{\delta_i,\xi_i}=1+O(\delta_i^2)\) and
\(q(y)(1+|y|)^2\le C(1+|y|)^{-4}\), we obtain
\eqref{e3.15}--\eqref{e3.17}.
We set \(L = 1\) for the rest of the proof.
By Lemma~\ref{t3.3} and \eqref{e3.10}, we obtain
\begin{equation}\label{huanxing}
\delta_i^3
\bigl(|\mathcal PV|+\varepsilon he^{24\pi^2V}\bigr)
\bigl(\exp_{\xi_i}(\delta_i y)\bigr)
\le C(1+|y|)^{-6}
\le C\delta_i(1+|y|)^{-4},
\qquad
\delta_i^{-\frac{1}{2}}\le |y|\le \frac{2r_0}{\delta_i}.
\end{equation}
On $\partial\Omega\setminus\bigcup_{i=1}^m B_{2r_0}(\xi_i)$, we have
\begin{equation}\label{xiaoqiuzhi}
|\mathcal R_\varepsilon(V)|\le C\delta_*^3.
\end{equation}
Since $
\omega\bigl(\exp_{\xi_i}(\delta_i y)\bigr)
\ge\delta_i^{-3}(1+|y|)^{-4}$ for $|y|\le2r_0/\delta_i,$ $\omega\ge1,$
equations \eqref{e3.15}--\eqref{e3.17}, \eqref{huanxing} and
\eqref{xiaoqiuzhi} imply \eqref{e3.14}.
To prove \eqref{e3.18}, fix $1\le k\le3$.
Estimate \eqref{huanxing} and
$|z_k(y)|\le C(1+|y|)^{-1}$ yield
\[
\int_{\partial\Omega\setminus B_{\sqrt{\delta_i}}(\xi_i)}
|\mathcal R_\varepsilon(V)Z_{i,k}|\,d\sigma
\le C\int_{|y|>\delta_i^{-\frac{1}{2}}}(1+|y|)^{-7}\,dy
\le C\delta_i^2.
\]
Since $qB_{\xi_i}$ is even and $z_k$ is odd, $\int_{|y|\le\delta_i^{-1/2}}
q(y)B_{\xi_i}(y)z_k(y)\,dy=0.$
The expansion \(\sqrt{\det g_{\xi_i}(\delta_i y)}=1+O(\delta_i^2|y|^2)\) yields
\[
\begin{aligned}
&\delta_i\left|\int_{|y|\le\delta_i^{-\frac{1}{2}}}
q(y)\bigl[B_{\xi_i}(y)+\langle\nabla_{\xi_i}\mathcal F_m(\xi),y\rangle\bigr]
z_k(y)\bigl[\sqrt{\det g_{\xi_i}(\delta_i y)}-1\bigr]\,dy\right|\\
&\qquad\le C\delta_i^3\int_{\mathbb R^3}
|y|^2q(y)(1+|y|)|z_k(y)|\,dy
\le C\delta_i^3.
\end{aligned}
\]
Thus \eqref{e3.15} and \eqref{e3.17} imply
\[
\begin{aligned}
\int_{\partial\Omega}\mathcal R_\varepsilon(V)Z_{i,k}\,d\sigma
={}&-24\pi^2\delta_i
\int_{|y|\le\delta_i^{-\frac{1}{2}}}
q(y)\langle\nabla_{\xi_i}\mathcal F_m(\xi),y\rangle
z_k(y)\,dy
+O\bigl(\delta_i^2|\log\varepsilon|\bigr).
\end{aligned}
\]
Since $
\delta_i\int_{|y|>\delta_i^{-1/2}}
q(y)|y|\,|z_k(y)|\,dy
\le C\delta_i^{5/2},$
we obtain \eqref{e3.18} by extending the integration domain
from $\{|y|\le\delta_i^{-1/2}\}$ to $\mathbb R^3$
and applying \eqref{e2.18}.
\end{proof}

\section{Finite-dimensional reduction}\label{s4}

\subsection{The projected linear problem}\label{s4.1}

For $\xi\in\mathcal K$, set $\Lambda_\varepsilon=|\log\varepsilon|$ and
$q_i=q_{\delta_i,\xi_i}$, $1\le i\le m$.
The linearization of \eqref{e1.4} at $V$ is
\[
L_V=\mathcal P-Q_V,\qquad Q_V=24\pi^2\varepsilon h e^{24\pi^2V}.
\]
The associated bilinear form on $H^{3/2}(\partial\Omega)$ is
\[
\mathcal B_V(u,v)=\langle\mathcal Pu,v\rangle
-\int_{\partial\Omega}Q_Vuv\,d\sigma.
\]
For the functions $Z_{i,k}$ defined in \eqref{e3.13}, set
\begin{equation}\label{e4.moment-space}
\begin{gathered}
\ell_{i,k}(u)=\int_{\partial\Omega}q_iZ_{i,k}u\,d\sigma,\\
X_\perp=\{u\in H^{\frac{3}{2}}(\partial\Omega):\ell_{i,k}(u)=0,\ 1\le i\le m,\ 0\le k\le3\}.
\end{gathered}
\end{equation}
A change of variables in normal coordinates, together with
\eqref{e2.16}--\eqref{e2.17} and Lemma~\ref{t3.1}, gives
\begin{equation}\label{e4.cutoff-moments}
\int_{\partial\Omega}q_iZ_{i,k}\,d\sigma=O(\delta_i^2),\qquad
\int_{\partial\Omega}q_iZ_{i,k}Z_{j,l}\,d\sigma
=\frac14\delta_{ij}\delta_{kl}+o(1),
\end{equation}
for $1\le i,j\le m$ and $0\le k,l\le3$.
Since $\mathcal PV=\sum_{i=1}^m q_i\le C\omega$, \eqref{e3.14} yields
\begin{equation}\label{e4.potential-bound}
0<Q_V=24\pi^2\biggl(\sum_{i=1}^m q_i-\mathcal R_\varepsilon(V)\biggr)
\le C\omega
\end{equation}
for sufficiently small $\varepsilon>0$.
\begin{lemma}\label{t4.1}
For \(p\ge2\) and \(u\in H^{3/2}(\partial\Omega)\),
\begin{equation}\label{e4.weighted-embedding}
 \|u\|_{L^{2p}}^2\le Cp\|u\|_{H^{\frac{3}{2}}}^2,\qquad
 \int_{\partial\Omega}\omega u^2\,d\sigma
 \le C\Lambda_\varepsilon\|u\|_{H^{\frac{3}{2}}}^2 .
\end{equation}
\end{lemma}
\begin{proof}
Assume $\|u\|_{H^{3/2}}=1$.
The critical exponential inequality~\cite{r3,r6} implies
$\int_{\partial\Omega}e^{cu^2}\,d\sigma\le C$ for some fixed $c>0$.
For every $p\ge2$, the pointwise bound
$|u|^{2p}\le (p/(ce))^p e^{cu^2}$ yields
\[
\|u\|_{L^{2p}}^2
\le \frac{p}{ce}\left(\int_{\partial\Omega}e^{cu^2}\,d\sigma\right)^{\frac{1}{p}}
\le Cp.
\]
Now take $p=\max\{2,\Lambda_\varepsilon\}$ and let
$p'=p/(p-1)$ be its H\"older conjugate.
Since $\delta_{\min}\asymp\varepsilon^{1/3}$ by \eqref{e3.10},
\[
\|\omega\|_{L^{p'}}
\le \|\omega\|_{L^1}^{\frac{1}{p'}}\|\omega\|_{L^\infty}^{\frac{1}{p}}
\le C\delta_{\min}^{-\frac{3}{p}}\le C.
\]
Applying H\"older's inequality, we obtain
\[
\int_{\partial\Omega}\omega u^2\,d\sigma
\le \|\omega\|_{L^{p'}}\|u\|_{L^{2p}}^2
\le Cp\le C\Lambda_\varepsilon.
\]
\end{proof}

\begin{lemma}\label{t4.2}
For \(1\le i\le m\), set
\(w_i=\mathcal P^{-1}q_i\) and \(W_i=w_i(\xi_i)\). Then
\begin{equation}\label{e4.bubble-norms}
 W_i=\frac{|\log\delta_i|}{4\pi^2}+O(1),\qquad
 \|w_i\|_\infty\le CW_i,\qquad w_i(\xi_j)=O(1)\quad(i\ne j).
\end{equation}
Let \(A\ge0\) and \(g_i\in L^\infty(\partial\Omega)\) satisfy $
 |g_i(x)|\le
 \frac{A\delta_i^{-3}}{(1+d(x,\xi_i)/\delta_i)^4}$ for a.e. $x\in\partial\Omega.$ It holds
\begin{equation}
 \|\mathcal P^{-1}g_i-\mu_iw_i\|_\infty\le CA,
 \qquad \mu_i=\int_{\partial\Omega}g_i\,d\sigma.
 \label{e4.1}
\end{equation}
If \(A=O(1)\) and $
 \delta_i^3g_i\bigl(\exp_{\xi_i}(\delta_i y)\bigr)\rightarrow0$ locally uniformly in $\mathbb R^3$ as $\varepsilon\to0,$ it follows
\begin{equation}
 \|\mathcal P^{-1}g_i-\mu_iw_i\|_\infty=o(1).
 \label{e4.2}
\end{equation}
\end{lemma}
\begin{proof}
By the definition of \(\mu_i\), we have
\[
\mathcal P^{-1}g_i(x)-\mu_iw_i(x)
=\int_{\partial\Omega}
\bigl[G_{\partial\Omega}(x,z)-w_i(x)\bigr]g_i(z)\,d\sigma(z).
\]
The local expansion \eqref{e2.6} and estimate \eqref{e3.7} imply $
\bigl|G_{\partial\Omega}(x,z)-w_i(x)\bigr|
\le C\left(1+\left|\log
\frac{d(x,z)}{\delta_i+d(x,\xi_i)}\right|\right),$ $x\ne z.$
For \(R\ge2\), we get
\begin{equation*}
\sup_{Y\in\mathbb R^3}\int_{|y|\le R}
\left(1+\left|\log\frac{|Y-y|}{1+|Y|}\right|\right)\,dy
\le C_R,\quad
\sup_{Y\in\mathbb R^3}\int_{|y|>R}
\frac{1+\left|\log\frac{|Y-y|}{1+|Y|}\right|}{(1+|y|)^4}\,dy
\le C\frac{1+\log R}{R}.
\end{equation*}
The change
of variables \(y=Y+(1+|Y|)s\) gives $
\frac{C}{1+|Y|}\int_{|s|<{1}/{4}}(1+|\log|s||)\,ds
\le\frac{C}{R}.$
For \(x\in B_{2r_0}(\xi_i)\) and \(z\in B_{r_0}(\xi_i)\), write
\(x=\exp_{\xi_i}(\delta_iY)\) and
\(z=\exp_{\xi_i}(\delta_i y)\). Then
\[
d(x,z)\asymp\delta_i|Y-y|,\qquad
\delta_i+d(x,\xi_i)=\delta_i(1+|Y|),\qquad
d\sigma(z)\le C\delta_i^3\,dy.
\]
If \(x\notin B_{2r_0}(\xi_i)\), then
\(|G_{\partial\Omega}(x,z)-w_i(x)|\le C\) for
\(z\in B_{r_0}(\xi_i)\).
On \(\partial\Omega\setminus B_{r_0}(\xi_i)\),
\(|g_i|\le CA\delta_i\); the logarithmic bound \eqref{e2.green-log-bound}
and \eqref{e3.7} therefore yield
\[
\sup_{x\in\partial\Omega}
\int_{\partial\Omega\setminus B_{r_0}(\xi_i)}
|G_{\partial\Omega}(x,z)-w_i(x)|\,|g_i(z)|\,d\sigma(z)
\le CA\delta_i(1+|\log\delta_i|).
\]
Splitting the integral over \(B_{r_0}(\xi_i)\) at \(|y|=R\),
 for \(2\le R<r_0/\delta_i\), we obtain
\begin{equation*}
\|\mathcal P^{-1}g_i-\mu_iw_i\|_\infty
\le C_R\operatorname*{ess\,sup}_{|y|\le R}
\delta_i^3\bigl|g_i(\exp_{\xi_i}(\delta_i y))\bigr|+CA\frac{1+\log R}{R}
+CA\delta_i(1+|\log\delta_i|).
\end{equation*}
Taking \(R=2\) proves \eqref{e4.1}.
If \(A=O(1)\) and
\(\delta_i^3g_i(\exp_{\xi_i}(\delta_i\,\cdot))\to0\) locally uniformly,
first letting \(\varepsilon\to0\) and then \(R\to\infty\)
proves \eqref{e4.2}.
\end{proof}

For $1\le i\le m$ and $0\le k\le3$, define $
Y_{i,k}=24\pi^2\mathcal P^{-1}(q_iZ_{i,k})$ and $
D_i=Y_{i,0}-\frac{w_i}{W_i}.$
Set $
R_\varepsilon=\delta_*^{-1/4}$ and $
\tau_\varepsilon=\delta_*^{1/4}\Lambda_\varepsilon^6.$

\begin{lemma}\label{t4.3}
For  sufficiently small $\varepsilon>0$ and $1\le i\le m$, it holds
\begin{equation}\label{e4.jacobi-norms}
\sum_{k=0}^3\bigl(\|Y_{i,k}\|_\infty+\|Y_{i,k}\|_{H^{\frac{3}{2}}}\bigr)
+\|D_i\|_\infty+\|D_i\|_{H^{\frac{3}{2}}}\le C
\end{equation}
and
\begin{equation}\label{e4.6}
\begin{aligned}
\|L_VY_{i,k}\|_*
+\Lambda_\varepsilon^{-\frac{1}{2}}\|L_VY_{i,k}\|_{H^{-\frac{3}{2}}}
&\le C\tau_\varepsilon,\qquad 1\le k\le3,\\
\|L_VD_i\|_*
+\Lambda_\varepsilon^{-\frac{1}{2}}\|L_VD_i\|_{H^{-\frac{3}{2}}}
&\le C\Lambda_\varepsilon^{-1}.
\end{aligned}
\end{equation}
\end{lemma}
\begin{proof}
For $0\le k\le3$, \eqref{e4.cutoff-moments}, \eqref{e4.1}, and
\eqref{e4.bubble-norms} give $\|Y_{i,k}\|_\infty\le C+C\delta_i^2W_i\le C$.
The coercivity estimate \eqref{e2.coercivity} yields
\[
\|Y_{i,k}\|_{H^{\frac{3}{2}}}^2
\le C\langle\mathcal PY_{i,k},Y_{i,k}\rangle\le C,\qquad
\left\|\frac{w_i}{W_i}\right\|_{H^{\frac{3}{2}}}^2
\le\frac{C}{W_i^2}\int_{\partial\Omega}q_iw_i\,d\sigma
\le\frac C{W_i}.
\]
Since $\|w_i/W_i\|_\infty\le C$, the $L^\infty$ and $H^{3/2}$ bounds follow from the definition of $D_i$.
Using $qz_0=\frac13\operatorname{div}(yq)$ and
$qz_k=-\frac1{3a}\partial_{y_k}q$ for $1\le k\le3$,
we integrate by parts in \eqref{e2.model-potential} to obtain
\begin{equation}\label{e4.model-jacobi-potential}
24\pi^2\int_{\mathbb R^3}\Gamma(y-s)q(s)z_k(s)\,ds
=\begin{cases}
1+z_0(y),&k=0,\\
z_k(y),&1\le k\le3.
\end{cases}
\end{equation}
For $|y|\le R_\varepsilon$, we apply \eqref{e2.6} and \eqref{e3.2}
to the representation \eqref{e2.green-representation} with $F=24\pi^2q_iZ_{i,k}$.
Since $|z_k|\le1$ and $\int_{\mathbb R^3}qz_k\,ds=0$,
\[
\begin{aligned}
&\left|Y_{i,k}(\exp_{\xi_i}(\delta_i y))
-24\pi^2\int_{\mathbb R^3}\Gamma(y-s)q(s)z_k(s)\,ds\right|\\
&\quad\le C\delta_i\int_{\mathbb R^3}(1+|y|+|s|)q(s)\,ds+C\delta_i^2\int_{\mathbb R^3}(1+|y|^2+|s|^2)
\bigl(\Lambda_\varepsilon+|\log|y-s||\bigr)q(s)\,ds\\
&\quad\le C\delta_i(1+|y|)
+C\delta_i^2(1+|y|^2)\Lambda_\varepsilon.
\end{aligned}
\]
By \eqref{e3.4} applied at $y$ and $0$, the same bound holds for
$|w_i(\exp_{\xi_i}(\delta_i y))-W_i-u_0(y)|$.
Since $\delta_i(1+R_\varepsilon)+\delta_i^2(1+R_\varepsilon^2)\Lambda_\varepsilon
\le C\tau_\varepsilon$, we obtain
\begin{equation}\label{e4.4}
\begin{aligned}
|w_i(\exp_{\xi_i}(\delta_i y))-W_i-u_0(y)|&\le C\tau_\varepsilon,\quad
|Y_{i,0}(\exp_{\xi_i}(\delta_i y))-1-z_0(y)|\le C\tau_\varepsilon,\\
|Y_{i,k}(\exp_{\xi_i}(\delta_i y))-z_k(y)|&\le C\tau_\varepsilon,
\qquad 1\le k\le3,
\end{aligned}
\end{equation}
for $|y|\le R_\varepsilon$.
For $x\in B_{r_0}(\xi_j)$, $j\ne i$, Taylor's formula for
$G_{\partial\Omega}(x,\cdot)$ and \eqref{e3.2} give
\[
|Y_{i,k}(x)|\le C\delta_i^2
+C\int_{\partial\Omega}d(z,\xi_i)q_i(z)\,d\sigma(z)
\le C\delta_i.
\]
Combining this estimate with \eqref{e3.6}, we obtain
\begin{equation}\label{e4.5}
Y_{i,k}(x)=O(\delta_*\Lambda_\varepsilon),\qquad
D_i(x)=-\frac{G_{\partial\Omega}(x,\xi_i)}{W_i}
+O(\delta_*\Lambda_\varepsilon),\qquad 0\le k\le3.
\end{equation}
For $1\le k\le3$, $L_VY_{i,k}
=24\pi^2q_i(Z_{i,k}-Y_{i,k})
-24\pi^2\sum_{j\ne i}q_jY_{i,k}
+24\pi^2\mathcal R_\varepsilon(V)Y_{i,k}.$
Moreover,
\[
\frac{q_j(x)}{\omega(x)}
\le C(1+d(x,\xi_j)/\delta_j)^{-2}
\le CR_\varepsilon^{-2}
\quad\text{for }d(x,\xi_j)\ge\delta_jR_\varepsilon.
\]
Using \eqref{e4.4}--\eqref{e4.5} on
$B_{\delta_jR_\varepsilon}(\xi_j)$, the preceding inequality
and \eqref{e3.14}, we obtain
\[
\|L_VY_{i,k}\|_*
\le C\bigl(\tau_\varepsilon+R_\varepsilon^{-2}
+\delta_*\Lambda_\varepsilon\bigr)
\le C\tau_\varepsilon.
\]
For $x=\exp_{\xi_i}(\delta_i y)$ with $|y|\le R_\varepsilon$,
using \eqref{e4.4} and \eqref{e3.14}, we have
\[
L_VD_i(x)=\frac{q_i(x)}{W_i}\bigl(24\pi^2u_0(y)-1\bigr)
+O\bigl(\tau_\varepsilon\omega(x)\bigr).
\]
Here $(1+|y|)^{-2}(1+|u_0(y)|)\le C$ and $W_i\asymp\Lambda_\varepsilon$.
On $B_{\delta_jR_\varepsilon}(\xi_j)$ with $j\ne i$,
\eqref{e4.5} implies $D_i=O(\Lambda_\varepsilon^{-1})$.
On $\partial\Omega\setminus\bigcup_{j=1}^m B_{\delta_jR_\varepsilon}(\xi_j)$, we have $q_j/\omega\le CR_\varepsilon^{-2}$ for every $j$. Thus,
\[
\|L_VD_i\|_*
\le C\bigl(\Lambda_\varepsilon^{-1}+\tau_\varepsilon
+R_\varepsilon^{-2}+\delta_*\Lambda_\varepsilon\bigr)
\le C\Lambda_\varepsilon^{-1}.
\]
Finally, \eqref{e4.weighted-embedding}, \eqref{e3.weight-bounds}, and the Cauchy--Schwarz inequality yield
\[
\int_{\partial\Omega}\omega|u|\,d\sigma
\le\|\omega\|_{L^1}^{\frac{1}{2}}
\left(\int_{\partial\Omega}\omega u^2\,d\sigma\right)^{\frac{1}{2}}
\le C\Lambda_\varepsilon^{\frac{1}{2}}\|u\|_{H^{\frac{3}{2}}},
\qquad u\in H^{\frac{3}{2}}(\partial\Omega).
\]
The $H^{-3/2}$ estimates in \eqref{e4.6} follow by duality.
\end{proof}
\begin{lemma}\label{t4.4}
Let $\varepsilon_n\to0^+$ and $\xi_n\in\mathcal K$.
Suppose $\phi_n\in H^{3/2}(\partial\Omega)\cap L^\infty(\partial\Omega)$
satisfies
\begin{equation}\label{pha}
L_{V_n}\phi_n=f_n+\sum_{j=1}^m\sum_{k=0}^3d_{j,k,n}q_{j,n}Z_{j,k,n},
\end{equation}
where $f_n\in L^\infty(\partial\Omega)$ and $d_{j,k,n}\in\mathbb R$.
Assume $\|\phi_n\|_\infty\le C$ and
$\|f_n\|_*+\sum_{j,k}|d_{j,k,n}|\to0$.
After passing to a subsequence,
\[
\widehat\phi_{i,n}(y):=\phi_n\bigl(\exp_{\xi_{i,n}}(\delta_{i,n}y)\bigr)
\rightarrow\widehat\phi_i(y)
\quad\text{in }C^{1,\alpha}(B_R(0))
\]
for every $1\le i\le m$, $R>0$, and $0<\alpha<1$, where
$\widehat\phi_i\in L^\infty(\mathbb R^3)$ satisfies
\begin{equation}\label{state}
2(-\Delta)^{\frac{3}{2}}\widehat\phi_i=24\pi^2q\widehat\phi_i
\quad\text{in }\mathcal S'(\mathbb R^3).
\end{equation}
For any fixed $i$, if $\int_{\partial\Omega}q_{i,n}Z_{i,k,n}\phi_n\,d\sigma=0$
for $0\le k\le3$ and every $n$, then $\widehat\phi_i\equiv0$.
\end{lemma}
\begin{proof}
By elliptic regularity, $\phi_n\in C^2(\partial\Omega)$.
Set $\omega_{j,n}(x)=\delta_{j,n}^{-3}(1+d(x,\xi_{j,n})/\delta_{j,n})^{-4}$,
$\omega_n=1+\sum_j\omega_{j,n}$, and
$g_{0,n}=\mathcal P\phi_n/\omega_n$, $g_{j,n}=\omega_{j,n}g_{0,n}$.
Then $\mathcal P\phi_n=g_{0,n}+\sum_jg_{j,n}$, and \eqref{pha} gives
\[
\|g_{0,n}\|_\infty\le C\|\phi_n\|_\infty+\|f_n\|_*+
C\sum_{j,k}|d_{j,k,n}|\le C,\qquad |g_{j,n}|\le C\omega_{j,n}.
\]
For $x=\exp_{\xi_{i,n}}(\delta_{i,n}y)$, $|y|\le R$,
using \eqref{e3.green-derivatives} for $s=1,2$ and the separation of the centers, we obtain
\[
\delta_{i,n}^s|\nabla_x^s\phi_n(x)|
\le C\int_{\mathbb R^3}|y-z|^{-s}(1+|z|)^{-4}\,dz
+C\delta_{i,n}^s\le C_R,\qquad s=1,2.
\]
Indeed, $\|g_{j,n}\|_{L^1}\le C$ and
$|g_{j,n}|\le C\delta_{j,n}$ in  $\partial\Omega\setminus B_{r_0}(\xi_{j,n})$.
These bounds also imply
$\nabla_x\mathcal P^{-1}(g_{0,n}+\sum_{j\ne i}g_{j,n})(x)=O(1)$.
Thus $\|\widehat\phi_{i,n}\|_{W^{2,\infty}(B_R(0))}\le C_R$,
and a diagonal argument yields convergence in $C^{1,\alpha}(B_R(0))$
for all $i$, $R>0$, and $0<\alpha<1$.

Fix $R_1>2R+2$ and take $n$ large enough that $\delta_{i,n}R_1<r_0$, it holds
\begin{equation*}
\widehat\phi_{i,n}(y)-\widehat\phi_{i,n}(0)
=\int_{B_{\delta_{i,n}R_1}(\xi_{i,n})}
\bigl[G_{\partial\Omega}(x,z)-G_{\partial\Omega}(\xi_{i,n},z)\bigr]
g_{i,n}(z)\,d\sigma(z)+O_R(R_1^{-2}+\delta_{i,n}).
\end{equation*}
Since $\omega_{i,n}/\omega_n\to1$ uniformly on $B_{\delta_{i,n}R_1}(\xi_{i,n})$,
\eqref{pha}, \eqref{e3.12}, and
$\|f_n\|_*+\sum_{j,k}|d_{j,k,n}|\to0$ imply
\[
\delta_{i,n}^3g_{i,n}\bigl(\exp_{\xi_{i,n}}(\delta_{i,n}z)\bigr)
\rightarrow24\pi^2q(z)\widehat\phi_i(z)
\quad\text{in }L^\infty(B_{R_1}(0)).
\]
By \eqref{e2.6}, for $|z|\le R_1$ and $z\ne0,y$, the rescaled Green kernel
difference converges to $\Gamma(y-z)-\Gamma(-z)$ and is bounded by
$C_{R_1}(1+|\log|y-z||+|\log|z||)$.
Dominated convergence on $|z|\le R_1$, followed by $R_1\to\infty$, yields
\begin{equation}\label{e4.7}
\widehat\phi_i(y)-\widehat\phi_i(0)
=24\pi^2\int_{\mathbb R^3}
[\Gamma(y-z)-\Gamma(-z)]q(z)\widehat\phi_i(z)\,dz.
\end{equation}
Applying $2(-\Delta)^{3/2}$ to \eqref{e4.7} gives \eqref{state}.
Under the orthogonality conditions, dominated convergence yields
$\int_{\mathbb R^3}qz_k\widehat\phi_i\,dy=0$ for $0\le k\le3$. Hence,
\cite[Theorem~1.1]{r9} and \eqref{e2.17} imply $\widehat\phi_i\equiv0$.
\end{proof}

\begin{proposition}\label{t4.5}
Let $f\in L^\infty(\partial\Omega)$.
For  sufficiently small $\varepsilon>0$ and every $\xi\in\mathcal K$, the problem
\begin{equation}\label{theeq}
 L_V\phi=f+\sum_{i=1}^m\sum_{k=0}^3d_{i,k}q_iZ_{i,k},
 \qquad \phi\in X_\perp
\end{equation}
has a unique weak solution $\phi$ and unique coefficients $d_{i,k}\in\mathbb R$.
Moreover,
\begin{equation}
 \|\phi\|_\infty+\|\phi\|_{H^{\frac{3}{2}}}
 +\sum_{i=1}^m\sum_{k=0}^3|d_{i,k}|\le C\|f\|_*,
 \label{e4.8}
\end{equation}
\end{proposition}
\begin{proof}
For fixed $\varepsilon>0$, the right-hand side of
$\mathcal P\phi=Q_V\phi+f+\sum_{i,k}d_{i,k}q_iZ_{i,k}$ belongs to $L^2(\partial\Omega)$.
Elliptic regularity therefore implies
$\phi\in H^3(\partial\Omega)\hookrightarrow C^0(\partial\Omega)$ for every weak solution.
By \eqref{e2.17}, \eqref{e4.4}, and \eqref{e4.jacobi-norms},
\[
\ell_{i,l}(Y_{j,k})=\frac14\delta_{ij}\delta_{lk}+o(1),\qquad
\ell_{i,l}(D_j)=\frac14\delta_{ij}\delta_{l0}+o(1),
\]
for $1\le i,j\le m$, $0\le l\le3$, and $1\le k\le3$.
Testing \eqref{theeq} against $Y_{j,k}$ and $D_j$, and using
\eqref{e4.6} and \eqref{e3.weight-bounds}, we obtain
\begin{equation}\label{e4.multiplier-bound}
\sum_{i,k}|d_{i,k}|
\le C\|f\|_*+C(\tau_\varepsilon+\Lambda_\varepsilon^{-1})\|\phi\|_\infty.
\end{equation}
We claim that $\|\phi\|_\infty\le C\|f\|_*$.
Otherwise, there exist $\varepsilon_n\to0$, $\xi_n\in\mathcal K$, and solutions
which, after normalization, satisfy
\[
\|\phi_n\|_\infty=1,\qquad
\|f_n\|_*+\sum_{i,k}|d_{i,k,n}|\longrightarrow0.
\]
By Lemma~\ref{t4.4}, we may pass to a subsequence such that
$\widehat\phi_{i,n}\to0$ locally uniformly in $\mathbb R^3$ for every $i$.
Since
\[
Q_{V_n}=24\pi^2\biggl(\sum_iq_{i,n}-\mathcal R_{\varepsilon_n}(V_n)\biggr),
\qquad
\frac{q_{i,n}(x)}{\omega_n(x)}
\le C\bigl(1+d(x,\xi_{i,n})/\delta_{i,n}\bigr)^{-2}.
\]
Using $
\mathcal P\phi_n
=Q_{V_n}\phi_n+f_n
+\sum_{i=1}^{m}\sum_{k=0}^{3}
d_{i,k,n}q_{i,n}Z_{i,k,n}$
 and \eqref{e3.14}, for every fixed $R\ge1$, we find
\[
\|\mathcal P\phi_n\|_*
\le C\max_{1\le i\le m}\|\widehat\phi_{i,n}\|_{L^\infty(B_R(0))}
+C(1+R)^{-2}+o(1).
\]
Letting $n\to\infty$ and then $R\to\infty$, we obtain
$\|\mathcal P\phi_n\|_*\to0$.
Set
\[
g_{0,n}=\frac{\mathcal P\phi_n}{\omega_n},\qquad
g_{i,n}(x)=\frac{\delta_{i,n}^{-3}g_{0,n}(x)}
{(1+d(x,\xi_{i,n})/\delta_{i,n})^4},\qquad
\mu_{i,n}=\int_{\partial\Omega}g_{i,n}\,d\sigma.
\]
Then $\mathcal P\phi_n=g_{0,n}+\sum_i g_{i,n}$,
$\|g_{0,n}\|_\infty=o(1)$, and $|\mu_{i,n}|\le C\|g_{0,n}\|_\infty=o(1)$.
Since $
\phi_n=\mathcal P^{-1}g_{0,n}
+\sum_{i=1}^{m}\mathcal P^{-1}g_{i,n},$
$\|\mathcal P^{-1}g_{0,n}\|_\infty
\le C\|g_{0,n}\|_\infty,$ \eqref{e4.1}  yields
\[
\left\|\phi_n-\sum_i\mu_{i,n}w_{i,n}\right\|_\infty
\le C\|g_{0,n}\|_\infty=o(1).
\]
Evaluating at $\xi_{j,n}$ and using
$\phi_n(\xi_{j,n})\to0$ and $w_{i,n}(\xi_{j,n})=O(1)$ for $i\ne j$, we find
$\mu_{j,n}W_{j,n}=o(1)$. Lemma~\ref{t4.2} now yields $
\|\phi_n\|_\infty
\le C\sum_j|\mu_{j,n}|W_{j,n}+o(1)=o(1),$
contrary to $\|\phi_n\|_\infty=1$.
Finally, testing \eqref{theeq} against $\phi\in X_\perp$ and using the coercivity
of $\mathcal P$ gives
\[
\|\phi\|_{H^{\frac{3}{2}}}^2
\le C\|\phi\|_\infty^2+C\|f\|_*\|\phi\|_\infty
\le C\|f\|_*^2.
\]
Together with \eqref{e4.multiplier-bound}, this proves \eqref{e4.8}.

For fixed $\varepsilon$, the coercivity of $\mathcal P$ and the compact embedding
$H^{3/2}(\partial\Omega)\hookrightarrow L^2(\partial\Omega)$ imply that
the operator associated with $\mathcal B_V$ on $X_\perp$ is Fredholm of index zero.
Since $X_\perp=\bigcap_{i,k}\ker\ell_{i,k}$, every solution of the variational
equation on $X_\perp$ satisfies the stated equation for some coefficients $d_{i,k}$.
Estimate~\eqref{e4.8} and the Fredholm alternative establish the existence and  the uniqueness.
\end{proof}

\begin{lemma}\label{t4.6}
For  sufficiently small $\varepsilon>0$ and every $\xi\in\mathcal K$,
there exist functions
$\widehat Z_{i,0}\in H^{3/2}(\partial\Omega)\cap L^\infty(\partial\Omega)$,
$1\le i\le m$, such that
\begin{equation}\label{e4.dilation-orthogonality}
\begin{gathered}
\ell_{j,k}(\widehat Z_{i,0})=0,
\qquad 1\le i,j\le m,\quad 1\le k\le3,\\
\mathcal B_V(\widehat Z_{i,0},\psi)=0,
\qquad 1\le i\le m,\quad \text{for all }\psi\in X_\perp.
\end{gathered}
\end{equation}
Moreover,
\begin{equation}\label{e4.dilation-approximation}
\begin{aligned}
\|\widehat Z_{i,0}-D_i\|_\infty
&+\|\widehat Z_{i,0}-D_i\|_{H^{\frac{3}{2}}}
\le C\Lambda_\varepsilon^{-1},
&&1\le i\le m,\\
\ell_{j,0}(\widehat Z_{i,0})
&=\frac14\delta_{ij}+O(\Lambda_\varepsilon^{-1}),
&&1\le i,j\le m,
\end{aligned}
\end{equation}
and
\begin{equation}\label{e4.9}
\mathcal B_V(\widehat Z_{i,0},\widehat Z_{j,0})
=\frac{\delta_{ij}}{W_i}+O(\Lambda_\varepsilon^{-2}),
\qquad 1\le i,j\le m.
\end{equation}
Then the matrix
$\bigl(\mathcal B_V(\widehat Z_{i,0},\widehat Z_{j,0})\bigr)_{i,j=1}^m$
is positive definite and its inverse has norm $O(\Lambda_\varepsilon)$.
\end{lemma}
\begin{proof}
By \eqref{e2.16}--\eqref{e2.17} and \eqref{e4.4}--\eqref{e4.5}, we obtain
\begin{equation}\label{e4.translation-moments}
\ell_{j,k}(Y_{i,l})
=\frac14\delta_{ij}\delta_{kl}+O(\tau_\varepsilon),\qquad
\ell_{j,k}(D_i)=O(\tau_\varepsilon),\qquad
\ell_{j,0}(D_i)=\frac14\delta_{ij}+O(\Lambda_\varepsilon^{-1}),
\end{equation}
for $1\le i,j\le m$ and $1\le k,l\le3$.
Hence, the translation moment matrix is invertible for $\varepsilon$ small.
For each $i$, choose constants $a_{j,k}^{(i)}=O(\tau_\varepsilon)$ so that
\[
\widetilde D_i=D_i+\sum_{j=1}^m\sum_{k=1}^3a_{j,k}^{(i)}Y_{j,k},
\qquad
\ell_{p,l}(\widetilde D_i)=0,\quad 1\le p\le m,\quad 1\le l\le3.
\]
By \eqref{e4.6}, $\|L_V\widetilde D_i\|_*=O(\Lambda_\varepsilon^{-1})$.
Applying Proposition~\ref{t4.5} with $f=-L_V\widetilde D_i$, we obtain
$\psi_i\in X_\perp$ with
$\|\psi_i\|_\infty+\|\psi_i\|_{H^{3/2}}=O(\Lambda_\varepsilon^{-1})$.
Set $\widehat Z_{i,0}=\widetilde D_i+\psi_i$.
In view of $\psi_i\in X_\perp$, we have
$\ell_{p,l}(\widehat Z_{i,0})=0$ for $1\le p\le m$ and $1\le l\le3$.
By the projected equation,
$\mathcal B_V(\widehat Z_{i,0},\psi)=0$ for every $\psi\in X_\perp$.
The norm estimate in \eqref{e4.dilation-approximation} follows from \eqref{e4.jacobi-norms}
and $\tau_\varepsilon=o(\Lambda_\varepsilon^{-1})$.
Since $\ell_{j,0}(\psi_i)=0$ and $|\ell_{j,0}(Y_{p,l})|\le C$,
we have $\ell_{j,0}(\widehat Z_{i,0})=\ell_{j,0}(D_i)+O(\tau_\varepsilon)
=\frac14\delta_{ij}+O(\Lambda_\varepsilon^{-1})$.

We next prove \eqref{e4.9}. Using \eqref{e4.4}--\eqref{e4.5}, we deduce
\begin{equation}\label{e4.dilation-expansions}
\begin{aligned}
D_i(\exp_{\xi_i}(\delta_i y))
&=z_0(y)-\frac{u_0(y)}{W_i}+O(\tau_\varepsilon),
&&|y|\le R_\varepsilon,\\
D_i(x)&=-\frac{G_{\partial\Omega}(\xi_j,\xi_i)}{W_i}
+O(\tau_\varepsilon),
&&x\in B_{\delta_jR_\varepsilon}(\xi_j),\quad j\ne i.
\end{aligned}
\end{equation}
Note that $\tau_\varepsilon=o(\Lambda_\varepsilon^{-2})$ and
\[
\int_{|y|>R_\varepsilon}q(y)(1+|u_0(y)|)\,dy
\le CR_\varepsilon^{-3}(1+\log R_\varepsilon)
=o(\Lambda_\varepsilon^{-2}).
\]
By \eqref{e2.16}, \eqref{e4.dilation-expansions}, and
the uniform bounds \eqref{e4.jacobi-norms}, we obtain
\[
\int_{\partial\Omega}q_jD_j\,d\sigma
=-\frac1{W_j}\int_{\mathbb R^3}qu_0\,dy+O(\tau_\varepsilon)
=O(\Lambda_\varepsilon^{-1}),\qquad
\int_{\partial\Omega}q_kD_iD_j\,d\sigma
=O(\Lambda_\varepsilon^{-2})\quad(k\ne i).
\]
For $k=j\ne i$, the second estimate follows from the first
and the expansion of $D_i$ near $\xi_j$ in \eqref{e4.dilation-expansions}.
For $k\ne i,j$, both $D_i$ and $D_j$ are $O(\Lambda_\varepsilon^{-1})$
on $B_{\delta_kR_\varepsilon}(\xi_k)$.
From the definition of $D_i$ and \eqref{e3.14}, we have
\[
\mathcal PD_i=24\pi^2q_iZ_{i,0}-\frac{q_i}{W_i},\qquad
Q_V=24\pi^2\sum_{k=1}^m q_k
+O(\delta_*\Lambda_\varepsilon\omega).
\]
Hence,
\[
\begin{aligned}
\mathcal B_V(D_i,D_j)
&=\int_{\partial\Omega}q_i
\left[24\pi^2(Z_{i,0}-D_i)-\frac1{W_i}\right]D_j\,d\sigma
-24\pi^2\sum_{k\ne i}\int_{\partial\Omega}q_kD_iD_j\,d\sigma
+O(\tau_\varepsilon)\\
&=\frac{\delta_{ij}}{W_i}
\int_{\mathbb R^3}qz_0(24\pi^2u_0-1)\,dy
+O(\Lambda_\varepsilon^{-2}),
\end{aligned}
\]
where the terms containing $u_0^2/W_i^2$ are
$O(\Lambda_\varepsilon^{-2})$ due to
$\int_{\mathbb R^3}qu_0^2\,dy<\infty$.
Integration by parts, combining $\operatorname{div}(yq)=3qz_0$
and \eqref{e2.16}, we infer
\[
\int_{\mathbb R^3}qz_0(24\pi^2u_0-1)\,dy
=-8\pi^2\int_{\mathbb R^3}q\,y\cdot\nabla u_0\,dy
=\int_{\mathbb R^3}q(1-z_0)\,dy=1.
\]
Moreover, it follows from \eqref{e4.6}, \eqref{e4.dilation-approximation},
and $\int_{\partial\Omega}(\omega+Q_V)\,d\sigma\le C$ that
\[
\begin{aligned}
|\mathcal B_V(D_i,\widehat Z_{j,0}-D_j)|
&\le C\|L_VD_i\|_*\|\widehat Z_{j,0}-D_j\|_\infty
\le C\Lambda_\varepsilon^{-2},\\
|\mathcal B_V(\widehat Z_{i,0}-D_i,\widehat Z_{j,0}-D_j)|
&\le C\|\widehat Z_{i,0}-D_i\|_{H^{\frac{3}{2}}}
\|\widehat Z_{j,0}-D_j\|_{H^{\frac{3}{2}}}+C\|\widehat Z_{i,0}-D_i\|_\infty
\|\widehat Z_{j,0}-D_j\|_\infty
\le C\Lambda_\varepsilon^{-2}.
\end{aligned}
\]
Expanding the bilinear form and using its symmetry,
\eqref{e4.9} holds. Since $W_i\asymp\Lambda_\varepsilon$, the matrix is
positive definite and its inverse has norm $O(\Lambda_\varepsilon)$.
\end{proof}

\begin{proposition}\label{t4.7}
For  sufficiently small $\varepsilon>0$, every $\xi\in\mathcal K$,
and every $f\in L^\infty(\partial\Omega)$, the problem
\begin{equation}\label{e4.projected-problem}
 L_V\phi=f+\sum_{i=1}^m\sum_{k=1}^3c_{i,k}q_iZ_{i,k},
\quad
 \ell_{i,k}(\phi)=0,\quad
 1\le i\le m,\quad 1\le k\le3,
\end{equation}
has a unique weak solution
$\phi\in H^{3/2}(\partial\Omega)\cap L^\infty(\partial\Omega)$
and uniquely determined coefficients $c_{i,k}\in\mathbb R$. Moreover,
\begin{equation}
 \|\phi\|_\infty+\|\phi\|_{H^{\frac{3}{2}}}
 +\sum_{i=1}^m\sum_{k=1}^3|c_{i,k}|
 \le C\Lambda_\varepsilon\|f\|_*.
 \label{e4.10}
\end{equation}
\end{proposition}
\begin{proof}
By \eqref{e4.dilation-orthogonality} and \eqref{e4.dilation-approximation}, every $\phi\in H^{3/2}(\partial\Omega)$ satisfying
$\ell_{i,k}(\phi)=0$ for $1\le i\le m$ and $1\le k\le3$
has a unique decomposition
\[
 \phi=\phi_\perp+\sum_{i=1}^m b_i\widehat Z_{i,0},
 \qquad \phi_\perp\in X_\perp,\quad b_i\in\mathbb R.
\]
Let $\phi_\perp$ be the solution from Proposition~\ref{t4.5}, and choose $b_i$ so that
\[
 \sum_{i=1}^m \mathcal B_V(\widehat Z_{i,0},\widehat Z_{j,0})b_i
 =\int_{\partial\Omega}f\widehat Z_{j,0}\,d\sigma,
 \qquad 1\le j\le m.
\]
Estimates \eqref{e4.jacobi-norms}, \eqref{e4.dilation-approximation}, and \eqref{e4.9} imply
\[
 \sum_{i=1}^m|b_i|\le C\Lambda_\varepsilon\|f\|_*,
 \qquad
 \|\phi\|_\infty+\|\phi\|_{H^{\frac{3}{2}}}
 \le C\Lambda_\varepsilon\|f\|_*.
\]
The choice of $\phi_\perp$ and $b_i$, together with the
$\mathcal B_V$-orthogonality of $\widehat Z_{i,0}$ to $X_\perp$, gives
\[
 \mathcal B_V(\phi,\psi)
 =\int_{\partial\Omega}f\psi\,d\sigma
\]
for every $\psi\in H^{3/2}(\partial\Omega)$ with zero translation moments.
Hence $\phi$ satisfies the stated equation for uniquely determined
coefficients $c_{i,k}$. Testing against $Y_{j,k}$, $1\le k\le3$, using \eqref{e4.6} and the translation moment estimate \eqref{e4.translation-moments}, we obtain
\[
 \sum_{i=1}^m\sum_{k=1}^3|c_{i,k}|
 \le C\bigl(\|f\|_*+\|\phi\|_\infty\bigr)
 \le C\Lambda_\varepsilon\|f\|_*.
\]
This yields \eqref{e4.10}. If $f=0$, by Proposition~\ref{t4.5} and the invertibility of the matrix
in \eqref{e4.9}, we get $\phi_\perp=0$ and $b_i=0$ for all $i$.
Thus $\phi=0$, and uniqueness is derived.
\end{proof}

\subsection{The nonlinear projected problem}\label{s5}\label{s5.1}

For $\phi\in L^\infty(\partial\Omega)$, define
\[
 \mathcal N_V(\phi)=\varepsilon he^{24\pi^2V}
 \bigl(e^{24\pi^2\phi}-1-24\pi^2\phi\bigr).
\]
We seek
$\phi\in H^{3/2}(\partial\Omega)\cap L^\infty(\partial\Omega)$
and  $c_{i,k}\in\mathbb R$ satisfying
\begin{equation}
 L_V\phi=-\mathcal R_\varepsilon(V)+\mathcal N_V(\phi)
 +\sum_{i=1}^m\sum_{k=1}^3c_{i,k}q_iZ_{i,k},
 \label{e5.1}
\end{equation}
\begin{equation}
 \ell_{i,k}(\phi)=0,
 \qquad 1\le i\le m,\quad1\le k\le3.
 \label{e5.2}
\end{equation}
If $c_{i,k}=0$ for all $1\le i\le m$ and $1\le k\le3$,
then $V+\phi$ solves \eqref{e1.4}.

\begin{proposition}\label{t5.2}
There exists $A>0$ such that, for  sufficiently small $\varepsilon>0$
and every $\xi\in\mathcal K$, problem \eqref{e5.1}--\eqref{e5.2}
has a unique solution with
$\|\phi\|_\infty\le A\delta_*\Lambda_\varepsilon^2$.
Moreover,
\begin{equation}
 \|\phi\|_\infty+\|\phi\|_{H^{\frac{3}{2}}}
 +\sum_{i=1}^m\sum_{k=1}^3|c_{i,k}|
 \le C\delta_*\Lambda_\varepsilon^2.
 \label{e5.3}
\end{equation}
\end{proposition}
\begin{proof}
For $\phi,\psi\in L^\infty(\partial\Omega)$ with
$\|\phi\|_\infty,\|\psi\|_\infty\le1$,
Taylor's formula and \eqref{e4.potential-bound} yield
\begin{equation}\label{e5.nonlinear-estimates}
 \|\mathcal N_V(\phi)\|_*\le C\|\phi\|_\infty^2,\quad
 \|\mathcal N_V(\phi)-\mathcal N_V(\psi)\|_*
 \le C(\|\phi\|_\infty+\|\psi\|_\infty)\|\phi-\psi\|_\infty.
\end{equation}
Let $\mathcal A(\phi)$ be the unique solution of \eqref{e4.projected-problem},
given by Proposition~\ref{t4.7} with
$f=-\mathcal R_\varepsilon(V)+\mathcal N_V(\phi)$.
Then \eqref{e5.1}--\eqref{e5.2} is equivalent to
$\phi=\mathcal A(\phi)$.
Estimates \eqref{e3.14}, \eqref{e4.10}, and \eqref{e5.nonlinear-estimates} imply
\[
 \|\mathcal A(\phi)\|_\infty
 \le C\Lambda_\varepsilon
 \bigl(\delta_*\Lambda_\varepsilon+\|\phi\|_\infty^2\bigr).
\]
Choose $A>2C$. Since $\delta_*\Lambda_\varepsilon^3\to0$
uniformly for $\xi\in\mathcal K$, for sufficiently small $\varepsilon$,
we have $A\delta_*\Lambda_\varepsilon^2\le1$ and
\begin{equation*}
 \|\mathcal A(\phi)\|_\infty
 \le \bigl(C+CA^2\delta_*\Lambda_\varepsilon^3\bigr)
 \delta_*\Lambda_\varepsilon^2
 \le A\delta_*\Lambda_\varepsilon^2,\quad
 \|\mathcal A(\phi)-\mathcal A(\psi)\|_\infty
 \le CA\delta_*\Lambda_\varepsilon^3\|\phi-\psi\|_\infty
 \le\tfrac12\|\phi-\psi\|_\infty
\end{equation*}
whenever $\|\phi\|_\infty,\|\psi\|_\infty
\le A\delta_*\Lambda_\varepsilon^2$.
Thus $\mathcal A$ is a contraction of this closed ball in
$L^\infty(\partial\Omega)$ into itself.
The contraction mapping theorem ensures a unique fixed point,
and Proposition~\ref{t4.7} determines the coefficients $c_{i,k}$ uniquely.
Finally, \eqref{e4.10} yields
\begin{equation*}
 \|\phi\|_\infty+\|\phi\|_{H^{\frac{3}{2}}}
 +\sum_{i=1}^m\sum_{k=1}^3|c_{i,k}|
 \le C\Lambda_\varepsilon
 \bigl(\delta_*\Lambda_\varepsilon+\|\phi\|_\infty^2\bigr)\le C\delta_*\Lambda_\varepsilon^2.
\end{equation*}
\end{proof}

Denote the correction  in Proposition~\ref{t5.2} by $\phi_{\varepsilon,\xi}$ and set
\[
v_{\varepsilon,\xi}=V_\xi+\phi_{\varepsilon,\xi},\qquad
\beta=24\pi^2,\qquad Q_{i,k}=q_iZ_{i,k}
\quad(1\le i\le m,\ 1\le k\le3).
\]
Then
\begin{equation}\label{e5.projected-residual}
\mathcal Pv_{\varepsilon,\xi}-\varepsilon he^{\beta v_{\varepsilon,\xi}}
=\sum_{i=1}^m\sum_{k=1}^3c_{i,k}Q_{i,k}.
\end{equation}
On a local chart of $\mathcal M_m$, let $\partial_p$ denote coordinate
derivatives and write $\partial_{p\kappa}=\partial_p\partial_\kappa$,
$1\le p,\kappa\le3m$. Using the frames fixed above, define
\[
g_{i,k}=d_i\mathcal F_m[e_{i,k}],\qquad
\mathbb A_{(i,k),p}=\langle e_{i,k},\partial_p\xi_i\rangle.
\]
Thus $
\partial_p\mathcal F_m=\sum_{i=1}^m\sum_{k=1}^3
g_{i,k}\mathbb A_{(i,k),p},$ $\det\mathbb A\ne0.$

Fix $\xi^0\in\mathcal K$. For $\xi$ near $\xi^0$, choose
diffeomorphisms $\Phi_\xi:\partial\Omega\to\partial\Omega$ satisfying
\begin{equation}\label{e5.transport-map}
\Phi_\xi\bigl(\exp_{\xi_i^0}(\delta_i(\varepsilon,\xi^0)y)\bigr)
=\exp_{\xi_i}(\delta_i y)
\end{equation}
on $B_{3r_0}(\xi_i^0)$, with $\Phi_\xi=\mathrm{id}$ in
$\partial\Omega\setminus\bigcup_{i=1}^m B_{4r_0}(\xi_i^0)$.
In the chosen frames, the local map is
\[
x\longmapsto\exp_{\xi_i}\left(
\frac{\delta_i(\varepsilon,\xi)}{\delta_i(\varepsilon,\xi^0)}
\exp_{\xi_i^0}^{-1}x\right).
\]
Extend its coordinate displacement from the identity by a fixed
smooth cutoff. By \eqref{e3.10}, the ratio of scales is independent
of $\varepsilon$ and has bounded parameter derivatives up to order two.
For $\xi$ close to $\xi^0$, the extension is $C^1$-close to the identity;
the maps and their inverses have uniform spatial $C^2$ bounds,
also after two parameter derivatives.
Thus \eqref{e2.transported-green}--\eqref{e2.transported-green-gradient}
apply. Define $J_\xi$ by $\Phi_\xi^*d\sigma=J_\xi\,d\sigma$.
Material derivatives act on $f\circ\Phi_\xi$ and
$(f\circ\Phi_\xi)J_\xi$ for functions and densities, respectively,
with the reference point fixed, hence with $y$ fixed in the rescaled
coordinates.

\begin{lemma}\label{t5.3}
Let $E_i$ be the remainder in \eqref{e3.15}.
For $1\le i\le m$, $1\le k\le3$, and $1\le p\le3m$, we have
\begin{equation}\label{e5.residual-derivative}
|\partial_pE_i(y)|
\le C\delta_i^2\Lambda_\varepsilon^2(1+|y|)^{-4},
\qquad |y|\le2\delta_i^{-\frac{1}{2}},
\end{equation}
\begin{equation}\label{e5.density-derivatives}
\left|\partial_p[(\mathcal R_\varepsilon(V)\circ\Phi_\xi)J_\xi]\right|
\le C\delta_*\Lambda_\varepsilon^2(\omega\circ\Phi_\xi)J_\xi,\quad
\left|\partial_p[((\varepsilon he^{\beta V})\circ\Phi_\xi)J_\xi]\right|
\le C(\omega\circ\Phi_\xi)J_\xi.
\end{equation}
The functions $Y_{i,k}$ satisfy
\begin{equation}\label{e5.jacobi-residual-moments}
\|L_VY_{i,k}\|_{L^1}
+\left\|\partial_p[((L_VY_{i,k})\circ\Phi_\xi)J_\xi]\right\|_{L^1}
\le C\delta_i\Lambda_\varepsilon^3,\quad
\left\|\int_{\partial\Omega}\mathcal R_\varepsilon(V)Y_{i,k}\,d\sigma
+\frac{2\pi^4}{a^4}\delta_i g_{i,k}\right\|_{C^1}
\le C\delta_i^2\Lambda_\varepsilon^4.
\end{equation}
For $\alpha=(i,k)$ and $\gamma=(j,l)$, $1\le i,j\le m$,
$1\le k,l\le3$, note $Q_\alpha=Q_{i,k}$,
$Y_\alpha=Y_{i,k}$, $\ell_\alpha=\ell_{i,k}$, and $c_\alpha=c_{i,k}$. Define $\mathbb G_{\alpha\gamma}
=\int_{\partial\Omega}Q_\gamma Y_\alpha\,d\sigma$ and $\mathbb G=(\mathbb G_{\alpha\gamma})\in\mathbb R^{3m\times3m}$,
then
\begin{equation}\label{e5.gram-matrix}
\left\|\mathbb G-\frac14I\right\|_{C^1}
\le C\delta_*\Lambda_\varepsilon^3.
\end{equation}
All estimates are uniform for $\xi\in\mathcal K$.
\end{lemma}
\begin{proof}
Write $j_\xi(z)=\sqrt{\det g_\xi(z)}$, $r=|y|$, and
$\lambda_{i,p}=\partial_p\log\delta_i$.
By \eqref{e3.10}, $|\lambda_{i,p}|+|\partial_\kappa\lambda_{i,p}|\le C$.
In normal coordinates, we have
\begin{equation}\label{e5.jacobian-derivatives}
\|j_{\xi_i}(\delta_i y)-1\|_{C^2}
\le C\delta_i^2r^2,\qquad r<\frac{2r_0}{\delta_i},
\end{equation}
and
\begin{equation}\label{e5.source-densities}
q_i\,d\sigma=b_{\delta_i,\xi_i}\chi(\delta_i r)j_{\xi_i}(\delta_i y)q(y)\,dy,\quad
Q_{i,k}\,d\sigma=b_{\delta_i,\xi_i}\chi(\delta_i r)^2j_{\xi_i}(\delta_i y)q(y)z_k(y)\,dy.
\end{equation}
The first two parameter derivatives of the cutoff are bounded and
supported at $r\asymp\delta_i^{-1}$, where the $q$-tail is $O(\delta_i^3)$.
Differentiating \eqref{e3.1} therefore yields
\begin{equation}\label{e5.normalization-derivatives}
\|b_{\delta_i,\xi_i}-1\|_{C^2}\le C\delta_i^2.
\end{equation}
Since $\int_{\mathbb R^3}s\chi(\delta_i|s|)q(s)\,ds=0$, the first moment in
\eqref{e3.2} has $C^2$ norm at most
\[
C\delta_i^3\int_{|s|<2r_0/\delta_i}|s|^3q(s)\,ds
\le C\delta_i^2\int_{\mathbb R^3}|s|^2q(s)\,ds\le C\delta_i^2.
\]
The second moment has the same bound by \eqref{e5.source-densities}.
For $r\le2\delta_i^{-1/2}$, \eqref{e2.kernel-dilation} gives
\begin{equation*}
\int_{|s|<2r_0/\delta_i}\|E_{\xi_i}(\delta_i y,\delta_i s)\|_{C^2}q(s)\,ds
\le C\delta_i^2\int_{\mathbb R^3}
(r+|s|)^2(1+\Lambda_\varepsilon+|\log|y-s||)q(s)\,ds
\le C\delta_i^2(1+r^2)\Lambda_\varepsilon.
\end{equation*}
By \eqref{n3}--\eqref{n4} and
\eqref{e5.jacobian-derivatives}--\eqref{e5.normalization-derivatives},
the errors in \eqref{defW} arising from the cutoffs, the normalization
factor $b_{\delta_i,\xi_i}$, and the surface density are bounded by
$C\delta_i^2(1+r^2)\Lambda_\varepsilon$ in $C^2$.
Since $H\in C^4$, Taylor expansion of $H$ and of
$W_{\delta_j,\xi_j}(\exp_{\xi_i}(\delta_i y))$, $j\ne i$,
together with \eqref{e3.2} and its first two parameter derivatives, yields the bound $
C\left[\delta_i^2(1+r^2)+\sum_{j\ne i}\delta_j^2\right]$
for the $C^2$ norms of the corresponding remainders.
Thus
\begin{equation}\label{e5.inner-remainder-derivatives}
\|E_{i,0}(y)\|_{C^2}\le C\delta_i^2(1+r^2)\Lambda_\varepsilon^2.
\end{equation}
Expanding \eqref{e3.12} and using
$\delta_i(1+r)+\delta_i^2(1+r^2)\Lambda_\varepsilon^2=o(1)$, we obtain
\[
\begin{aligned}
\beta^{-1}\delta_i^3Q_V(\exp_{\xi_i}(\delta_i y))
&=q(y)\left\{1+\beta\delta_i
\left[B_{\xi_i}(y)+\sum_lg_{i,l}y_l\right]
+O_{C^1}\bigl(\delta_i^2(1+r^2)\Lambda_\varepsilon^2\bigr)\right\},\\
\delta_i^3\mathcal PV(\exp_{\xi_i}(\delta_i y))&=b_{\delta_i,\xi_i}q(y).
\end{aligned}
\]
Subtracting and using $q(y)(1+r)^2\le C(1+r)^{-4}$ proves
\eqref{e5.residual-derivative}; the density transformation
\eqref{e5.source-densities} implies \eqref{e5.density-derivatives} for $d(x,\xi_i)\le 2\delta_i^{1/2}$, $1\le i\le m$.

For the remaining region, use
\[
(V\circ\Phi_\xi)(x)=\sum_{j=1}^m\int_{\partial\Omega}
G_{\partial\Omega}(\Phi_\xi x,\Phi_\xi z)
(q_j\circ\Phi_\xi)(z)J_\xi(z)\,d\sigma(z).
\]
Equations \eqref{e5.source-densities}--\eqref{e5.normalization-derivatives}
and $q(y)\le C(1+|y|)^{-6}$ imply
\[
\left|\partial_\xi^\alpha[(q_j\circ\Phi_\xi)J_\xi]\right|
\le C(\omega\circ\Phi_\xi)J_\xi,\qquad |\alpha|\le2.
\]
The cutoff derivatives are supported at $|y|\asymp\delta_j^{-1}$
and satisfy the same bound. Hence
\eqref{e2.transported-green} and the weighted logarithmic estimate
in \eqref{e2.14} yield
\begin{equation}\label{e5.material-V-bound}
\sum_{|\alpha|\le2}\|\partial_\xi^\alpha(V\circ\Phi_\xi)\|_\infty
\le C\Lambda_\varepsilon.
\end{equation}
For $\sqrt{\delta_i}\le r_x=d(x,\xi_i)\le r_0$,
\begin{equation}\label{e5.annular-density}
q_i(x)+\beta^{-1}Q_V(x)\le C\delta_i^3r_x^{-6},
\left|\partial_p[(q_i\circ\Phi_\xi)J_\xi](\Phi_\xi^{-1}x)\right|
+\left|\partial_p[((\beta^{-1}Q_V)\circ\Phi_\xi)J_\xi]
(\Phi_\xi^{-1}x)\right|
\le C\Lambda_\varepsilon\delta_i^3r_x^{-6}.
\end{equation}
Here $\Lambda_\varepsilon\delta_i^3r_x^{-6}
\le C\delta_i\Lambda_\varepsilon\omega(x)$ and $\delta_i\Lambda_\varepsilon=o(1)$.
For $x\in\partial\Omega\setminus\bigcup_{i=1}^m B_{r_0}(\xi_i)$,
we have
\begin{equation*}
\sum_{i=1}^m q_i(x)+\beta^{-1}Q_V(x)+\sum_{i=1}^m
\left|\partial_p[(q_i\circ\Phi_\xi)J_\xi]
(\Phi_\xi^{-1}x)\right|
+\left|\partial_p[((\beta^{-1}Q_V)\circ\Phi_\xi)J_\xi]
(\Phi_\xi^{-1}x)\right|
\le C\delta_*^3\Lambda_\varepsilon.
\end{equation*}
Since $\mathcal R_\varepsilon(V)=\sum_iq_i-\beta^{-1}Q_V$
and $\omega\ge1$, this proves \eqref{e5.density-derivatives}
on $\partial\Omega\setminus\bigcup_iB_{r_0}(\xi_i)$.
The bounds already established for $d(x,\xi_i)\le2\sqrt{\delta_i}$
and the annular estimate \eqref{e5.annular-density}
complete the proof of \eqref{e5.density-derivatives}
on $\partial\Omega$.

We next estimate $Y_{i,k}$, $1\le k\le3$.
Since $\int qz_k=0$ and $K_{-4,\xi_i}$ is Lipschitz, we obtain
\[
\left|\int_{\mathbb R^3}[K_{-4,\xi_i}(y-s)-K_{-4,\xi_i}(y)]q(s)z_k(s)\,ds\right|
\le C\int_{\mathbb R^3}|s|q(s)|z_k(s)|\,ds\le C.
\]
Using \eqref{e4.model-jacobi-potential} and subtracting the smooth
kernel value at $s=0$, \eqref{e2.kernel-dilation} and the density
estimates \eqref{e5.jacobian-derivatives}--\eqref{e5.normalization-derivatives} yield
\begin{equation}\label{e5.jacobi-core}
\|Y_{i,k}(\exp_{\xi_i}(\delta_i y))-z_k(y)\|_{C^1}
\le C\delta_i\Lambda_\varepsilon^2
+C\delta_i^2(1+r^2)\Lambda_\varepsilon^2,
\qquad r\le2\delta_i^{-\frac{1}{2}}.
\end{equation}
For $x=\exp_{\xi_i}(\delta_i y)$ with $\sqrt{\delta_i}\le r_x\le r_0$,
split the Green representation at $|s|=r_x/(2\delta_i)$.
For $|s|\le r_x/(2\delta_i)$,
$d(x,\exp_{\xi_i}(t\delta_i s))\ge r_x/2$, $0\le t\le1$.
The path $t\mapsto\exp_{\xi_i^0}(t\delta_i(\varepsilon,\xi^0)s)$
is fixed in the reference coordinates and has speed at most
$C\delta_i|s|$. Integrating
\eqref{e2.transported-green-gradient} along this path yields
\[
\|G_{\partial\Omega}(x,\exp_{\xi_i}(\delta_i s))
-G_{\partial\Omega}(x,\xi_i)\|_{C^1}
\le C\delta_i\frac{|s|}{r_x},\qquad |s|\le \frac{r_x}{2\delta_i}.
\]
Using $\int_{\mathbb R^3}qz_k\,ds=0$,
$q(s)|z_k(s)|\le C(1+|s|)^{-7}$, and \eqref{e5.source-densities}, we obtain
\begin{equation}\label{e5.jacobi-annulus}
|Y_{i,k}(x)|+|Z_{i,k}(x)|
+|\partial_p(Y_{i,k}\circ\Phi_\xi)(\Phi_\xi^{-1}x)|
\le C\Lambda_\varepsilon^2\frac{\delta_i}{r_x+\delta_i}.
\end{equation}
On $\partial\Omega\setminus B_{r_0}(\xi_i)$, the same terms are bounded by
$C\delta_i\Lambda_\varepsilon^2$.
Since $\mathcal PY_{i,k}=\beta q_iZ_{i,k}$, \eqref{e4.potential-bound},
\eqref{e5.density-derivatives}, and
\eqref{e5.annular-density}--\eqref{e5.jacobi-annulus} yield
\begin{equation*}
\|L_VY_{i,k}\|_{L^1}
+\|\partial_p[((L_VY_{i,k})\circ\Phi_\xi)J_\xi]\|_{L^1}
\le C\delta_i\Lambda_\varepsilon^3
+C\delta_i\delta_*\Lambda_\varepsilon^4
+C\delta_i^4\Lambda_\varepsilon^3\int_{\sqrt{\delta_i}}^{r_0}r_x^{-5}\,dr_x
\le C\delta_i\Lambda_\varepsilon^3.
\end{equation*}
Combining \eqref{e2.17}, \eqref{e5.source-densities}
and \eqref{e5.jacobi-core}--\eqref{e5.jacobi-annulus}, we find
\[
\mathbb G_{(i,k),(j,l)}
=\delta_{ij}\int_{\mathbb R^3}qz_kz_l\,dy+O_{C^1}(\delta_*\Lambda_\varepsilon^3)
=\tfrac14\delta_{ij}\delta_{kl}+O_{C^1}(\delta_*\Lambda_\varepsilon^3),
\]
which proves \eqref{e5.gram-matrix}.

Finally, we split the integral in \eqref{e5.jacobi-residual-moments}
with the cutoff $\chi(r_0\sqrt{\delta_i}r)$.
By \eqref{e3.15}--\eqref{e3.17}, \eqref{e5.residual-derivative}, and
\eqref{e5.jacobi-core}, we may replace $Y_{i,k}$ by $z_k$ in the localized
integral with an $O_{C^1}(\delta_i^2\Lambda_\varepsilon^4)$ error.
Using also \eqref{e5.jacobian-derivatives}, we obtain
\[
\begin{aligned}
&\left\|
\int_{|y|<2\delta_i^{-1/2}}
\chi(r_0\sqrt{\delta_i}r)E_i(y)z_k(y)j_{\xi_i}(\delta_i y)\,dy
\right\|_{C^1}
\le C\delta_i^2\Lambda_\varepsilon^2,\\[2mm]
&\Biggl\|
\beta\delta_i\int_{|y|<2\delta_i^{-1/2}}
\chi(r_0\sqrt{\delta_i}r)q(y)\Bigl[B_{\xi_i}(y)+\textstyle\sum_{l=1}^3g_{i,l}y_l\Bigr]
 z_k(y)\bigl[j_{\xi_i}(\delta_i y)-1\bigr]\,dy
\Biggr\|_{C^1}
\le C\delta_i^3,\\[2mm]
&\left\|
\beta\delta_i\sum_{l=1}^3g_{i,l}
\int_{\mathbb R^3}
\bigl[1-\chi(r_0\sqrt{\delta_i}r)\bigr]q(y)y_lz_k(y)\,dy
\right\|_{C^1}
\le C\delta_i^{\frac{5}{2}}.
\end{aligned}
\]
The $B_{\xi_i}$ term in \eqref{e3.15} integrates to zero, also after parameter
differentiation, since $\chi(r_0\sqrt{\delta_i}r)qB_{\xi_i}$ is even and $z_k$ is odd.
For the $1-\chi(r_0\sqrt{\delta_i}r)$ part of \eqref{e5.jacobi-residual-moments},
\eqref{e3.14}, \eqref{e5.density-derivatives}, and
\eqref{e5.annular-density}--\eqref{e5.jacobi-annulus} yield
\[
C\bigl(\delta_i^2\Lambda_\varepsilon^3
+\delta_i\delta_*\Lambda_\varepsilon^4
+\delta_i\delta_*^3\Lambda_\varepsilon^3\bigr)
\le C\delta_i^2\Lambda_\varepsilon^4
\]
in $C^1$. Thus \eqref{e2.18} yields
\begin{equation*}
\int_{\partial\Omega}\mathcal R_\varepsilon(V)Y_{i,k}\,d\sigma
=-\beta\delta_i\sum_{l=1}^3g_{i,l}\int_{\mathbb R^3}q(y)y_lz_k(y)\,dy
+O_{C^1}(\delta_i^2\Lambda_\varepsilon^4)=-\frac{2\pi^4}{a^4}\delta_i g_{i,k}
+O_{C^1}(\delta_i^2\Lambda_\varepsilon^4),
\end{equation*}
which completes the proof.
\end{proof}

\begin{proposition}\label{t5.4}
Let $c=(c_{i,k})_{1\le i\le m,\,1\le k\le3}$ be the coefficients
in \eqref{e5.1}. For each sufficiently small $\varepsilon>0$, the map $
\xi\mapsto(\phi_{\varepsilon,\xi},c)$
is of class $C^2$ with values in
$H^3(\partial\Omega)\times\mathbb R^{3m}$.
Moreover, $\widetilde\phi=\phi_{\varepsilon,\xi}\circ\Phi_\xi$ satisfies
\begin{equation}\label{e5.material-correction}
\|\partial_p\widetilde\phi\|_\infty+|\partial_pc|
\le C\delta_*\Lambda_\varepsilon^5,\qquad 1\le p\le3m,
\end{equation}
uniformly for $\xi\in\mathcal K$.
\end{proposition}
\begin{proof}
By \eqref{e2.transported-green} and the weighted logarithmic estimate
used in \eqref{e2.14}, we obtain
\begin{equation}\label{e5.transported-kernel}
\sup_{x\in\partial\Omega}\int_{\partial\Omega}
|\partial_pG_{\partial\Omega}(\Phi_\xi x,\Phi_\xi z)|
(\omega\circ\Phi_\xi)(z)J_\xi(z)\,d\sigma(z)
\le C\Lambda_\varepsilon.
\end{equation}
Define
\[
\begin{aligned}
(\mathscr G_\xi f)(x)
&=\int_{\partial\Omega}G_{\partial\Omega}(\Phi_\xi x,\Phi_\xi z)
(Q_V\circ\Phi_\xi)(z)J_\xi(z)f(z)\,d\sigma(z),\\
\mathfrak A_\xi(f,c)
&=\Bigl(f-\mathscr G_\xi f-\sum_\alpha c_\alpha
[(\mathcal P^{-1}Q_\alpha)\circ\Phi_\xi],
\bigl(\ell_\alpha(f\circ\Phi_\xi^{-1})\bigr)_\alpha\Bigr).
\end{aligned}
\]
On $C^0(\partial\Omega)\times\mathbb R^{3m}$, we have
\begin{equation}\label{e5.transported-inverse}
\|\partial_p\mathfrak A_\xi\|
+\|\mathfrak A_\xi^{-1}\|\le C\Lambda_\varepsilon.
\end{equation}
The derivative bound follows from \eqref{e5.density-derivatives},
\eqref{e5.source-densities}, \eqref{e2.14}, and \eqref{e5.transported-kernel}.
To estimate $\|\mathfrak A_\xi^{-1}\|$, we compose $f$
and the first component of the data with $\Phi_\xi^{-1}$, where
$\Phi_\xi$ is defined in \eqref{e5.transport-map}.
Given $(F,b)\in C^0(\partial\Omega)\times\mathbb R^{3m}$,
choose $B_{\mathrm{tr}}\in\operatorname{span}\{Y_\alpha\}$ with
$\ell_\alpha(B_{\mathrm{tr}})=b_\alpha-\ell_\alpha(F)$.
By \eqref{e5.gram-matrix} and \eqref{e4.6}, we obtain
\[
\|B_{\mathrm{tr}}\|_\infty+\|L_VB_{\mathrm{tr}}\|_*
\le C(\|F\|_\infty+|b|).
\]
Writing $f=F+B_{\mathrm{tr}}+\psi_{\mathrm{lin}}$, we solve
\[
L_V\psi_{\mathrm{lin}}
=Q_VF-L_VB_{\mathrm{tr}}+\sum_\alpha d_\alpha Q_\alpha,
\qquad \ell_\alpha(\psi_{\mathrm{lin}})=0.
\]
Estimate \eqref{e4.10} yields the inverse bound in
\eqref{e5.transported-inverse}.
For fixed $\varepsilon>0$, applying $\mathcal P^{-1}$ to \eqref{e5.1}
yields, together with \eqref{e5.2}, a $C^2$ system in $(\xi,\phi,c)$.
The derivative with respect to $(\phi,c)$ at $\phi=0$ is
\[
(\psi,d)\longmapsto
\left(
\psi-\mathcal P^{-1}(Q_V\psi)
-\sum_\alpha d_\alpha\mathcal P^{-1}Q_\alpha,\;
\bigl(\ell_\alpha(\psi)\bigr)_\alpha
\right).
\]
The change of variables \eqref{e5.transport-map} transforms this operator
into $\mathfrak A_\xi$.
For $\|\phi\|_\infty\le1$, \eqref{e2.14} and $Q_V\le C\omega$ give
\begin{equation*}
\bigl\|\mathcal P^{-1}[D_\phi\mathcal N_V(\phi)\psi]\bigr\|_\infty
=\bigl\|\mathcal P^{-1}[Q_V(e^{\beta\phi}-1)\psi]\bigr\|_\infty\le C\Lambda_\varepsilon\|\phi\|_\infty\|\psi\|_\infty.
\end{equation*}
Combining the preceding estimate with \eqref{e5.3} and
\eqref{e5.transported-inverse}, we obtain
\begin{equation}\label{simaL}
C\Lambda_\varepsilon^2\|\phi\|_\infty
\le C\delta_*\Lambda_\varepsilon^4=o(1).
\end{equation}
Thus the derivative with respect to $(\phi,c)$ at
$(\phi_{\varepsilon,\xi},c)$ is invertible for sufficiently small
$\varepsilon$. The implicit function theorem therefore yields
$C^2$ dependence of $(\phi_{\varepsilon,\xi},c)$ on $\xi$ with values in
$C^0(\partial\Omega)\times\mathbb R^{3m}$.
To obtain the same dependence in $H^3(\partial\Omega)\times\mathbb R^{3m}$,
we rewrite \eqref{e5.1} as
\[
\phi=\mathcal P^{-1}\left(
Q_V\phi-\mathcal R_\varepsilon(V)+\mathcal N_V(\phi)
+\sum_\alpha c_\alpha Q_\alpha
\right).
\]
The expression inside $\mathcal P^{-1}$ is $C^2$ in $\xi$ with values in
$L^2(\partial\Omega)$, so the conclusion follows from the boundedness of
$\mathcal P^{-1}:L^2(\partial\Omega)\to H^3(\partial\Omega)$.
 We apply $\mathcal P^{-1}$ to
\eqref{e5.1}, use \eqref{e5.2}, and differentiate with respect to $p$
after the change of variables \eqref{e5.transport-map}. It holds
\begin{equation*}
\mathfrak A_\xi(\partial_p\widetilde\phi,\partial_pc)
=-(\partial_p\mathfrak A_\xi)(\widetilde\phi,c)
+\Bigl(-\partial_p\bigl[(\mathcal P^{-1}\mathcal R_\varepsilon(V))
\circ\Phi_\xi\bigr]+\beta^{-1}(\partial_p\mathscr G_\xi)
(e^{\beta\widetilde\phi}-1-\beta\widetilde\phi)+\mathscr G_\xi\bigl[(e^{\beta\widetilde\phi}-1)
\partial_p\widetilde\phi\bigr],\,0\Bigr).
\end{equation*}
By \eqref{e3.14}, \eqref{e5.density-derivatives}, \eqref{e2.14},
and \eqref{e5.transported-kernel}, we get
\begin{equation*}
\left\|\partial_p\bigl[(\mathcal P^{-1}\mathcal R_\varepsilon(V))
\circ\Phi_\xi\bigr]\right\|_\infty\le C\delta_*\Lambda_\varepsilon^3,\quad
\|\mathscr G_\xi f\|_\infty+
\|(\partial_p\mathscr G_\xi)f\|_\infty
\le C\Lambda_\varepsilon\|f\|_\infty,
\quad f\in C^0(\partial\Omega).
\end{equation*}
Applying $\mathfrak A_\xi^{-1}$ to the differentiated equation and
using Taylor's formula, \eqref{e5.3}, and \eqref{e5.transported-inverse},
we have
\[
\|\partial_p\widetilde\phi\|_\infty+|\partial_pc|
\le C\delta_*\Lambda_\varepsilon^5
+C\delta_*\Lambda_\varepsilon^4\|\partial_p\widetilde\phi\|_\infty,
\]
which together with \eqref{simaL} yields \eqref{e5.material-correction}.
\end{proof}

\begin{proposition}\label{t5.5}
As $\varepsilon\to0$,
\begin{equation}\label{e5.multiplier-expansion}
 \frac{c_{i,k}}{\delta_i}
 =-\frac{8\pi^4}{a^4}g_{i,k}+o_{C^1}(1),
 \qquad1\le i\le m,\quad1\le k\le3,
\end{equation}
uniformly for $\xi\in\mathcal K$.
\end{proposition}
\begin{proof}
Taking \(Y_\alpha\) as a test function in \eqref{e5.1}
 and using that \(\mathcal L_V\) is self-adjoint, we get
\[
 \sum_\gamma \mathbb G_{\alpha\gamma}c_\gamma
 =\int_{\partial\Omega}\mathcal R_\varepsilon(V)Y_\alpha\,d\sigma
 +\int_{\partial\Omega}\phi L_VY_\alpha\,d\sigma
 -\int_{\partial\Omega}\mathcal N_V(\phi)Y_\alpha\,d\sigma.
\]
Differentiating after the change of variables $x=\Phi_\xi(z)$,
we use  \eqref{e5.3}, \eqref{e5.jacobi-residual-moments} and
\eqref{e5.material-correction} to obtain
\begin{equation}\label{no1}
 \left\|\int_{\partial\Omega}\phi L_VY_{i,k}\,d\sigma\right\|_{C^1}
 \le C\delta_i\delta_*\Lambda_\varepsilon^8.
\end{equation}
The density estimates \eqref{e4.potential-bound},
\eqref{e5.density-derivatives}, and \eqref{e3.weight-bounds}, together with
 \eqref{e4.jacobi-norms},
\eqref{e5.jacobi-core}, and \eqref{e5.jacobi-annulus}, yield
\begin{equation}\label{no2}
 \left\|\int_{\partial\Omega}\mathcal N_V(\phi)Y_{i,k}\,d\sigma\right\|_{C^1}
 \le C\left(\Lambda_\varepsilon^2\|\phi\|_\infty^2
 +\|\phi\|_\infty\max_{1\le p\le3m}
 \|\partial_p\widetilde\phi\|_\infty\right)\le C\delta_*^2\Lambda_\varepsilon^8,
\end{equation}
by Taylor's formula, \eqref{e5.3}, and \eqref{e5.material-correction}.
From the second estimate in \eqref{e5.jacobi-residual-moments}, we deduce
\begin{equation}\label{no3}
 \int_{\partial\Omega}\mathcal R_\varepsilon(V)Y_{i,k}\,d\sigma
 =-\frac{2\pi^4}{a^4}\delta_i g_{i,k}
 +O_{C^1}(\delta_i^2\Lambda_\varepsilon^4).
\end{equation}
Moreover, \eqref{e5.gram-matrix} implies
$\mathbb G^{-1}=4I+O_{C^1}(\delta_*\Lambda_\varepsilon^3)$.
Combining \eqref{no1}--\eqref{no3}, and using $\delta_i\asymp\delta_*$
together with $|\partial_p\log\delta_i|\le C$
from \eqref{e3.10}, we obtain
\[
 \left\|\frac{c_{i,k}}{\delta_i}
 +\frac{8\pi^4}{a^4}g_{i,k}\right\|_{C^1}
 \le C\delta_*\Lambda_\varepsilon^8\longrightarrow0.
\]
\end{proof}

\subsection{The reduced energy}\label{s5.2}

The boundary functional is
\begin{equation}
 I_\varepsilon(v)
 =\frac12\langle \mathcal{P}v,v\rangle
 -\frac{\varepsilon}{24\pi^2}\int_{\partial\Omega}he^{24\pi^2v}\,d\sigma,
 \qquad v\in H^{3/2}(\partial\Omega).
 \label{e5.4}
\end{equation}
 This is \eqref{e1.2} restricted to
\(v\mapsto\mathcal Ev\). Set
\(\mathcal{K}_\varepsilon(\xi)=I_\varepsilon(v_{\varepsilon,\xi})\).

\begin{proposition}\label{t5.7}
There exists $C_{\varepsilon,m}$, independent of $\xi$, such that
\begin{equation}\label{e5.energy-expansion}
 \mathcal K_\varepsilon=C_{\varepsilon,m}-\mathcal F_m+o_{C^2}(1)
\end{equation}
uniformly for $\xi\in\mathcal K$ as $\varepsilon\to0$. Moreover, for sufficiently small $\varepsilon$,
\begin{equation}\label{e5.8}
 d\mathcal K_\varepsilon=0
 \quad\Longleftrightarrow\quad
 c_{i,k}=0\quad(1\le i\le m,\ 1\le k\le3).
\end{equation}
\end{proposition}
\begin{proof}
We have $I_\varepsilon(V)
 =\frac12\sum_{i=1}^m\int_{\partial\Omega}Vq_i\,d\sigma
 -\frac{\varepsilon}{\beta}\int_{\partial\Omega}he^{\beta V}\,d\sigma.$
Changing variables $x=\exp_{\xi_i}(\delta_i y)$ and applying
\eqref{e3.11} and \eqref{e5.source-densities}, we obtain
\begin{equation}\label{e5.5}
\begin{aligned}
 \int_{\partial\Omega}Vq_i\,d\sigma
 ={}&\int_{\mathbb R^3}q(y)\left[
 u_0(y)-\frac{\log\delta_i}{4\pi^2}+\frac{\log a}{4\pi^2}
 +\mathcal R_i(\xi)+\delta_iB_{\xi_i}(y)
 +\delta_i\langle A_i,y\rangle\right]dy
 +O(\delta_i^{3/2}\Lambda_\varepsilon)\\
 ={}&\int_{\mathbb R^3}u_0q\,dy
 -\frac{\log\delta_i}{4\pi^2}+\frac{\log a}{4\pi^2}
 +\mathcal R_i(\xi)+O(\delta_i\Lambda_\varepsilon).
\end{aligned}
\end{equation}
On the other hand, $\mathcal PV=\sum_iq_i$,
$\int_{\partial\Omega}q_i\,d\sigma=1$, and
\eqref{e3.14} imply
\begin{equation}\label{e5.6}
 \varepsilon\int_{\partial\Omega}he^{\beta V}\,d\sigma
 =m-\int_{\partial\Omega}\mathcal R_\varepsilon(V)\,d\sigma
 =m+O(\delta_*\Lambda_\varepsilon).
\end{equation}
Combining \eqref{e5.5}--\eqref{e5.6} with the choice of scales in
\eqref{e3.10}, we obtain
\begin{equation}\label{e5.7}
 I_\varepsilon(V)
 =C_{\varepsilon,m}-\frac12\sum_i\mathcal R_i(\xi)
   -\frac1{24\pi^2}\sum_i\log h(\xi_i)+o(1)
   =C_{\varepsilon,m}-\mathcal F_m(\xi)+o(1),
\end{equation}
where $C_{\varepsilon,m}
 =\frac m2\int_{\mathbb R^3}u_0q\,dy
 -\frac{m\log a}{8\pi^2}
 -\frac{m(\log\varepsilon+1)}{24\pi^2}.$
Taylor expansion and \eqref{e5.1}--\eqref{e5.2} yield
\[
\begin{aligned}
 I_\varepsilon(V+\phi)-I_\varepsilon(V)
 &=\frac12\int_{\partial\Omega}\mathcal R_\varepsilon(V)\phi\,d\sigma
   +\frac12\int_{\partial\Omega}\mathcal N_V(\phi)\phi\,d\sigma
   +O(\|\phi\|_\infty^3).
\end{aligned}
\]
Hence \eqref{e4.potential-bound}, \eqref{e3.weight-bounds}, and
\eqref{e5.3} imply
\begin{equation}\label{e5.energy-correction}
 |I_\varepsilon(V+\phi)-I_\varepsilon(V)|
 \le C\delta_*\Lambda_\varepsilon\|\phi\|_\infty
      +C\|\phi\|_\infty^3
 \le C\delta_*^2\Lambda_\varepsilon^5=o(1).
\end{equation}
Let us now compute the derivatives of $\mathcal K_\varepsilon$.
For $1\le i\le m$, $1\le k\le3$, and $1\le p\le3m$, set $
 \mathbb M_{(i,k),p}
 =\int_{\partial\Omega}Q_{i,k}\partial_pv_{\varepsilon,\xi}\,d\sigma.$
We claim that
\begin{equation}\label{e5.translation-matrix}
 \delta_i\mathbb M_{(i,k),p}
 =\frac{a}{32\pi^2}\mathbb A_{(i,k),p}+o_{C^1}(1).
\end{equation}
Write $C_i=-\frac{\log\delta_i}{4\pi^2}
     +\frac{\log a}{4\pi^2}+\mathcal R_i(\xi).$
By \eqref{e3.11}, we have
\[
 V(\exp_{\xi_i}(\delta_i y))-C_i
 =u_0(y)+\delta_iB_{\xi_i}(y)
   +\delta_i\langle A_i,y\rangle+E_{i,0}(y).
\]
Note that $qz_k$ is odd, while $u_0$ and $B_{\xi_i}$ are even.
Using \eqref{e5.material-V-bound} for the tail and a smooth radial
cutoff at $|y|\asymp\delta_i^{-1/2}$,
\eqref{e5.jacobian-derivatives}--\eqref{e5.inner-remainder-derivatives}
imply
\begin{equation}\label{e5.scalar-moment}
 \left\|\int_{\partial\Omega}Q_{i,k}(V-C_i)\,d\sigma\right\|_{C^2}
 \le C\delta_i+C\delta_i^2\Lambda_\varepsilon^2
       +C\Lambda_\varepsilon\int_{|y|>\delta_i^{-\frac{1}{2}}}|qz_k|\,dy
       \le C\delta_*\Lambda_\varepsilon^3.
\end{equation}
Differentiating the density in \eqref{e5.source-densities}, we get
\begin{equation}\label{e5.density-derivative}
 \delta_i(\partial_pQ_{i,k})(\exp_{\xi_i}(\delta_i y))
 \delta_i^3j_{\xi_i}(\delta_i y)
 =-\sum_{l=1}^3\mathbb A_{(i,l),p}\partial_{y_l}(qz_k)(y)
 +E^Q_{i,k,p}(y).
\end{equation}
Here we take $\partial_pQ_{i,k}$ at fixed $x$ and set the left-hand side of
\eqref{e5.density-derivative} equal to zero for $|y|\ge 2r_0/\delta_i$,
thereby defining $E^Q_{i,k,p}$ on $\mathbb R^3$.
Using $
 \left.\partial_py\right|_x
 =-\delta_i^{-1}(\mathbb A_{(i,l),p})_{l=1}^3+O_{C^1}(1+|y|),$
\eqref{e5.jacobian-derivatives}--\eqref{e5.normalization-derivatives},
and $|\nabla(qz_k)|\le C(1+|y|)^{-8}$, we find
\[
 |E^Q_{i,k,p}(y)|+|\partial_\kappa E^Q_{i,k,p}(y)|
 \le C\delta_i(1+|y|)^{-7}
      +C\delta_i^2(1+|y|)^{-5},
 \qquad1\le\kappa\le3m,
\]
where derivatives of rescaled quantities are taken at fixed $y$.
This estimate includes the cutoff terms at $|y|\asymp\delta_i^{-1}$
and the extension of the principal term. In particular,
\begin{equation}\label{e5.density-error}
 \int_{\mathbb R^3}(1+\log(2+|y|))
 \bigl(|E^Q_{i,k,p}|+|\partial_\kappa E^Q_{i,k,p}|\bigr)\,dy
 \le C\delta_i\Lambda_\varepsilon^2.
\end{equation}
It holds
\begin{equation*}
 \delta_i\int_{\partial\Omega}Q_{i,k}\partial_pV\,d\sigma
 =\delta_i\partial_p\int_{\partial\Omega}Q_{i,k}(V-C_i)\,d\sigma-\delta_i\int_{\partial\Omega}(\partial_pQ_{i,k})(V-C_i)\,d\sigma
 +\delta_i(\partial_pC_i)\int_{\partial\Omega}Q_{i,k}\,d\sigma.
\end{equation*}
The first term is $o_{C^1}(1)$ by \eqref{e5.scalar-moment}.
The last is also $o_{C^1}(1)$, since
$\partial_pC_i=O_{C^1}(1)$ and
$\|\int_{\partial\Omega}Q_{i,k}\,d\sigma\|_{C^2}\le C\delta_i^2$
by \eqref{e5.source-densities} and oddness.
Thus \eqref{e5.density-derivative}--\eqref{e5.density-error},
\eqref{e3.11}, and \eqref{e5.inner-remainder-derivatives} give
\begin{equation}\label{e5.approximate-translation-moment}
 \delta_i\int_{\partial\Omega}Q_{i,k}\partial_pV\,d\sigma
 =\sum_{l=1}^3\mathbb A_{(i,l),p}
   \int_{\mathbb R^3}\partial_{y_l}(qz_k)u_0\,dy+o_{C^1}(1)
   =\frac{a}{32\pi^2}\mathbb A_{(i,k),p}+o_{C^1}(1).
\end{equation}
For the remaining term, differentiation of \eqref{e5.2} yields
\[
 \int_{\partial\Omega}Q_{i,k}\partial_p\phi\,d\sigma
 =-\int_{\partial\Omega}(\partial_pQ_{i,k})\phi\,d\sigma.
\]
It follows from \eqref{e5.density-derivative}--\eqref{e5.density-error},
\eqref{e5.3}, and \eqref{e5.material-correction} that
\[
 \left\|\delta_i\int_{\partial\Omega}Q_{i,k}\partial_p\phi\,d\sigma\right\|_{C^1}
 \le C\delta_*\Lambda_\varepsilon^5=o(1),
\]
which proves \eqref{e5.translation-matrix}.
Differentiating $\mathcal K_\varepsilon$ and using
\eqref{e5.projected-residual}, we obtain
\begin{equation}\label{e5.reduced-gradient}
 \partial_p\mathcal K_\varepsilon
 =\sum_{i=1}^m\sum_{k=1}^3c_{i,k}\mathbb M_{(i,k),p}.
\end{equation}
By \eqref{e5.multiplier-expansion}, \eqref{e5.translation-matrix},
and $a^3=\pi^2/4$, we get
\[
 \partial_p\mathcal K_\varepsilon
 =-\sum_{i=1}^m\sum_{k=1}^3g_{i,k}\mathbb A_{(i,k),p}+o_{C^1}(1)
 =-\partial_p\mathcal F_m+o_{C^1}(1).
\]
Together with \eqref{e5.7}--\eqref{e5.energy-correction}, this proves
\eqref{e5.energy-expansion}.
Finally, the matrix $(\delta_i\mathbb M_{(i,k),p})_{(i,k),p}$
converges uniformly to $\frac{a}{32\pi^2}\mathbb A$, which is invertible.
Hence $\mathbb M$ is invertible for small $\varepsilon$, and \eqref{e5.reduced-gradient} implies \eqref{e5.8}.
\end{proof}

\section{Existence of multi-bubble solutions}\label{s6}

\begin{proof}[Proof of Theorem~\ref{t1.1}]
Let $D$ be as in \eqref{e1.9}. It follows from \eqref{e5.energy-expansion} that
\[
 C_{\varepsilon,m}-\mathcal K_\varepsilon
 =\mathcal F_m+o_{C^1(\overline D)}(1).
\]
By the $C^1$-stability of $K$, $\mathcal K_\varepsilon$ has a critical point
$\xi_\varepsilon\in D$ for sufficiently small $\varepsilon>0$. Moreover,
$\nabla\mathcal F_m(\xi_\varepsilon)=o(1)$.
In view of compactness of $\overline D$ and \eqref{e1.9}, every limit point
of $\xi_\varepsilon$ belongs to $K$, and hence
$\dist(\xi_\varepsilon,K)\to0$.

Equations \eqref{e5.8} and \eqref{e5.projected-residual} imply that
$u_\varepsilon=\mathcal E v_{\varepsilon,\xi_\varepsilon}$ solves
\eqref{e1.1}. The required decomposition and estimates follow from
\eqref{e3.10} and \eqref{e5.3}:
\begin{equation*}
\delta_{i,\varepsilon}=\delta_i(\varepsilon,\xi_\varepsilon)
\asymp\varepsilon^{\frac{1}{3}},\quad 1\le i\le m,\quad
\|\phi_{\varepsilon,\xi_\varepsilon}\|_{L^\infty(\partial\Omega)}
+\|\phi_{\varepsilon,\xi_\varepsilon}\|_{H^{\frac{3}{2}}(\partial\Omega)}
\le C\varepsilon^{\frac{1}{3}}|\log\varepsilon|^2\rightarrow0.
\end{equation*}
With $V$ and $q_i$ evaluated at $\xi_\varepsilon$, we obtain
\[
\varepsilon he^{\beta v_{\varepsilon,\xi_\varepsilon}}-\sum_iq_i
=\varepsilon he^{\beta V}
\bigl(e^{\beta\phi_{\varepsilon,\xi_\varepsilon}}-1\bigr)
-\mathcal R_\varepsilon(V).
\]
By \eqref{e3.weight-bounds}, \eqref{e3.14},
\eqref{e4.potential-bound} and \eqref{e5.3}, we get
\begin{equation}\label{e6.mass-comparison}
\left\|\varepsilon he^{\beta v_{\varepsilon,\xi_\varepsilon}}
-\sum_iq_i\right\|_{L^1}
\le C\bigl(\|\phi_{\varepsilon,\xi_\varepsilon}\|_\infty
+\|\mathcal R_\varepsilon(V)\|_*\bigr)
\le C\delta_*\Lambda_\varepsilon^2\rightarrow0.
\end{equation}
Moreover, \eqref{e3.1} and \eqref{e3.2} imply, for each fixed $r>0$,
\[
\int_{\partial\Omega}q_i\,d\sigma=1,\qquad
\int_{\partial\Omega\setminus B_r(\xi_i)}q_i\,d\sigma
\le Cr^{-2}\delta_i^2\rightarrow0.
\]
Thus $q_i\,d\sigma\rightharpoonup\delta_{\xi_i^*}$ whenever
$\xi_\varepsilon\to\xi^*$. Together with \eqref{e6.mass-comparison},
this proves \eqref{e1.mass-convergence}.
For a nondegenerate critical point $\xi^*$, take $K=\{\xi^*\}$.
\end{proof}

\begin{proof}[Proof of Corollary~\ref{t1.2}]
Fix $c>\min_{\mathcal M_m}\mathcal F_m$ and set
\[
 S_c=\{\xi\in\mathcal M_m:\mathcal F_m(\xi)\le c\}.
\]
For $m\ge2$, since $\mathcal F_m(\xi)\to+\infty$ as
$\min_{i\ne j}d(\xi_i,\xi_j)\to0$, the set $S_c$ is closed in
the compact space $(\partial\Omega)^m$ and contained in $\mathcal M_m$;
hence $S_c\subset\subset\mathcal M_m$.
For $m=1$, $S_c$ is compact as a closed subset of $\mathcal M_1=\partial\Omega$.
Choose the compact set $\mathcal K$  so that
$S_c\subset\operatorname{int}\mathcal K$, and let $\xi_\varepsilon$
maximize $\mathcal K_\varepsilon$ on $S_c$. By \eqref{e5.energy-expansion}, we have
\[
0\le\mathcal F_m(\xi_\varepsilon)-\min_{\mathcal M_m}\mathcal F_m
\le2\|\mathcal K_\varepsilon-C_{\varepsilon,m}
       +\mathcal F_m\|_{C^0(S_c)}=o(1).
\]
Thus $\mathcal F_m(\xi_\varepsilon)<c$ for small $\varepsilon$, so
$\xi_\varepsilon$ is an interior point of $S_c$ and
$d\mathcal K_\varepsilon(\xi_\varepsilon)=0$.
The solution and correction estimate follow from \eqref{e5.3}, \eqref{e5.projected-residual}, and \eqref{e5.8}.
Equation \eqref{e6.mass-comparison} also proves \eqref{e1.mass-convergence} along each subsequence with convergent centers.
\end{proof}

In the rest of the paper, we write $\xi_i=\xi_{i,\varepsilon}$,
$\delta_i=\delta_{i,\varepsilon}$, and
$v_\varepsilon=v_{\varepsilon,\xi_\varepsilon}$, and set
$\rho_\varepsilon=\varepsilon h e^{24\pi^2v_\varepsilon}$.
With this notation, \eqref{e3.12} and \eqref{e5.3} yield
\[
v_\varepsilon(\exp_{\xi_i}(\delta_i y))
+\frac{\log\delta_i}{8\pi^2}
+\frac{\log(\varepsilon h(\xi_i))}{24\pi^2}
\longrightarrow u_0(y)
\]
locally uniformly in $\mathbb R^3$ as $\varepsilon\to0$,
for $1\le i\le m$.
\section{Qualitative properties of the concentrating solutions}\label{s7}

\subsection{Local uniqueness and  nondegeneracy}\label{s7.1}

Let $\xi^*$ be a nondegenerate critical point of $\mathcal F_m$.
We work in fixed normal coordinates on $\mathcal M_m$
centered at $\xi^*$ and write $D^2$ for the matrix of second
coordinate derivatives in this subsection.
At critical points, this matrix represents the Riemannian Hessian.
Set $A_0=-D^2\mathcal F_m(\xi^*)$ and choose concentric coordinate balls
$\xi^*\in\mathcal U\subset\subset\mathcal U^+\subset\subset\mathcal M_m$ such that
\[
\sup_{\mathcal U^+}\|A_0^{-1}(-D^2\mathcal F_m-A_0)\|<\tfrac14.
\]
We take $\Omega_\star\subset\mathcal U$ in Theorem~\ref{t1.3}
and apply the reduction in a neighborhood of $\overline{\mathcal U^+}$.
For fixed $\varepsilon$, we write
$\delta_i(\xi)=\delta_i(\varepsilon,\xi)$.

\begin{lemma}\label{t7.2}
Let $A>0$ be fixed. For  sufficiently small $\varepsilon>0$,
if $v\in L^\infty(\partial\Omega)$ and $\xi\in\mathcal U$ satisfy
\begin{equation}\label{e7.1}
\|v-V_\xi\|_\infty\le A\delta_*(\xi)\Lambda_\varepsilon^2,
\end{equation}
then there exists $\eta\in\mathcal U^+$ such that
\begin{equation}\label{le5.1}
\begin{gathered}
\int_{\partial\Omega}Q_{i,k}(\eta)(v-V_\eta)\,d\sigma=0,
\qquad 1\le i\le m,\quad1\le k\le3,\\
\max_{1\le i\le m}d(\eta_i,\xi_i)
\le C_A\delta_*(\xi)^2\Lambda_\varepsilon^2,
\qquad
\|v-V_\eta\|_\infty\le C_A\delta_*(\xi)\Lambda_\varepsilon^2.
\end{gathered}
\end{equation}
\end{lemma}
\begin{proof}
For $t=(t_{i,k})_{1\le i\le m,\,1\le k\le3}\in\mathbb R^{3m}$, set
\[
\theta(t)=(\theta_1(t),\ldots,\theta_m(t)),\qquad
\theta_i(t)=\exp_{\xi_i}\!\left(
\delta_i(\xi)\sum_{k=1}^3t_{i,k}e_{i,k}\right),
\]
where $e_{i,k}$ is evaluated at $\xi_i$.
For $|t|$ smaller than a fixed constant and $\varepsilon$ sufficiently small,
$\theta(t)\in\mathcal U^+$. Define $F=(F_{i,k})_{i,k}$ by
\[
F_{i,k}(t)=\int_{\partial\Omega}
Q_{i,k}(\theta(t))(v-V_{\theta(t)})\,d\sigma.
\]
In the estimates below, $q_i$ and $Q_{i,k}$ are evaluated at
$\theta(t)$, with scales prescribed by \eqref{e3.10}, and derivatives
with respect to $t$ are taken at fixed $x$.
Differentiating \eqref{e3.1} and using
$\int_{\partial\Omega}q_i\,d\sigma=1$, we obtain
\begin{equation*}
\int_{\partial\Omega}D_tq_i\,d\sigma=0,\quad
\int_{\partial\Omega}D_t^2q_i\,d\sigma=0,\quad
|D_tq_i(x)|+|D_t^2q_i(x)|
\le\frac{C\delta_i(\theta(t))^{-3}}
{(1+d(x,\theta_i(t))/\delta_i(\theta(t)))^4}.
\end{equation*}
Then combining \eqref{e4.1} and  $Q_{i,k}=q_iZ_{i,k}$, we infer
\begin{equation}\label{bounda}
\|D_tV_{\theta(t)}\|_\infty+\|D_t^2V_{\theta(t)}\|_\infty
+\sum_{i=1}^m\sum_{k=1}^3
\bigl(\|D_tQ_{i,k}\|_{L^1}+\|D_t^2Q_{i,k}\|_{L^1}\bigr)\le C.
\end{equation}
Since $\|Q_{i,k}\|_{L^1}\le C$, \eqref{e7.1} yields
$|F(0)|\le C_A\delta_*(\xi)\Lambda_\varepsilon^2$.
Moreover,
\[
\partial_{t_{j,l}}F_{i,k}(t)
=\int_{\partial\Omega}(\partial_{t_{j,l}}Q_{i,k})
(v-V_{\theta(t)})\,d\sigma
-\int_{\partial\Omega}Q_{i,k}\partial_{t_{j,l}}V_{\theta(t)}\,d\sigma.
\]
At $t=0$, \eqref{e5.approximate-translation-moment} implies
$DF(0)=-\frac{a}{32\pi^2}I+o(1)$. Using
\eqref{bounda} and
$\|v-V_{\theta(t)}\|_\infty\le C_A\delta_*(\xi)\Lambda_\varepsilon^2+C|t|$, we obtain
\[
\|DF(t)-DF(0)\|\le C|t|.
\]
For a sufficiently large constant $C_A$ and small $\varepsilon$,
the map $t\mapsto t+\frac{32\pi^2}{a}F(t)$ is a contraction of
$\{t:|t|\le C_A\delta_*(\xi)\Lambda_\varepsilon^2\}$ into itself.
Its fixed point satisfies $F(t)=0$.
Taking $\eta=\theta(t)$, we obtain
\[
\max_i d(\eta_i,\xi_i)\le\delta_*(\xi)|t|,
\qquad
\|v-V_\eta\|_\infty
\le\|v-V_\xi\|_\infty+C|t|,
\]
which proves \eqref{le5.1}.
\end{proof}

\begin{proof}[Proof of Theorem~\ref{t1.3}: local uniqueness]
Let $v=\widetilde u_\varepsilon|_{\partial\Omega}$.
Since $\varepsilon^{1/3}\le C\delta_*(\xi)$, \eqref{e1.uniqueness-class}
gives $
\|v-V_\xi\|_\infty
\le CA\delta_*(\xi)\Lambda_\varepsilon^2.$
Applying Lemma~\ref{t7.2} with $CA$ in place of $A$ and using
 \eqref{e1.4}, we obtain $\eta\in\mathcal U^+$
such that $\psi=v-V_\eta$ satisfies
\[
L_{V_\eta}\psi=-\mathcal R_\varepsilon(V_\eta)
+\mathcal N_{V_\eta}(\psi),\qquad
\int_{\partial\Omega}Q_{i,k}(\eta)\psi\,d\sigma=0.
\]
By Lemma~\ref{t7.2} and \eqref{e5.3}, we have $
\|\psi\|_\infty+\|\phi_{\varepsilon,\eta}\|_\infty
\le C_A\delta_*(\eta)\Lambda_\varepsilon^2.$
Applying \eqref{e4.10} to $\psi-\phi_{\varepsilon,\eta}$
and using Taylor's formula and \eqref{e4.potential-bound},
 we obtain
\begin{equation*}
\|\psi-\phi_{\varepsilon,\eta}\|_\infty
\le C\Lambda_\varepsilon
\bigl(\|\psi\|_\infty+\|\phi_{\varepsilon,\eta}\|_\infty\bigr)
\|\psi-\phi_{\varepsilon,\eta}\|_\infty\le C_A\delta_*(\eta)\Lambda_\varepsilon^3
\|\psi-\phi_{\varepsilon,\eta}\|_\infty.
\end{equation*}
Due to $\delta_*(\eta)\Lambda_\varepsilon^3\to0$,
$\psi=\phi_{\varepsilon,\eta}$.
We derive from Proposition~\ref{t4.7} and \eqref{e5.8} that $c_{i,k}(\eta)=0$ and
$d\mathcal K_\varepsilon(\eta)=0$.
Using \eqref{e5.energy-expansion} and the choice of $\mathcal U^+$, we get
\[
\sup_{\mathcal U^+}
\|A_0^{-1}(D^2\mathcal K_\varepsilon-A_0)\|<\tfrac12.
\]
For small $\varepsilon$, the center $\xi_\varepsilon$ of the constructed
solution belongs to $\mathcal U$ and is also a critical point of
$\mathcal K_\varepsilon$.
In the fixed coordinates, integration along the segment in the convex set $\mathcal U^+$ yields
\[
0=\int_0^1D^2\mathcal K_\varepsilon
\bigl(\xi_\varepsilon+s(\eta-\xi_\varepsilon)\bigr)
(\eta-\xi_\varepsilon)\,ds,
\qquad
|\eta-\xi_\varepsilon|\le\tfrac12|\eta-\xi_\varepsilon|.
\]
Thus $\eta=\xi_\varepsilon$ and $v=v_{\varepsilon,\xi_\varepsilon}$.
Uniqueness of the extension problem \eqref{e1.3} yields
$\widetilde u_\varepsilon=u_\varepsilon$ in $\Omega$.
\end{proof}

\begin{proof}[Proof of Theorem~\ref{t1.3}:  nondegeneracy]
Let $\Phi\in\mathcal H$ be a weak solution of
\eqref{e1.linearized-problem}, and set
$f=\Phi|_{\partial\Omega}$.
Following uniqueness for the extension problem \eqref{e1.3}, $\Phi=\mathcal Ef$. Then
\[
L_\varepsilon f=0,\qquad
L_\varepsilon=\mathcal P-\beta\rho_\varepsilon.
\]
For each fixed $\varepsilon$, we have
$f\in L^2(\partial\Omega)$ and
$\rho_\varepsilon\in L^\infty(\partial\Omega)$, so
$\mathcal Pf=\beta\rho_\varepsilon f\in L^2(\partial\Omega)$.
Elliptic regularity and Sobolev embedding therefore give
$f\in H^3(\partial\Omega)\hookrightarrow L^\infty(\partial\Omega)$.
For $f,g\in H^{3/2}(\partial\Omega)$, define
\begin{equation}\label{e7.2}
\mathcal B_\varepsilon(f,g)=\langle L_\varepsilon f,g\rangle
=\langle\mathcal Pf,g\rangle
-\beta\int_{\partial\Omega}\rho_\varepsilon fg\,d\sigma.
\end{equation}
Set $
T_p=\left.\partial_pv_{\varepsilon,\xi}\right|_{\xi=\xi_\varepsilon},$ $1\le p\le3m,$
and
\begin{equation}\label{matr}
X_{\mathrm{tr},\varepsilon}
=\{g\in H^{\frac{3}{2}}(\partial\Omega):\ell_{i,k}(g)=0,
\ 1\le i\le m,\ 1\le k\le3\},
\end{equation}
where the moments are evaluated at $\xi_\varepsilon$.
From \eqref{e5.translation-matrix}, the matrix
$(\ell_{i,k}(T_p))_{(i,k),p}$ is invertible. Consequently,
\begin{equation}\label{e7.translation-splitting}
H^{\frac{3}{2}}(\partial\Omega)=X_{\mathrm{tr},\varepsilon}
\oplus\operatorname{span}\{T_p:1\le p\le3m\}.
\end{equation}
In view of  $c_{i,k}(\xi_\varepsilon)=0$,
\eqref{e5.projected-residual} and \eqref{e5.reduced-gradient}, we deduce
\begin{equation}\label{e7.translation-hessian}
L_\varepsilon T_p
=\sum_{i=1}^m\sum_{k=1}^3(\partial_pc_{i,k})Q_{i,k},\quad
\mathcal B_\varepsilon(T_p,w_\partial)=0
\quad(w_\partial\in X_{\mathrm{tr},\varepsilon}),\quad
\mathcal B_\varepsilon(T_p,T_\kappa)
=\partial_{p\kappa}^2\mathcal K_\varepsilon(\xi_\varepsilon).
\end{equation}
Here all terms on the right-hand side are evaluated at $\xi_\varepsilon$,
and $1\le p,\kappa\le3m$.
Write $f=w_\partial+\sum_{p=1}^{3m}a_pT_p$, with
$w_\partial\in X_{\mathrm{tr},\varepsilon}$ and
$\mathbf a=(a_p)_{p=1}^{3m}$.
Testing $L_\varepsilon f=0$ with $T_\kappa$ yields
\[
D^2\mathcal K_\varepsilon(\xi_\varepsilon)\mathbf a=0.
\]
By \eqref{e5.energy-expansion} and $\xi_\varepsilon\to\xi^*$,
$D^2\mathcal K_\varepsilon(\xi_\varepsilon)\to A_0$.
Thus $\mathbf a=0$ and $f\in X_{\mathrm{tr},\varepsilon}$.
Moreover, $L_\varepsilon f=0$ implies
\[
L_{V_{\xi_\varepsilon}}f
=\beta\varepsilon h
\bigl(e^{\beta v_\varepsilon}-e^{\beta V_{\xi_\varepsilon}}\bigr)f.
\]
From \eqref{e4.potential-bound}, \eqref{e4.10} and \eqref{e5.3},
 we obtain
\[
\|f\|_\infty
\le C\Lambda_\varepsilon
\|\phi_{\varepsilon,\xi_\varepsilon}\|_\infty\|f\|_\infty
\le C\delta_*(\xi_\varepsilon)\Lambda_\varepsilon^3\|f\|_\infty,
\]
so that $f=0$ and $\Phi=\mathcal Ef=0$.
\end{proof}

\subsection{The Morse index}

Throughout this subsection,  we write $V=V_{\xi_\varepsilon}$.
The quadratic form on $\mathcal H$ associated with the linearized problem is
\begin{equation*}
\mathcal Q_\varepsilon(\Phi)
=\int_\Omega\bigl((\Delta\Phi)^2+\Phi^2\bigr)\,dx
-\beta\int_{\partial\Omega}\rho_\varepsilon\Phi^2\,d\sigma.
\end{equation*}

\begin{lemma}\label{t7.3}
The negative indices of $\mathcal B_\varepsilon$ on $H^{3/2}(\partial\Omega)$
and $\mathcal Q_\varepsilon$ on $\mathcal H$ coincide.
\end{lemma}
\begin{proof}
Every $\Phi\in\mathcal H$ has the unique decomposition
$\Phi=\mathcal Ef+w$, where
$f=\Phi|_{\partial\Omega}$ and $w\in H_0^2(\Omega)$.
Thus
\begin{equation}\label{51}
\mathcal H
=\mathcal E H^{\frac{3}{2}}(\partial\Omega)\oplus H_0^2(\Omega).
\end{equation}
Using the weak formulation of \eqref{e1.3}, $
\int_\Omega
\bigl(\Delta\mathcal Ef\,\Delta w+\mathcal Ef\,w\bigr)\,dx=0.$
Since $w|_{\partial\Omega}=0$, we obtain
\begin{equation}\label{52}
\mathcal Q_\varepsilon(\mathcal Ef+w)
=\mathcal B_\varepsilon(f,f)
+\int_\Omega\bigl((\Delta w)^2+w^2\bigr)\,dx.
\end{equation}
Moreover,
\begin{equation}\label{53}
\mathcal Q_\varepsilon(w)
=\int_\Omega\bigl((\Delta w)^2+w^2\bigr)\,dx>0
\qquad
\text{for }0\ne w\in H_0^2(\Omega).
\end{equation}
Combining \eqref{51}--\eqref{53}, we have
\[
\ind(\mathcal Q_\varepsilon)
=\ind(\mathcal B_\varepsilon).
\]
\end{proof}

Equip $\mathcal H_\partial=H^{3/2}(\partial\Omega)$ with the inner product $\langle u,v\rangle_{\mathcal H_\partial}=\langle\mathcal Pu,v\rangle,$
whose norm is equivalent to the Sobolev norm by \eqref{e2.coercivity}, and set
\[
\mathcal H_0=\bigoplus_{i=1}^mL^2(\mathbb R^3,q\,dy),
\qquad
\langle f,g\rangle_{\mathcal H_0}
=\sum_{i=1}^m\int_{\mathbb R^3}qf_ig_i\,dy.
\]
In the  normal coordinates
$x=\exp_{\xi_i}(\delta_i y)$, we have
\begin{equation}\label{e7.rescaled-density}
q_i\,d\sigma=q(y)g_{i,\varepsilon}(y)\,dy,
\qquad
g_{i,\varepsilon}(y)
=b_{\delta_i,\xi_i}\chi(\delta_i|y|)j_{\xi_i}(\delta_i y),
\end{equation}
where we set $g_{i,\varepsilon}(y)=0$ for $|y|\ge 2r_0/\delta_i$.
Define $B_\varepsilon:\mathcal H_\partial\to\mathcal H_0$ by
\[
(B_\varepsilon u)_i(y)=
\begin{cases}
\sqrt{\beta}\,g_{i,\varepsilon}(y)^{\frac{1}{2}}
u(\exp_{\xi_i}(\delta_i y)),
& |y|<2r_0/\delta_i,\\
0, & |y|\ge 2r_0/\delta_i,
\end{cases}
\qquad 1\le i\le m.
\]
Then, by \eqref{e7.rescaled-density}, we obtain $
\|B_\varepsilon u\|_{\mathcal H_0}^2
=\beta\sum_{i=1}^m\int_{\partial\Omega}q_i u^2\,d\sigma.$
 For fixed $\varepsilon>0$, we have
$q_i\in L^\infty(\partial\Omega)$ by \eqref{e3.1}, and hence
\[
\|B_\varepsilon u\|_{\mathcal H_0}^2
\le\beta\left\|\sum_{i=1}^m q_i\right\|_{L^\infty(\partial\Omega)}
\|u\|_{L^2(\partial\Omega)}^2.
\]
Combining this estimate with the compact embedding
$H^{3/2}(\partial\Omega)\hookrightarrow L^2(\partial\Omega)$,
we conclude that $B_\varepsilon:\mathcal H_\partial\to\mathcal H_0$
is compact.
We now denote by $B_\varepsilon^*:\mathcal H_0\to\mathcal H_\partial$
the adjoint of $B_\varepsilon$ with respect to the inner products of
$\mathcal H_\partial$ and $\mathcal H_0$. Consequently,
$T_\varepsilon=B_\varepsilon B_\varepsilon^*$ is a compact, positive
self-adjoint operator on $\mathcal H_0$.
Moreover, using \eqref{e2.green-representation} and
\eqref{e7.rescaled-density}, we obtain
\begin{equation*}
(T_\varepsilon f)_i(y)
=\beta g_{i,\varepsilon}(y)^{\frac{1}{2}}\sum_{j=1}^m\int_{\mathbb R^3}
G_{\partial\Omega}\bigl(\exp_{\xi_i}(\delta_i y),
\exp_{\xi_j}(\delta_j z)\bigr)
 g_{j,\varepsilon}(z)^{\frac{1}{2}}q(z)f_j(z)\,dz.
\end{equation*}
The product involving $G_{\partial\Omega}$ is extended by zero
when $\delta_i|y|\ge2r_0$ or $\delta_j|z|\ge2r_0$.
To state the decomposition of $T_\varepsilon$, let
$e_i\in\mathcal H_0$, $1\le i\le m$, equal $1$ in the
$i$th component and $0$ in the others, and write $
(e_i\otimes e_j)f=\langle f,e_j\rangle_{\mathcal H_0}e_i.$
The operators $\mathscr K$ and $\Pi$ are defined in \eqref{e2.model-operators}.

\begin{lemma}\label{t7.4}
For $\xi$ in a fixed compact subset of $\mathcal M_m$ and
$\delta_i\asymp\delta_*\to0$ for all $i$, we have
\begin{equation}\label{e7.operator-expansion}
T_\varepsilon
=\sum_{i,j=1}^mA_{ij,\varepsilon}(e_i\otimes e_j)
+\bigoplus_{i=1}^m\mathscr K+E_\varepsilon,
\quad
A_{ij,\varepsilon}
=\begin{cases}
6|\log\delta_i|+\beta r(\xi_i),&i=j,\\
\beta G_{\partial\Omega}(\xi_i,\xi_j),&i\ne j.
\end{cases}
\end{equation}
where
\begin{equation}\label{e7.operator-remainder}
\|E_\varepsilon\|_{\mathcal L(\mathcal H_0)}
\le C\delta_*^{\frac{1}{4}}(1+|\log\delta_*|)^2.
\end{equation}
\end{lemma}

\begin{proof}
Set $R=\delta_*^{-1/4}$. By \eqref{e3.2} and
\eqref{e5.jacobian-derivatives},
$g_{i,\varepsilon}^{1/2}=1+O(\delta_*^2R^2)$ for $|y|\le R$.
For $|y|,|z|\le R$, \eqref{e2.6} and \eqref{e2.10} give
\begin{equation}\label{e7.diagonal-kernel-expansion}
\beta G_{\partial\Omega}\bigl(\exp_{\xi_i}(\delta_i y),
\exp_{\xi_i}(\delta_i z)\bigr)
=A_{ii,\varepsilon}-6\log|y-z|+O(\delta_*R)+O\!\left(\delta_*^2R^2
\bigl(1+|\log\delta_*|+|\log|y-z||\bigr)\right).
\end{equation}
For $i\ne j$, separation of the centers implies
\begin{equation}\label{e7.offdiagonal-kernel-expansion}
\beta G_{\partial\Omega}\bigl(\exp_{\xi_i}(\delta_i y),
\exp_{\xi_j}(\delta_j z)\bigr)
=A_{ij,\varepsilon}+O(\delta_*R).
\end{equation}
Since $q(y)\le C(1+|y|)^{-6}$, we have
\begin{equation}\label{e7.log-square-bound}
\int_{\mathbb R^3}q(z)|\log|y-z||^2\,dz
\le C\bigl(1+\log(2+|y|)\bigr)^2
\end{equation}
and
\begin{equation}\label{e7.log-tail}
\int_{|y|>R}q(y)\bigl(1+\log(2+|y|)\bigr)^2\,dy
\le CR^{-3}(1+\log R)^2.
\end{equation}
By \eqref{e7.operator-expansion}, we arrive at
\begin{equation}\label{e7.remainder-integral}
\begin{aligned}
(E_\varepsilon f)_i(y)
&=\sum_{j=1}^m\int_{\mathbb R^3}(E_\varepsilon)_{ij}(y,z)q(z)f_j(z)\,dz,\\
(E_\varepsilon)_{ij}(y,z)
&=\beta g_{i,\varepsilon}(y)^{\frac{1}{2}}g_{j,\varepsilon}(z)^{\frac{1}{2}}
G_{\partial\Omega}\bigl(\exp_{\xi_i}(\delta_i y),\exp_{\xi_j}(\delta_j z)\bigr)-A_{ij,\varepsilon}+6\delta_{ij}\log|y-z|.
\end{aligned}
\end{equation}
For $|y|,|z|\le R$, we know
$g_{i,\varepsilon}(y)^{1/2}g_{j,\varepsilon}(z)^{1/2}
=1+O(\delta_*^2R^2)$. Combining
\eqref{e7.diagonal-kernel-expansion}--\eqref{e7.offdiagonal-kernel-expansion}
and \eqref{e7.remainder-integral},
\[
|(E_\varepsilon)_{ij}(y,z)|
\le C\delta_*R+C\delta_*^2R^2
\bigl(1+|\log\delta_*|+\bigl|\log|y-z|\bigr|\bigr).
\]
For $\delta_i|y|,\delta_j|z|<2r_0$,
\begin{equation*}
d\bigl(\exp_{\xi_i}(\delta_i y),\exp_{\xi_j}(\delta_j z)\bigr)
\asymp\delta_i|y-z|\quad (i=j),\quad
d\bigl(\exp_{\xi_i}(\delta_i y),\exp_{\xi_j}(\delta_j z)\bigr)
\ge C^{-1} \quad (i\ne j).
\end{equation*}
Since $0\le g_{i,\varepsilon}\le C$, \eqref{e2.green-log-bound}
and \eqref{e7.remainder-integral} imply
\[
|(E_\varepsilon)_{ij}(y,z)|
\le C\bigl(1+|\log\delta_*|+\bigl|\log|y-z|\bigr|\bigr),
\qquad \max\{|y|,|z|\}>R.
\]
By Cauchy--Schwarz and \eqref{e7.log-square-bound}--\eqref{e7.log-tail},
\[
\begin{aligned}
\|E_\varepsilon\|
&\le\left(\sum_{i,j=1}^m\iint_{\mathbb R^3\times\mathbb R^3}
|(E_\varepsilon)_{ij}(y,z)|^2q(y)q(z)\,dy\,dz\right)^{\frac{1}{2}}\\
&\le C\left[\delta_*R+
\bigl(\delta_*^2R^2+R^{-\frac{3}{2}}\bigr)(1+|\log\delta_*|)^2\right]
\le C\delta_*^{\frac{1}{4}}(1+|\log\delta_*|)^2.
\end{aligned}
\]
\end{proof}

Define $\mathbb J_0,\mathbb J_\varepsilon:\mathbb R^{4m}\to\mathcal H_0$
by their columns
\[
(\mathbb J_0)_{i,k}=2z_ke_i,
\qquad
(\mathbb J_\varepsilon)_{i,k}
=2g_{i,\varepsilon}^{\frac{1}{2}}\chi(\delta_i|y|)z_ke_i,
\qquad 1\le i\le m,\quad 0\le k\le3.
\]
Combining \eqref{e4.moment-space} and \eqref{e7.rescaled-density}, we obtain
\[
\mathbb J_\varepsilon^*B_\varepsilon u
=2\sqrt\beta\,\bigl(\ell_{i,k}(u)\bigr)_{i,k},
\qquad
X_\perp=\ker(\mathbb J_\varepsilon^*B_\varepsilon).
\]
Moreover, \eqref{e2.17} implies
$\mathbb J_0^*\mathbb J_0=I_{\mathbb R^{4m}}$.
The estimates for $g_{i,\varepsilon}^{1/2}$ in the proof of
Lemma~\ref{t7.4} and \eqref{e7.log-tail} yield
\begin{equation}\label{e7.4}
\|\mathbb J_\varepsilon-\mathbb J_0\|\le C\delta_*^{\frac{1}{4}}.
\end{equation}

\begin{proposition}\label{t7.5}
Define $
\mathcal B_0(u,v)
=\langle\mathcal Pu,v\rangle
-\beta\sum_{i=1}^m\int_{\partial\Omega}q_iuv\,d\sigma
=\langle u,v\rangle_{\mathcal H_\partial}
-\langle B_\varepsilon u,B_\varepsilon v\rangle_{\mathcal H_0}.$
There exists $c\in(0,1)$ such that
for  sufficiently small $\varepsilon$,
\[
\ind(\mathcal B_0|_{X_\perp})=m,
\qquad
\spec(\mathcal B_0|_{X_\perp})
\subset(-\infty,-c\Lambda_\varepsilon]\cup[c,1].
\]
Here $\spec(\mathcal B_0|_{X_\perp})$ denotes the spectrum of the
corresponding operator with respect to the $\mathcal H_\partial$ inner product.
The index identity and spectral inclusion also hold for
$\mathcal B_V|_{X_\perp}$ and $\mathcal B_\varepsilon|_{X_\perp}$.
\end{proposition}
\begin{proof}
Set $T^0_\varepsilon
:=\sum_{i,j=1}^m A_{ij,\varepsilon}(e_i\otimes e_j)
+\bigoplus_{i=1}^m\mathscr K.$
By \eqref{e2.16}--\eqref{e2.17}
and \eqref{e4.model-jacobi-potential}, we have
\[
\mathbb J_0^*T^0_\varepsilon\mathbb J_0=I,
\qquad \|T^0_\varepsilon\mathbb J_0\|\le C,
\qquad \|T^0_\varepsilon\|\le C\Lambda_\varepsilon.
\]
It follows from \eqref{e7.operator-remainder} and \eqref{e7.4} that
\begin{equation}\label{e7.gram-error}
\|\mathbb J_\varepsilon^*T_\varepsilon\mathbb J_\varepsilon-I\|
\le C\delta_*^{\frac{1}{4}}\Lambda_\varepsilon^2=o(1).
\end{equation}
Hence the $\mathcal H_\partial$-orthogonal projection onto $X_\perp$ is
\[
P_X=I-B_\varepsilon^*\mathbb J_\varepsilon
(\mathbb J_\varepsilon^*T_\varepsilon\mathbb J_\varepsilon)^{-1}
\mathbb J_\varepsilon^*B_\varepsilon.
\]
The compact positive operators $P_XB_\varepsilon^*B_\varepsilon|_{X_\perp}$
and $T_{X,\varepsilon}=B_\varepsilon P_XB_\varepsilon^*
=T_\varepsilon-T_\varepsilon\mathbb J_\varepsilon
(\mathbb J_\varepsilon^*T_\varepsilon\mathbb J_\varepsilon)^{-1}
\mathbb J_\varepsilon^*T_\varepsilon$
have the same nonzero eigenvalues, with the same multiplicities.
Set $
T^0_{X,\varepsilon}
=T^0_\varepsilon-T^0_\varepsilon\mathbb J_0\mathbb J_0^*T^0_\varepsilon.$
Using \eqref{e7.operator-remainder}, \eqref{e7.4}, and
\eqref{e7.gram-error},  on the fixed space $\mathcal H_0$ we obtain
\begin{equation}\label{e7.5}
\|T_{X,\varepsilon}-T^0_{X,\varepsilon}\|
\le C\delta_*^{\frac{1}{4}}\Lambda_\varepsilon^4=o(1).
\end{equation}
For $f\perp\operatorname{span}\{e_1,\ldots,e_m\}$,
\eqref{e4.model-jacobi-potential} gives
$\mathbb J_0^*T^0_\varepsilon f=\mathbb J_0^*f$. Therefore,
\eqref{e2.25} implies
\[
\langle T^0_{X,\varepsilon}f,f\rangle_{\mathcal H_0}
=\sum_{i=1}^m\langle\mathscr Kf_i,f_i\rangle_{L^2(q\,dy)}
-|\mathbb J_0^*f|^2
\le\tfrac14\|f\|_{\mathcal H_0}^2.
\]
On the other hand, from \eqref{e7.operator-expansion} and
$\|T^0_\varepsilon\mathbb J_0\|\le C$, we obtain
\[
\langle T^0_{X,\varepsilon}f,f\rangle_{\mathcal H_0}
\ge\bigl(6\min_i|\log\delta_i|-C\bigr)\|f\|_{\mathcal H_0}^2,
\qquad f\in\operatorname{span}\{e_1,\ldots,e_m\}.
\]
Let $\lambda_j$ denote eigenvalues in nonincreasing order, counted
with multiplicity. Applying the min--max principle~\cite{r30} to the
quadratic-form estimates for $T^0_{X,\varepsilon}$ on
$\operatorname{span}\{e_1,\ldots,e_m\}$ and its orthogonal
complement, and using \eqref{e7.5}, we obtain
\begin{equation}\label{boundt}
\lambda_m(T_{X,\varepsilon})\ge6\min_i|\log\delta_i|-C,
\qquad
\lambda_{m+1}(T_{X,\varepsilon})\le\tfrac14+o(1).
\end{equation}
For $u,v\in X_\perp$, $
\mathcal B_0(u,v)
=\langle u-P_XB_\varepsilon^*B_\varepsilon u,v\rangle_{\mathcal H_\partial}.$
Consequently, \eqref{boundt} and
$|\log\delta_i|\asymp\Lambda_\varepsilon$ imply
\[
\ind(\mathcal B_0|_{X_\perp})=m,
\qquad
\spec(\mathcal B_0|_{X_\perp})
\subset(-\infty,-c\Lambda_\varepsilon]\cup[\frac{3}{4}-o(1),1].
\]

From \eqref{e3.14}, \eqref{e5.3} and \eqref{e4.potential-bound}, we get $
\left|\varepsilon he^{\beta V}-\sum_{i=1}^m q_i\right|
+\left|\rho_\varepsilon-\sum_{i=1}^m q_i\right|
\le C\delta_*\Lambda_\varepsilon^2\omega.$
Thus \eqref{e4.weighted-embedding} implies
\begin{equation}\label{e7.potential-perturbation}
|\mathcal B_V(u,u)-\mathcal B_0(u,u)|
+|\mathcal B_\varepsilon(u,u)-\mathcal B_0(u,u)|\le C\delta_*\Lambda_\varepsilon^3
\|u\|_{\mathcal H_\partial}^2
=o(1)\|u\|_{\mathcal H_\partial}^2.
\end{equation}
Since $Q_V,\rho_\varepsilon>0$, we obtain
\[
\max\{\mathcal B_V(u,u),\mathcal B_\varepsilon(u,u)\}
\le\|u\|_{\mathcal H_\partial}^2,
\qquad u\in X_\perp.
\]
The min--max principle and \eqref{e7.potential-perturbation} therefore yield
\begin{equation*}
\ind(\mathcal B_V|_{X_\perp})
=\ind(\mathcal B_\varepsilon|_{X_\perp})=m,\quad
\spec(\mathcal B_V|_{X_\perp})
\cup\spec(\mathcal B_\varepsilon|_{X_\perp})
\subset(-\infty,-c\Lambda_\varepsilon]\cup[\tfrac34-o(1),1].
\end{equation*}
\end{proof}

\begin{corollary}\label{t7.6}
For sufficiently small $\varepsilon$,
$\mathcal B_V|_{X_{\mathrm{tr},\varepsilon}}$ and
$\mathcal B_\varepsilon|_{X_{\mathrm{tr},\varepsilon}}$
are nondegenerate and have index $m$.
\end{corollary}
\begin{proof}
By \eqref{e4.dilation-orthogonality} and \eqref{e4.dilation-approximation},
\[
X_{\mathrm{tr},\varepsilon}
=X_\perp\oplus\operatorname{span}\{\widehat Z_{i,0}:1\le i\le m\}
\]
is a $\mathcal B_V$-orthogonal decomposition. Since
$\mathcal B_V|_{\operatorname{span}\{\widehat Z_{i,0}:1\le i\le m\}}>0$
by \eqref{e4.9}, this orthogonal decomposition and
Proposition~\ref{t7.5} yield
\begin{equation*}
\ker(\mathcal B_V|_{X_{\mathrm{tr},\varepsilon}})
=\ker(\mathcal B_V|_{X_\perp})=\{0\},\quad
\ind(\mathcal B_V|_{X_{\mathrm{tr},\varepsilon}})
=\ind(\mathcal B_V|_{X_\perp})=m.
\end{equation*}
Fix $t\in[0,1]$ and
$u\in\ker\bigl(((1-t)\mathcal B_V+t\mathcal B_\varepsilon)
|_{X_{\mathrm{tr},\varepsilon}}\bigr)$.
By the definition of $X_{\mathrm{tr},\varepsilon}$ (see \eqref{matr}),
\[
L_Vu-t\beta\varepsilon h
\bigl(e^{\beta v_\varepsilon}-e^{\beta V}\bigr)u
\in\operatorname{span}\{Q_{i,k}:1\le i\le m,\ 1\le k\le3\}.
\]
For fixed $\varepsilon$, $\mathcal Pu\in L^2(\partial\Omega)$, so the
regularity argument in the proof of Theorem~\ref{t1.3} implies
$u\in L^\infty(\partial\Omega)$.
Applying \eqref{e4.10} to $u$ and using
\eqref{e4.potential-bound} and \eqref{e5.3},  uniformly in $t\in[0,1]$, we obtain
\begin{equation*}
\|u\|_\infty
\le C\Lambda_\varepsilon
\bigl\|t\beta\varepsilon h(e^{\beta v_\varepsilon}-e^{\beta V})u\bigr\|_*
\le C\Lambda_\varepsilon\|\phi_{\varepsilon,\xi_\varepsilon}\|_\infty
\|u\|_\infty
\le C\delta_*\Lambda_\varepsilon^3\|u\|_\infty.
\end{equation*}
Since $\delta_*\Lambda_\varepsilon^3\to0$, we conclude
\[
\ker\bigl(((1-t)\mathcal B_V+t\mathcal B_\varepsilon)
|_{X_{\mathrm{tr},\varepsilon}}\bigr)=\{0\},
\qquad t\in[0,1].
\]
For $t\in[0,1]$, define
$\mathscr A_t:X_{\mathrm{tr},\varepsilon}\to X_{\mathrm{tr},\varepsilon}$ by
\[
\langle\mathscr A_tu,v\rangle_{\mathcal H_\partial}
=(1-t)\mathcal B_V(u,v)+t\mathcal B_\varepsilon(u,v),
\qquad u,v\in X_{\mathrm{tr},\varepsilon}.
\]
Then $
\mathscr A_t^*=\mathscr A_t,$
$\mathscr A_t-I_{X_{\mathrm{tr},\varepsilon}}\ \text{is compact},$
$\|\mathscr A_t-\mathscr A_s\|
=|t-s|\,$ $\|\mathscr A_1-\mathscr A_0\|,$
$s,t\in[0,1].$
Consequently,
\[
\ind\bigl(((1-t)\mathcal B_V+t\mathcal B_\varepsilon)
|_{X_{\mathrm{tr},\varepsilon}}\bigr)
=\ind(\mathcal B_V|_{X_{\mathrm{tr},\varepsilon}})=m,
\qquad t\in[0,1].
\]
Setting $t=1$, we deduce that $
\ker(\mathcal B_\varepsilon|_{X_{\mathrm{tr},\varepsilon}})=\{0\},$
$\ind(\mathcal B_\varepsilon|_{X_{\mathrm{tr},\varepsilon}})=m.$
\end{proof}

\begin{proof}[Proof of Theorem~\ref{t1.4}]
By \eqref{e5.energy-expansion} and $\xi_\varepsilon\to\xi^*$, we have $
D^2\mathcal K_\varepsilon(\xi_\varepsilon)
\rightarrow-D^2\mathcal F_m(\xi^*).$
Since $\det D^2\mathcal F_m(\xi^*)\ne0$, for sufficiently small $\varepsilon$,
\[
\ind\bigl(D^2\mathcal K_\varepsilon(\xi_\varepsilon)\bigr)
=\ind\bigl(-D^2\mathcal F_m(\xi^*)\bigr).
\]
Thus, applying \eqref{e7.translation-splitting}--\eqref{e7.translation-hessian},
Lemma~\ref{t7.3}, and Corollary~\ref{t7.6}, we obtain
\begin{equation*}
\ind(u_\varepsilon)
=\ind(\mathcal B_\varepsilon)
=m+\ind\bigl(D^2\mathcal K_\varepsilon(\xi_\varepsilon)\bigr)=m+\ind\bigl(-D^2\mathcal F_m(\xi^*)\bigr)
=4m-\ind\bigl(D^2\mathcal F_m(\xi^*)\bigr).
\end{equation*}
\end{proof}

\subsection{Generic nondegeneracy of critical points}\label{s8}

For $\mathfrak h\in C^5(\partial\Omega)$, let
$\mathcal F_{m,\mathfrak h}$ denote $\mathcal F_m$ with $h=e^{\mathfrak h}$.
By Proposition~\ref{t2.3}, $\mathcal F_{m,\mathfrak h}\in C^4(\mathcal M_m)$.

\begin{lemma}\label{t8.1}
Suppose $m\ge1$ and  $K\subset\mathcal M_m$ is compact.
Then \[
\mathcal O_{m,K}
:=
\left\{
\mathfrak h\in C^5(\partial\Omega):
\ \forall\,\xi\in K,\quad
d\mathcal F_{m,\mathfrak h}(\xi)=0
\ \Longrightarrow\
\det D^2\mathcal F_{m,\mathfrak h}(\xi)\ne0
\right\}
\]
 is open and dense in \(C^5(\partial\Omega)\).
\end{lemma}
\begin{proof}
By compactness of \(K\) and the continuous dependence of \(\mathcal F_{m,\mathfrak h}\) on \(\mathfrak h\) in \(C^2\), the set \(\mathcal O_{m,K}\) is open in \(C^5(\partial\Omega)\).
For density, fix $\mathfrak h_0\in C^5(\partial\Omega)$.
Using local coordinates and compactness of \(K\), choose \(\psi_1,\ldots,\psi_{N_0}\in C^\infty(\partial\Omega)\) and an open set \(K\subset\mathcal U_K\subset\subset\mathcal M_m\) such that
\[
 \operatorname{span}\left\{
 \bigl(d\psi_\ell(\xi_i)\bigr)_{i=1}^m:1\le\ell\le N_0
 \right\}=T_\xi^*\mathcal M_m,
 \qquad \xi\in\mathcal U_K.
\]
Set $\mathfrak h_t=\mathfrak h_0+\sum_{\ell=1}^{N_0}t_\ell\psi_\ell,$  $t\in\mathbb R^{N_0}.$
Since
\begin{equation*}
 \partial_{t_\ell}\bigl(d\mathcal F_{m,\mathfrak h_t}(\xi)\bigr)
 =\frac1{24\pi^2}\bigl(d\psi_\ell(\xi_i)\bigr)_{i=1}^m,\quad
 \operatorname{rank}D_t\bigl(d\mathcal F_{m,\mathfrak h_t}(\xi)\bigr)=3m,
\end{equation*}
the implicit function theorem yields a $C^3$ manifold
\[
 \mathcal Z=\bigl\{(\xi,t)\in\mathcal U_K\times\mathbb R^{N_0}:
 d\mathcal F_{m,\mathfrak h_t}(\xi)=0\bigr\},
 \qquad \dim\mathcal Z=N_0.
\]
For $\operatorname{pr}_t:\mathcal Z\to\mathbb R^{N_0}$,
$(\xi,t)\mapsto t$, we have $
 \ker D\operatorname{pr}_t(\xi,t)
 =\ker D^2\mathcal F_{m,\mathfrak h_t}(\xi)\times\{0\}.$
Thus, by Sard's theorem~\cite{r23}, there exist regular values
$t_n\to0$ of $\operatorname{pr}_t$ such that
\[
 \mathfrak h_{t_n}\in\mathcal O_{m,K},
 \qquad
 \|\mathfrak h_{t_n}-\mathfrak h_0\|_{C^5(\partial\Omega)}\rightarrow0,
\]
which proves density.
\end{proof}

\begin{proof}[Proof of Corollary~\ref{t1.5}]
For $m\ge2$, consider
\begin{equation*}
 K_{m,n}=\bigl\{\xi\in(\partial\Omega)^m:
 \min_{i\ne j}|\xi_i-\xi_j|\ge n^{-1}\bigr\},\quad n\ge1,\quad
 \mathcal M_m=\bigcup_{n\ge1}K_{m,n}.
\end{equation*}
By Lemma~\ref{t8.1}, Baire's theorem, and the homeomorphism
$C^5(\partial\Omega)\ni\mathfrak h\mapsto e^{\mathfrak h}
\in\mathcal W_+^5(\partial\Omega)$, we obtain
\[
 \mathcal W_{\mathrm{reg}}
 =\left\{e^{\mathfrak h}:\mathfrak h\in
 \mathcal O_{1,\partial\Omega}\cap
 \bigcap_{m\ge2}\bigcap_{n\ge1}\mathcal O_{m,K_{m,n}}\right\}
\]
is residual and dense in $\mathcal W_+^5(\partial\Omega)$.
For $h\in\mathcal W_{\mathrm{reg}}$ and $m\ge1$,
\[
 d\mathcal F_m(\xi)=0
 \quad\Longrightarrow\quad
 \det D^2\mathcal F_m(\xi)\ne0,
 \qquad \xi\in\mathcal M_m.
\]
In particular,
\[
 \xi^{(m)}\in\operatorname*{arg\,min}_{\mathcal M_m}\mathcal F_m
 \quad\Longrightarrow\quad D^2\mathcal F_m(\xi^{(m)})>0.
\]
Thus, Theorems~\ref{t1.1} and~\ref{t1.4} provide $m$-bubble solutions
associated with $\xi^{(m)}$ satisfying
\[
 \ind(u_\varepsilon)
 =m+\ind\bigl(-D^2\mathcal F_m(\xi^{(m)})\bigr)=4m.
\]
\end{proof}

\bigskip
\noindent\textbf{Acknowledgments}
The research was supported by the National Natural Science Foundation of China, grants 12371121 and 12671143, and by the Fundamental Research Funds for the Central Universities, grant SWU-KF26002.

\vspace{.1cm}
\noindent\textbf{Data availability.}
No data were used for the research described in this article.

\vspace{.1cm}
\noindent\textbf{Conflict of interest.}
The authors declare that they have no conflict of interest.

\end{document}